\documentclass{HMS}

\begin{document} \title{Arithmetic mirror symmetry for the 2-torus}
\begin{abstract} 
	This paper explores a refinement of homological mirror symmetry which relates exact symplectic topology 
	to arithmetic algebraic geometry. We establish a derived equivalence of the Fukaya category of
	the 2-torus, relative to a basepoint, with the category of perfect
	complexes of coherent sheaves on the Tate curve over the formal disc $\spec \Z\series{q}$. It
	specializes to a derived equivalence, over $\Z$, of the Fukaya category of the
	punctured torus with perfect complexes on the curve $y^2+xy=x^3$ over $\spec \Z$, 
	the central fibre of the Tate curve; and, over the `punctured disc' $\spec \Z\laurent{q}$, 
	to an integral refinement of the known statement of homological mirror
	symmetry for the 2-torus.  We also prove that the wrapped Fukaya category 
	of the punctured torus is  derived-equivalent over $\Z$ to coherent sheaves 
	on the central fiber of the Tate curve.

\end{abstract}
\author{Yank{\i} Lekili and Timothy Perutz}
\address{Author affiliations}
\address{YL: University of Cambridge}
\email{yl319@cam.ac.uk}
\address{TP: University of Texas at Austin}
\email{perutz@math.utexas.edu}
\maketitle

%\begin{center}
%\today\\
%\currenttime
%\end{center}

\section{Introduction}
This paper explores a basic case of what we believe is a general connection between exact Lagrangian submanifolds in the complement to an ample divisor $D$ in a complex Calabi--Yau manifold $X$---we view $X\setminus D$ as an exact symplectic manifold---and coherent sheaves on a scheme \emph{defined over $\spec \Z$}, the `mirror' to $X\setminus D$. We take $X$ to be an elliptic curve; its complex structure is irrelevant, so it is really a 2-torus $T$. We take $D$ to be a point $z$. The mirror is the Weierstrass cubic $Y^2Z + XYZ = X^3$, the restriction to $q=0$ of the Tate curve $\EuT\to \spec \Z\series{q}$.

Kontsevich's 1994 homological mirror symmetry (HMS) conjecture \cite{Kon} claims that the Fukaya $A_\infty$-category $\EuF(X)$ of a polarized Calabi--Yau manifold should  have a formal enlargement---precisely formulated a little later as the closure $\tw^\pi \EuF(X)$ under taking mapping cones and passing to idempotent summands---which is $A_\infty$-quasi-equivalent to a dg enhancement for the derived category of coherent sheaves on the `mirror' $\check{X}$, a Calabi--Yau variety over the field of complex Novikov series.\footnote{Beware: the circumstances under which one expects to find such an $\check{X}$ are more subtle than those claimed by our one-sentence pr\'ecis of Kontsevich's conjecture.} The HMS conjecture has inspired a great deal of work in symplectic geometry, algebraic geometry and mathematical physics; the HMS paradigm has been adapted so as to apply not only to varieties whose canonical bundle $\EuK$ is trivial, but also to those where either $\EuK^{-1}$ or $\EuK$ is ample, with such varieties playing either symplectic or algebro-geometric roles. Meanwhile, progress on the original case of Calabi--Yau manifolds has been slow. There are currently complete mirror-symmetric descriptions of the Fukaya category only for the 2-torus $\R^2/\Z^2$ and of the square 4-torus $\R^4/\Z^4$ \cite{AS}. The case of Calabi--Yau hypersurfaces in projective space has been solved up to a certain ambiguity in identifying the mirror variety \cite{SeiQuartic, She}. There are significant partial results for linear symplectic tori of arbitrary dimension \cite{KS}. 

Our contention is that even in the solved cases, there is more to be said about HMS. The Fukaya category for the 2-torus has a natural model which is defined over $\Z\series{q}$, a subring of the complex Novikov field. This model has a mirror-symmetric description as the perfect complexes on the Tate curve $\EuT$ over $\Z\series{q}$. The symplectic geometry of the torus is thereby connected to the arithmetic algebraic geometry of $\EuT$. Establishing this connection is the task of this article. 

Experts have certainly been aware that, in principle, homological mirror symmetry should have an arithmetic-geometric dimension (cf. Kontsevich's lecture \cite{KonLec}, for instance), but we believe that this article is the first to treat this idea in detail. Whilst refining existing proofs of HMS for the 2-torus might   be a viable option, our method is also new: we identify a generating subalgebra $\EuA$ of the Fukaya category, and show that Weierstrass cubic curves precisely parametrize the possible $A_\infty$-structures on it (Theorem \ref{Ainf comparison}). The mirror to $(T,z)$ is then the unique Weierstrass curve corresponding to the $A_\infty$-structure belonging to the Fukaya category. Our identification of this mirror parallels an argument of Gross \cite{Gro} but also has a novel aspect, relating the multiplication rules for theta-functions on the Tate curve to \emph{counts of lattice points} in triangles (not areas of triangles). Our identification of the \emph{wrapped} Fukaya category of the punctured torus with coherent complexes on $\EuT|_{q=0}$ appears to be a basic case of an unexplored aspect of mirror symmetry for Calabi--Yau manifolds.

\subsection{Statement}
Let $T$ be a closed, orientable surface of genus 1; $\omega$ a symplectic form on $T$; $z\in T$ a basepoint; $T_0=T\setminus \{z\}$; and $\theta$ a primitive for $\omega$ on $T_0$. Fix also a grading for the symplectic manifold $T$, that is, an unoriented line-field $\ell$. These data suffice to specify the relative Fukaya category $\EuF(T,z)$ up to quasi-isomorphism. It is an $A_\infty$-category linear over $\Z\series{q}$ whose objects are embedded circles $\gamma\subset T_0$ which are exact ($\int_\gamma \theta =0$) and are equipped with orientations, double covers $\tilde{\gamma}\to \gamma$ and gradings (a grading is a homotopy from $\ell|_\gamma$ to $T\gamma$ in $T(T_0)|_\gamma$).  

Let  $\EuT\to \spec \Z\series{q}$ denote the Tate curve, the cubic curve in $\mathbb{P}^2(\Z\series{q})$ with equation
\begin{equation}\label{Tate curve} 
Y^2Z + XY Z= X^3 +a_4(q) XZ^2 + a_6(q) Z^3,  
\end{equation}
where
\begin{equation}  \label{Tate coeffs}
a_4(q) =  -5\sum_{n>0}{\frac{n^3q^n }{1-q^n}},\quad
a_6(q) = - \frac{1}{12} \sum_{n>0}{ \frac{(5 n^3 + 7 n^5)q^n } {1-q^n}}
\end{equation}
(note that $n^2(5 + 7n^2)$ is always divisible by 12).

Let $\VB(\EuT)$ denote the $\Z\series{q}$-linear differential graded (dg) category whose objects are locally free sheaves of finite rank over $\EuT$, and whose morphism spaces are \v{C}ech complexes with respect to a fixed affine open cover: $\hom_{\VB(\EuT)}(\EuE,\EuF) = \check{C}^\bullet(\SheafHom(\EuE,\EuF))$. 

\begin{figure}[ht!]
\centering
\labellist
\small\hair 2pt 
\pinlabel  {$L_0^\#$ (slope $\beta$)} at 195 10
\pinlabel  {$L_\infty^\# $ (slope $\alpha$)} at -20 180
\pinlabel  $L_{(1,-5)}^\#$  at 185 180
\pinlabel  $z$ at 60 60
\pinlabel  {$L_\infty^\# \longleftrightarrow \EuO_\sigma$} at 420 240
\pinlabel  {$L_0^\# \longleftrightarrow \EuO$} at 420 200
\pinlabel  {$L_{(1,-n)}^\# \longleftrightarrow \EuO(np)$} at 420 160
\pinlabel  {\tiny{horizontal line field}} at 300 240
\pinlabel {\tiny{grades $T$}} at 300 230
\pinlabel {\tiny{rotate line field to grade a Lagrangian}} at 200 295
\pinlabel {\tiny{stars indicate non-triviality}} at -30 300
\pinlabel {\tiny{of the double cover}} at -30 290
\endlabellist
\includegraphics[width=0.6\textwidth]{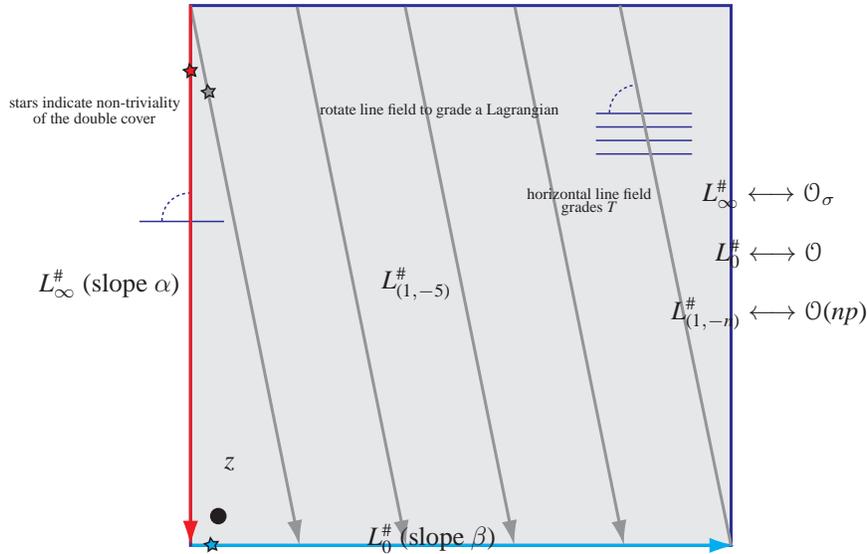}
\caption{The torus $T$ and the mirror correspondence $\psi$, for one possible choice of the line field $\ell$.}\label{HMScorrespondence}
\end{figure}

\begin{main}\label{mainth}
A choice of basis $(\alpha,\beta)$ for $H_1(T)$, with $\alpha\cdot \beta =1$, determines, canonically up to an overall shift and up to natural quasi-equivalence, a $\Z\series{q}$-linear $A_\infty$-functor
\[  \psi \colon \EuF(T,z) \to  \tw(\VB(\EuT))  \]
from the relative Fukaya category to the dg category of twisted complexes in $\VB(\EuT)$. Moreover,
\begin{enumerate}
\item [(i)]
the functor $\psi$ maps an object $L_0^\#$ representing $\beta$ to the structure sheaf $\EuO$. It maps an object $L_\infty^\#$ representing $\alpha$ to the complex $[ \EuO\to\EuO(\sigma)]$, where $\sigma=[0:1:0]$ is the section at infinity of $\EuT$, and the map is the inclusion. (This complex is quasi-isomorphic to the skyscraper sheaf $\EuO_\sigma = \sigma_*\EuO_{\spec \Z [[q]]}$ at the section at infinity.) It is an embedding on the full subcategory $\EuA$ on  $\{ L_0^\#, L_\infty^\#\}$; and is characterized, up to natural equivalence, by its restriction to $\EuA$. See Figure \ref{HMScorrespondence}.
\item [(ii)]
$\psi$ extends to an equivalence
\[ \dersplit \EuF(T,z) \to  \perf(\EuT)\simeq H^0(\tw \VB(\EuT))  \]
from the idempotent-closed derived Fukaya category to the triangulated category of perfect complexes on $\EuT$.
\item [(iii)]
The specialization of $\psi$ to $q=0$ is a $\Z$-linear functor
\[ \psi_0\colon \EuF(T_0)\to \tw \VB(\EuT|_{q=0}) \] 
from the exact Fukaya category of $(T_0,\theta)$ to the category of perfect complexes on the central fiber of the Tate curve, 
inducing an equivalence on derived categories 
\[ \der \psi_0\colon \der \EuF(T_0)\to \perf(\EuT|_{q=0}) \] 
(both of these derived categories are already idempotent-closed).
\item [(iv)]
$\der \psi_0$ extends to an equivalence of triangulated categories
\[ \der \EuW(T_0) \to \der^b \coh(\EuT|_{q=0}) \]
from the derived wrapped Fukaya category to the bounded derived category of coherent sheaves on $\EuT|_{q=0}$ (these derived categories are again idempotent-closed).
\end{enumerate}
\end{main}
\begin{rmk}
The functor $\psi$ has an additional property, which is that it is `trace-preserving', in a sense to be discussed later. 
\end{rmk}
Clause (ii) has the following corollary:
\begin{cor}
There is an $A_\infty$ quasi-equivalence $\Modules \EuF(T,z) \to \QC(\EuT)$ from the category of cohomologically unital $\EuF(T,z)$-modules to a DG enhancement of the derived category of unbounded quasi-coherent complexes on the Tate curve.
\end{cor}
Indeed, $\QC(\EuT)$ is quasi-equivalent to $\Modules \VB(\EuT)$ as an instance of the general theory of \cite{Toe} or \cite{BFN}.

\paragraph{Comparison to the standard formulation.}
The $A_\infty$-structure in the `relative' Fukaya category $\EuF(T,z)$ is based on counting holomorphic polygons weighted by powers $q^s$, where $s$ counts how many times the polygon passes through the basepoint $z$. The `absolute' Fukaya category $\EuF(T)$, in the version most popular for mirror symmetry, has as objects Lagrangian branes $L^\#$ in $T$ equipped with $U(1)$ local systems $E\to L$. In the latter version, holomorphic polygons are weighted by $(\mathrm{holonomy})\,q^{\mathrm{area}}$. The coefficient-ring for $\EuF(T)$ is usually taken to be $\Lambda_\C$, the field of complex Novikov series $\sum_{k>0}{a_k q^{r_k}}$: here $a_k\in \C$, $r_k \in \R$, and $r_k\to \infty$. 

To explain the relation between the relative and absolute versions, note first that there is an equation of currents $\omega=\delta_D + d\Theta$, where $\Theta$ is a 1-current. We take $\theta$ to be the (smooth) restriction of $\Theta$ to $M$.  
\begin{lem}
There is a fully faithful `inclusion' functor 
\[ e\colon \EuF(T,z)\otimes_{\Z[[q]]} \Lambda_{\C} \to \EuF(T), \]
linear over $\Lambda_{\C}$ and acting as the identity on objects.  For each exact Lagrangian $L$, select a function $K_L\in C^\infty(L)$ such that $dK_L = \theta|_L$.
Then define $e$ on morphism-spaces $\hom(L_0^\#,L_1^\#)= CF(\phi(L_0^\#),L_1^\#)$ by
\[ e(x) = q^{A(x)} x ,\quad x\in \phi(L_0) \cap L_1,   \] 
where $A(x)=A_{\phi(L_0),L_1}(x)$ is the symplectic action, defined via the $K_L$, and $\phi$ is the exact symplectomorphism used to obtain transversality. The higher $A_\infty$-terms for $e$ are identically zero.
\end{lem} 
\begin{pf}
The symplectic action is defined as follows. For a path $\gamma\colon ([0,1]; 0, 1) \to (M; L_0,L_1)$ (for instance, a constant path at an intersection point) we put 
\[A_{L_0,L_1}(\gamma)= - \int_0^1 {\gamma^*\theta} - K_{L_0}(\gamma(0)) + K_{L_1}(\gamma(1)).\]
For any disc $u\colon (D,\partial D)\to (X,L) $, we have
\[  \int_D{u^*\omega} - D\cdot u = \int_D{u^*(\omega-\delta_D)}= \int_{\partial D} {u|_{\partial D}^*\theta} 
= \int_{\partial D}{ d (u|_{\partial D} ^*K_L )} = 0.  \]
Similarly, if $u\colon D\to X$ is a polygon attached to a sequence of Lagrangians $(L_0,L_1,\dots, L_d)$ (where $d\geq 1$) at corners $x_1\in L_0\cap L_1,\dots, x_{d+1}\in L_d\cap L_0$, then
\[  \int_D{u^*\omega} - D\cdot u = \int_D{u^*(\omega-\delta_D)}= A_{L_{d+1}, L_0}(x_{d+1})+\sum_{i=1}^{d+1}{A_{L_{i-1},L_i}(x_i)}.  \]
From this it follows that $e \circ \mu^d_{\EuF(T,z)}(x_1,\dots,x_d) = \mu^d_{\EuF(T)} \circ (ex_1,\dots, ex_d)$, which proves that $e$ is a functor. Note that the perturbations that are used to define hom-spaces in $\EuF(T,z)$ serve equally well in $\EuF(T)$.
It is clear that $e$ is fully faithful.
\end{pf}

The `standard' statement of mirror symmetry is as follows. Let $\EuT_{\Lambda_\C}= \EuT \times_{\Z[[q]]} \Lambda_\C$; it is an elliptic curve over the field $\Lambda_\C$.
When $\omega$ is normalized so that $\int_T{\omega}=1$, there is a functor
\[ \Phi\colon \EuF(T)\to  \widetilde{\der}^b\coh(\EuT_{\Lambda_\C}),  \]  
where $\widetilde{\der}^b \coh$ is the unique dg enhancement of the bounded derived category $\der^b\coh$ \cite{LO}, inducing a derived equivalence; and that this functor is again canonically determined by a choice of basis for $H_1(T)$: see \cite{PZ, Pol00, Pol03, AS} for one proof; \cite{Gro} for an expository account of another, occasionally missing technical details (e.g. certain signs); and \cite{SeiFlux} for yet another.  Our result is an arithmetic refinement of this standard one:
\begin{thm}\label{comparing functors}
The diagram  
\[\xymatrix{
\EuF(T,z)\otimes \Lambda_\C  \ar[d]_{e} \ar[r]^{\psi\otimes 1}&  \tw \VB(\EuT) \otimes\Lambda_\C\ar^{i}[d] \\
\EuF(T) \ar[r]^{\Phi}  & \widetilde{\der}^b\coh(\EuT_{\Lambda_\C}).
} \]
is homotopy-commutative under composition of $A_\infty$-functors. 
\end{thm}
Since $\EuT\times_{\Z[[q]]} \Lambda_\C$ is a non-singular variety over the field $\Lambda_\C$, we may take $\tw \VB(\EuT_{\Lambda_\C})$ as our dg enhancement of $\der\coh(\EuT_{\Lambda_\C})$. Then $i$ is the base-change functor $\tw \VB(\EuT)\to \tw \VB(\EuT_{\Lambda_\C})$. For this theorem to make sense, $\psi$ and $\Phi$ must be set up so that $i\circ (\psi\otimes 1)$ and $\Phi\circ e$ agree precisely (not just up to quasi-isomorphism) on objects.

\subsection{The Tate curve}
Useful references for this material include \cite{Tat, Hid, Con, Gro}.
The \emph{Tate curve} is the plane projective curve $\EuT$ over  $\Z \series{q}$ whose affine equation is the Weierstrass cubic
\begin{equation} y^2 + xy = x^3 +a_4 x + a_6,  
\end{equation}
where $a_4$ and $a_6$ are as at (\ref{Tate coeffs}). So $\EuT$ is a projective curve in $\mathbb{P}^2(\Z\series{q})$. Like any Weierstrass curve $w(x,y)=0$, $\EuT$ comes with a canonical differential with poles at the singularities,
\[ \Omega=dx/w_y= -dy/w_x = dx/(2y+x) = -dy / (y-3x^2 -a_4) .\]

Notation: 
\begin{equation}
\hat{\EuT} = \text{$\EuT$ specialized to $\Z\laurent{q} (=\Z\series{q}[q^{-1}])$}.
\end{equation}

The analytic significance of the Tate curve is the following. Consider the Riemann surface $E_\tau=\C/\langle 1,\tau \rangle$, where $\imag \tau>0$. The exponential map $z\mapsto q := \exp(2\pi iz)$ identifies $E_\tau$ with $\C^*/  q^{\Z} $. As $q$ varies over the punctured unit disc $\mathbb{D}^*$, the Riemann surfaces $\C^*/  q^{\Z} $ form a holomorphic family $\mathbf{E} \to \mathbb{D}^*$. The Weierstrass function $\wp_q$ for the modular parameter $q$ defines an embedding 
\[  \mathbf{E} \to \mathbb{CP}^2 \times \mathbb{D}^*;\quad 
(z,q) \mapsto ([(2\pi i)^{-2}\wp_q(z): (2\pi i)^{-3}\wp'_q(z) : 1 ], q). \]
This embedding is cut out by an equation $y^2 = 4x^3 - g_2(q)x  - g_3(q)$, which is a Weierstrass cubic in $(x,y)$ varying holomorphically with $q$. The functions $g_2$ and $g_3$ are holomorphic at $q=0$, and moreover are defined over $\Z[\frac{1}{6}] \series{q}$ (making this so is the purpose of the powers of $2\pi i$ in the definition of the embedding). We can change coordinates, writing $x'=x- \frac{1}{12}$ and $2y' +x' =y$, so as to put the equation in the form $y'^2+x'y'= x'^3+a_4(q)x' + a_6(q)$. The benefit of the coordinate-change is that the coefficients now lie in $\Z\series{q}$. The series $a_4$ and $a_6$ are those given above---so the algebraic curve $y'^2+x'y'= x'^3+a_4(q)x' + a_6(q)$ is the Tate curve $\EuT$.

We conclude, then, that the specialized Tate curve $\hat{\EuT}$ is an elliptic curve, analytically isomorphic over $\C$ to the family $\Z\laurent{q}^*/  q^{\Z}$ when $0<|q|<1$. 

Its integrality is one interesting feature of $\EuT$, but another is that the absence of negative powers of $q$. One can therefore specialize $\EuT$ to $q=0$. The result is the curve $\EuT_0= \EuT|_{q=0}$ in $\mathbb{P}^2(\Z)$ given by
\begin{equation} \label{Tate eq q=0}
y^2+xy = x^3. 
\end{equation}

We can characterize this Weierstrass curve as follows:

\begin{lem}\label{Tate at q=0}
The curve $\EuT_0 \to \spec \Z$ has a section $s = [0:0:1]$ which is a node of $\EuT_0 \times_{\Z} \mathbb{F}_p$, the mod $p$ reduction of $\EuT_0$, for every prime $p$. Any Weierstrass curve $C\to \spec \Z$ possessing a section $s$ with this property can be transformed by integral changes of variable to
$\EuT_0$.
\end{lem}

\begin{pf}
Consider a general Weierstrass curve $C=[a_1,a_2,a_3,a_4,a_6]$, given as the projective closure of
\begin{equation}\label{Weier}
 y^2 + a_1 xy + a_3 y = x^3 + a_2 x^2 + a_4 x + a_6,\quad a_i \in \Z.  
\end{equation}
Integral points of $C\subset \mathbb{P}^2_{\Z}$, other than $[0:1:0]$, can represented as rational points on the affine curve. The point $[0:1:0]$ is regular over any field, and is the unique point of $C$ with $Z=0$. Suppose $[X:Y:Z]$ is an integral point that is nodal mod $p$ for all primes $p$. Then $Z$ must be non-zero mod $p$ for every prime $p$, and hence $Z$ is a unit of $\Z$.
Consider the $\Z$-point $(x_0,y_0)=(X/Z,Y/Z)$ of the affine curve. The partial derivatives vanish, since they vanish mod $p$ for all $p$:
\begin{equation}\label{crit} 
2y_0 + a_1 x_0 + a_3 = 0,\quad a_1 y_0  = 3x_0^2 + 2a_2 x_0 + a_4. 
\end{equation}
The nodal condition is that the Hessian is non-singular, that is,
\begin{equation}
\label{nodal}  a_1^2 +  2(6x_0+2a_2) \neq 0 \mod p.   
\end{equation}
(We note in passing that conditions (\ref{crit}, \ref{nodal}) hold for the point $[0:0:1]$ of $\EuT|_0$ at all primes $p$.) Since (\ref{nodal}) holds for all $p$, we have 
\begin{equation}
\label{criterion} 
a_1^2 + 12 x_0 + 4 a_2 = \pm 1.\end{equation} 
We shall use the criterion (\ref{criterion}) to make three changes of variable, successively making $a_1$, $a_2$ and $a_3$ equal to their counterparts for $\EuT_0$. 

First, (\ref{criterion}) tells us that $a_1$ is odd. Hence by a change of variable $x=x'$, $y = y' + c$, we may assume that $a_1=1$, whereupon $6x_0 + 2a_2 $ is either $0$ or $-1$. The latter possibility is absurd, so $3x_0 + a_2=0$. Being divisible by 3, $a_2$ can be removed altogether by a change of variable $x=x'+dy$, $y=y'$ without interfering with $a_1$. Thus we can assume additionally that $a_2=0$. 
We now find from (\ref{criterion}) that $x_0=0$. Hence $2y_0 + a_3=0$, so $a_3$ is even. It follows that $a_3$ can be set to zero by a change of variable $x=x'$, $y=y'+ e$, leaving $a_1$ and $a_2$ untouched.  Equations (\ref{crit}) now tell us that $y_0=0=a_4$, while the equation (\ref{Weier}) for $C$ tells us that $a_6=a_4^2=0$.
\end{pf}

More abstractly, if we define a curve $\pi\colon C \to \spec \Z$ by taking $\mathbb{P}^1_{\Z}$ and identifying the sections $[0:1]$ and $[1:1]$, so as to make every geometric fiber nodal, then the parametrization $\mathbb{P}^1_{\Z}\to \mathbb{P}^2_{\Z}$ given by $[s:t] \mapsto [st(s-t): s(s-t)^2 : t^3]$ identifies $C$ with $\EuT_0$.

\paragraph{Outline of method and algebraic results.}
This article is long partly because it contains rather more than a single proof of Theorem \ref{mainth}, and partly because working over $\Z$ presents significant technicalities beyond those that would be present if one worked over fields (or in some cases, of fields in which 6 is invertible).  Part I---a large chunk---is purely algebraic; it refines and elaborates the method of \cite{LP}. The basic point is that for any Weierstrass curve $C$, one has a 2-object subcategory $\EuB_C$ of $\perf C$---the dg category of perfect complexes of coherent sheaves---with objects $\EuO$ (the structure sheaf) and $\EuO_{p}$ (the skyscraper sheaf at the point at infinity), and this subcategory split-generates $\perf C$.  The cohomology category $A=H^*\EuB_C$ is independent of $C$, but the dg structure of $\EuB_C$ knows $C$. One can transfer the dg structure to a minimal $A_\infty$-structure on $A$. This procedure defines a functor from the category of Weierstrass curves to the category of minimal $A_\infty$-structures on $A$. We prove in Theorem \ref{Ainf comparison} that this functor is an equivalence. A slightly coarsened statement of Theorem \ref{Ainf comparison} is as follows:
\begin{thm}
Let $R$ be an integral domain which is either noetherian and normal of characteristic zero, or an arbitrary field. Let $(\EuB,\mu^*_{\EuB})$ be an $R$-linear $A_\infty$-category together with a Calabi--Yau structure of dimension 1. Assume that $\EuB$ is minimal, has just two objects $a$ and $b$, both spherical of dimension 1 and forming an $A_2$-chain (i.e. $\hom(a,a) \cong \Lambda^* (R[-1]) \cong \hom(b,b)$ as graded $R$-algebras; and $\hom(a,b)\cong R$,  $\hom(b,a)\cong R[-1]$ as  graded $R$-modules; and $\mu^1_\EuB=0$). Then $\EuB$ is trace-preservingly quasi-equivalent to $\EuB_C$ for a unique Weierstrass curve $C\to \spec R$, where $\EuB_C$ has the Calabi--Yau structure arising from its Weierstrass differential $\Omega \in  \Omega^1_{C/\spec R}$. 
\end{thm}
The proof of Theorem \ref{Ainf comparison} invokes the Hochschild cohomology $\Hoch^*(A,A)$. We computed this cohomology additively in \cite{LP}, but here we give a complete calculation, as a Gerstenhaber algebra, by interpreting $\Hoch^*(A,A)$ as the Hochschild cohomology $\Hoch^*(\Ccusp)$ of a cuspidal Weierstrass curve $\Ccusp$ (Theorem \ref{HH for cusp curve}). 

In Part II, we identify the unique curve $C_{\mathsf{mirror}}$ for which $\EuA_{C_{\mathsf{mirror}}}$ is quasi-isomorphic to the 2-object subcategory $\EuA_{\mathsf{symp}}$ of the Fukaya category $\EuF(T,z)$ on objects of slopes $0$ and $-\infty$, equipped with non-trivial double coverings. In \cite{LP}, we used Abouzaid's plumbing model \cite{AboPlumb} to prove that $\EuA_{\mathsf{symp}}|_{q=0}$ is not formal, which implies that $C_{\mathsf{mirror}}$ is not cuspidal. Here we identify $C_{\mathsf{mirror}}$ precisely. In fact, we identify the specialization ${C_{\mathsf{mirror}}}|_{q=0}$ in three independent ways: (i) by eliminating the possibility that $C_{\mathsf{mirror}}$ is smooth or cuspidal after reduction to an arbitrary prime $p$, by means of the `closed open string map' from symplectic cohomology to Hochschild cohomology of the Fukaya category; (ii) by calculating ``Seidel's mirror map'' \cite{Zas}, or more precisely, by determining the affine coordinate ring of ${C_{\mathsf{mirror}}}|_{q=0}$ via a calculation in the exact Fukaya category; and (iii) via theta-functions.  The third proof extends to a proof of mirror symmetry for $\EuF(T,z)$, not just its restriction to $q=0$.  We use an intrinsic model for the Tate curve, and the integral theta-functions for this curve which played a major role in Gross's proof \cite{Gro}. The nub is the multiplication rule for these theta-functions and its relation to counts of lattice-points in triangles.  The proof of mirror symmetry for the wrapped category is a rather formal extension of that for the exact category.

We should perhaps make one more remark about exposition. The authors' background is in symplectic topology. We imagine that typical readers will have an interest in mirror symmetry, perhaps with a bias towards the symplectic, algebro-geometric or physical aspects, but, like us, will not be expert in arithmetic geometry. We would be delighted to have readers who \emph{do} come from an arithmetic geometry background, but ask for their patience in an exposition which we fear belabors what is obvious to them and rushes through what is not.

\paragraph{Higher dimensions?}
We believe that there should be an arithmetic refinement to homological mirror symmetry for Calabi--Yau manifolds in higher dimensions, but will leave the formulation of such conjectures for elsewhere; the 2-torus is, we think, far from being an isolated case. The case of 2-tori with several basepoints can be treated rather straightforwardly starting from the one-pointed case, but we shall also leave that for another article.

\paragraph{Acknowledgements.}
YL was partly supported by the Herchel Smith Fund and by Marie Curie grant EU-FP7-268389;  TP by NSF grant DMS 1049313.  Paul Seidel provided notes outlining the algebro-geometric approach to computing Hochschild cohomology for $A$.  Conversations with several mathematicians proved useful as this work evolved, and we would particularly like to thank Mohammed Abouzaid, David Ben-Zvi, Mirela \c{C}iperiani, Brian Conrad, Kevin Costello, and Paul Seidel. We thank the Simons Center for Geometry and Physics for its generous hospitality.

\part{Algebraic aspects}

\section{Background material}
\subsection{\texorpdfstring{Derived categories and $A_\infty$-categories}{Derived categories and A-infinity categories}}
Our conventions and definitions are those of \cite[chapter 1]{SeiBook}; see \cite{Gro} for an informal introduction.  For now, we work over a ground \emph{field} $\BbK$ (commutative and unital), but we shall discuss presently more general ground rings. All our $A_\infty$-categories and functors are cohomologically unital. 

\paragraph{Triangulated envelopes.} Any $A_\infty$-category $\EuC$ has a triangulated envelope, a minimal formal enlargement that is a triangulated $A_\infty$-category, i.e., every morphism in $\EuC$ has a mapping cone in $\EuC$. The \emph{twisted complexes} $\tw \EuC$ of an $A_\infty$-category $\EuC$ form a model for the triangulated envelope. The cohomological category $H^0(\tw \EuC)$ is known as the derived category and denoted $\der \EuC$.

\paragraph{Split closure.}
One can formally enlarge $\tw\EuC$ further to another triangulated $A_\infty$-category $\twsplit \EuC$ which is additionally \emph{split-closed} (also known as idempotent-closed or Karoubi complete). An \emph{idempotent} in the $A_\infty$-category $\tw\EuC$ is defined to be an $A_\infty$-functor $\pi \colon \BbK \to \tw \EuC$ from the trivial $A_\infty$-category $\BbK$, which has one object $\star$ and $\hom(\star,\star)=\BbK$ (the ground field). For example, if the object  $X$ is the direct sum of objects $X_1$ and $X_2$, meaning that $\hom(\cdot, X) \cong \hom(\cdot,X_1)\oplus \hom(\cdot, X_2)$ in the category of $\tw \EuC$-modules, then $X_1$ defines an idempotent $\pi$ in $X_1\oplus X_2$, with $\pi(\star)=X$ and, on morphisms, $\pi(1) = \id_{X_1} \oplus 0_{X_2}$. The module $\hom(\cdot , X_1)$ is actually intrinsic to the idempotent $\pi$ (it can be constructed as the `abstract image of $\pi$' \cite[Chapter 1, (4b)]{SeiBook}); the object $X_1$ \emph{represents} the abstract image. Split-closed means that the abstract image of an arbitrary idempotent is represented by an object. We write $\dersplit \EuC$ for the triangulated category $H^0(\twsplit \EuC)$. It is useful to note that $\tw\EuC$ is split-closed as an $A_\infty$-category if and only if $H^0(\tw\EuC)$ is split-closed as an ordinary $\BbK$-linear category. 

\paragraph{Thomason's theorem.}
By \cite{Tho}, a necessary and sufficient condition for an $A_\infty$-functor which is a quasi-embedding to be a quasi-equivalence is that (a) it should induce a quasi-isomorphism after split-closure, and (b) that it should induce an isomorphism of Grothendieck groups $K_0$. Thus, clauses (iii) and (iv) from Theorem \ref{mainth}, which assert derived equivalence without split-closure, are partly statements about $K_0$.

\paragraph{$A_\infty$-categories over rings.}
Our Fukaya categories will be $A_\infty$-categories over unital commutative rings $\mathbb{L}$. The usual definition of an $A_\infty$-category $\EuC$ makes sense over such rings: the morphism spaces are arbitrary graded $\mathbb{L}$-modules. Let's call such an object a \emph{naive $A_\infty$-category}.  The basic notions carry through. For instance, the twisted complexes $\tw \EuC$, defined as usual (the multiplicity spaces are finite-rank free modules), form a triangulated envelope for $\EuC$, as in \cite[chapter 1]{SeiBook}. However, some of the naive constructions do not have the homotopical significance one might wish for. An example is that the Hochschild cohomology $\Hoch^*(A,A)$ of a $\mathbb{L}$-algebra $A$, defined through the bar complex, does not compute the bimodule-Ext module $\ext^*_{A^e}(A,A)$, but rather, relative $\ext$ for the map of $\mathbb{L}$-algebras $\mathbb{L} \to A^e$ \cite{Wei}. 

Over fields, $A_\infty$-constructions are automatically `derived'. To retain this property, we define a \emph{projective} $A_\infty$-category to be a naive $A_\infty$-category in which the morphism spaces are projective graded $\mathbb{L}$-modules. Fukaya categories are projective because the hom-spaces come with finite bases. 
Projective graded modules satisfy $\ext(V_1,V_2)=\Hom(V_1,V_2)$ and $\mathsf{Tor}(V_1,V_2)=V_1\otimes V_2$. The naive definitions of $A_\infty$-functors and their natural transformations, and of Hochschild homology and cohomology, work well for projective $A_\infty$-categories. 

\paragraph{DG categories over rings.} Differential graded (dg) categories over commutative rings have been well studied \cite{Kel}, and the theory does not depend on such ad hoc arrangements as having projective hom-spaces. There is a self-contained theory in which derived categories are defined via localization, not via twisted complexes. 

\paragraph{Calabi--Yau structures.} When $\BbK$ is a field, a \emph{Serre functor} of a $\BbK$-linear category $\EuC$ with finite-dimensional hom-spaces is an equivalence $S\colon \EuC\to \EuC$, together with isomorphisms
 \[ \phi_{A,B} \colon \Hom_{\EuC}(A,B) \simeq  \Hom_{\EuC}(B,SA)^\vee ,  \]
natural in both inputs, such that $S^\vee \circ \phi_{A,B} = \phi_{SA,SB}\circ S$ as maps $\Hom_\EuC(A,B)\to \Hom(SB,S^2 A)^\vee$ \cite{BK, BO}. A Serre functor $S_X$ for the (bounded) derived category $\der^b \coh(X)$ of a smooth projective variety $X$ over a field $\BbK$ is given by $S_X= \cdot \otimes \EuK_X [\dim X]$.  The maps $\phi_{X,Y}$ are Serre duality isomorphisms. 

A \emph{Calabi--Yau (CY) structure} of dimension $n$ on $\EuC$ is a Serre functor $(S,\phi)$ in which $S$ is the shift functor $Z\mapsto Z[n]$. If $X$ is a smooth projective Calabi--Yau variety, equipped with an $n$-form $\Omega$ trivializing $\EuK_X$, then its derived category has a CY structure induced by the isomorphism $\Omega \colon \EuO\to \EuK_X$. The role of $\Omega$ is to normalize the CY structure.

The cohomological (not derived) Fukaya category $H\EuF(M)$ of a (compact or exact) symplectic manifold $M^{2n}$, with coefficients in a field $\BbK$, comes with a natural Calabi--Yau structure: $\phi_{L_0,L_1}$ is the Floer-theoretic Poincar\'e duality isomorphism $HF(L_0,L_1)\cong HF(L_1,L_0)^\vee[n]$. It is a subtler matter to obtain a Calabi--Yau structure on the \emph{derived} Fukaya category $\der \EuF(M)$. It is expected that such a structure does exist, and is canonical, and arises from a cyclic symmetry defined on the $A_\infty$-level (see \cite{Fuk} for a construction of such a cyclic symmetry over $\R$, and \cite{KS} for an account of the relevant homological algebra), but this more refined structure will play no role in our considerations. 

A CY structure gives a canonical `trace map' $\tr_X = \phi_{X,X}(\id_X) \colon \Hom_\EuC^n(X,X)\to \BbK$. From the trace maps, one can reconstruct all the maps $\phi_{X,Y}$. In this article we think of CY structures in terms of their trace maps; a functor preserving CY structures will be called \emph{trace-preserving}. 

We shall need to say what we mean by a CY structure for a category over a commutative ring $\mathbb{L}$. The categories in question are of form $H^0\EuC$, where $\EuC$ is an $A_\infty$-category, and this permits us to make an expedient (but not fully satisfactory) definition:
\begin{defn} 
A \emph{CY structure} consists on the $\mathbb{L}$-linear $A_\infty$-category $\EuC$ consists of cochain-level maps
\[ \phi_{A,B} \colon \hom_{\EuC}(A,B) \simeq  \hom_{\EuC}(B,A[n])^\vee \]
such that the induced maps on cohomology
\[  [\phi_{A,B}\otimes 1_{\mathbb{F}}] \colon \Hom_{H^0(\EuC\times_{\mathbb{L}} \mathbb{F})}(A,B) \simeq  \Hom_{H^0(\EuC\times_{\mathbb{L}} \mathbb{F})}(B,A[n])^\vee  \]
form a CY structure for each residue field $\mathbb{L} \to \mathbb{F}\to 0$. (Note that since the hom-space in $\EuC$ are projective modules, they are also flat, so tensoring them with $\mathbb{F}$ commutes with $H^0$.) If $\EuC$ and $\EuD$ have CY structures, an $A_\infty$-functor $\psi\colon \EuC\to \EuD$ is called \emph{trace-preserving} if the induced functors $H^0(\EuC\otimes \mathbb{F})\to H^0(\EuD\otimes \mathbb{F})$ are all trace-preserving.
\end{defn}
With this definition, Fukaya categories have CY structures over arbitrary rings $\mathbb{L}$, since Poincar\'e duality is defined at cochain level but our demands on the maps are all at cohomology-level.\footnote{This does \emph{not} apply to wrapped Fukaya categories.} 
 
\paragraph{Perfect complexes.}  Let $X$ be a scheme. A \emph{strictly perfect complex} is a bounded complex of locally free, finite rank $\EuO_X$-modules. A \emph{perfect complex} is a cohomologically bounded complex $\EuP^\bullet$ of coherent sheaves of $\EuO_X$-modules which is locally quasi-isomorphic to a strictly perfect complex. Inside the bounded derived category of coherent sheaves $\der ^b\coh(X)$, one has a full triangulated subcategory $\perf (X)$ of perfect complexes.

We will need to consider dg enhancements of $\perf(X)$; that is, we want a pre-triangulated dg category $\EuC$ and an equivalence of triangulated categories $\varepsilon\colon H^0(\EuC)\to \perf(X)$.  When $X$ is a projective scheme over a field $\BbK$, $\perf(X)$ has a dg enhancement $(\EuC,\varepsilon)$ which is unique: if $(\EuC',\varepsilon')$ is another then there is a \emph{quasi-functor} $\phi\colon \EuC\to \EuC'$ such that $\varepsilon' \circ H^0(\phi) = \varepsilon$ \cite{LO}.  Since we wish to work over more general base rings, and for computational purposes, we specify a dg enhancement of $\perf(X)$, valid for $X$ a projective noetherian scheme, as follows. 

Assume $X$ is separated and noetherian. Fix an affine open covering $\EuU$ of $X$. Define a dg category $\VB(X)$ whose objects are locally free sheaves (=vector bundles)  of finite rank, and whose hom-spaces, denoted $\RR\hom^\bullet(\EuE,\EuF)$, are  \v{C}ech complexes: 
\[  \RR\hom^\bullet(\EuE,\EuF) = \left( \check{C}^\bullet(\EuU; \SheafHom(\EuE,\EuF)),\delta\right), \]
with $\delta$ the \v{C}ech differential. The cohomology of the \v{C}ech complex is
\[  \RR\Hom^\bullet(\EuE,\EuF) =\check{H}^\bullet(\EuU; \SheafHom(\EuE,\EuF)) \cong \ext^\bullet(\EuE,\EuF) \] 
by \cite[Theorem III.4.5]{Har} and the fact that $\ext^\bullet(\EuE,\EuF)\cong H^\bullet(\EuE^\vee\otimes \EuF)$. Composition combines the shuffle product of \v{C}ech cochains with the composition of  sheaf-morphisms. Whilst $\VB(X)$ depends on the open covering, different choices lead to quasi-isomorphic dg categories (take the union of the two coverings). We now pass to the pre-triangulated dg category $\tw \VB(X)$ of twisted complexes. There is an embedding $H^0 (\tw \VB(X))\to  \perf(X)$, mapping a twisted complex to its total complex. This embedding is a quasi-equivalence, because every perfect complex is quasi-isomorphic to a strictly perfect complex \cite[Prop. 2.3.1(d)]{TT}. 

Another approach to dg enhancement is to use injective resolutions; the equivalence of the injective and \v{C}ech approaches is shown in \cite[Lemma 5.1]{SeiQuartic}, over fields; the proof remains valid over rings.

Grothendieck--Serre duality defines a CY structure for $\tw \VB(X)$ when $X$ is equipped with a trivialization of the relative dualizing sheaf $\omega_{X/\mathbb{L}}$. In Theorem \ref{mainth}, the functor $\psi$ is trace-preserving.

\subsection{Geometry of Weierstrass curves} 

\subsubsection{Genus-one curves} We shall need to work
with curves over the rings $\Z$ and $\Z\series{q}$, and to this end we note
some terminology for curves over schemes (cf. for example \cite{DR, ConEll}). A
\emph{curve} over a noetherian scheme $S$ is a morphism of schemes $\pi\colon C
\to S$ that is separated, flat and finitely presented, such that for every
closed point $s \in S$ the fiber $C_s$ is non-empty of pure dimension 1. The
Euler characteristic $\chi(C_s,\EuO_{C_s})$ is then locally constant; when it
is constant and equal to $1-g$, and the geometric fibers are connected, we say that $C$ has arithmetic genus $g$.  

We shall always apply the restrictions that curves are to be proper,  and
that the fibres $C_s$ are Cohen--Macaulay. This implies that one has  a dualizing sheaf $\omega_{C/S}$, and where $C\to S$ is regular it coincides with the sheaf of differentials $\Omega^1_{C/S}$.
A reminder on duality \cite{Con}: there is an intrinsic residue isomorphism of
sheaves on $S$ 
\[  \mathsf{res} \colon \RR^1f_*(\omega_{C/S}) \to \EuO_S.  \] 

With the Yoneda (composition) product $\smile$, this defines the Serre duality
pairing, 
\[   \RHom_S^{1-i}( \EuF, \omega^1_{C/S}) \otimes  \RR^i f_*(\EuF)
\stackrel{\smile}{\to}  \RR^1f_*(\omega^1_{C/S}) \stackrel{\mathsf{res}}{\to} \EuO_S,
\] 
for any coherent sheaf $\EuF$.

A curve has arithmetic genus one if and only if $\EuO_C \cong \omega_{C/S}$, i.e., if and only if $\EuO_C$ is a dualizing sheaf. If $\omega \colon \EuO_C\to \omega_{C/S}$ is an isomorphism then it composes with the residue map to give an isomorphism
\[  \tr_\omega \colon \RR^1f_*(\EuO_C) \to \EuO_S,  \]
and a Serre duality pairing
\[   \RHom_S^{1-i}( \EuF, \EuO_C) \otimes  \RR^i f_*(\EuF) \stackrel{\smile}{\to}  \RR^1f_*(\EuO_C) 
\stackrel{\tr_\omega}{\to} \EuO_S \]
which induces a perfect pairing on stalks at any closed point $s\in S$.

\subsubsection{Weierstrass curves: definitions} 
\begin{defn}
An \emph{abstract Weierstrass curve} $(C,\sigma ,\omega)$ over $S$ is a curve $C\to S$ of arithmetic genus one, such that each geometric fiber $C_s$ is irreducible, equipped with a section $\sigma \colon S\to C$ of $\pi$ and a specific isomorphism $\omega \colon \EuO_C\to \omega_{C/S}$. An \emph{isomorphism} of abstract Weierstrass curves $(C_1,\sigma_1,\omega_1)$ and $(C_2,\sigma_2,\omega_2)$ over $S$ is an isomorphism $f\colon C_1\to C_2$ of $S$-schemes such that $f\circ \sigma_1 = \sigma_2$, and such that the map $f^*\omega_2 \colon f^*\EuO_{C_2}\to f^*\omega_{C_2/S}$ coincides with $\omega_1$ under the identifications $f^* \EuO_{C_2} \cong \EuO_{C_1}$ and $f^*\omega_{C_2/S} \cong \omega_{C_1/S}$ induced by $f$.
\end{defn}

\begin{defn}
An \emph{embedded Weierstrass curve} over $S=\spec R$ is a curve $C \subset \mathbb{P}^2_S$ embedded as a cubic 
\begin{equation} \label{Weierstrass equation}
y^2 + a_1 xy + a_3 y= x^3 + a_2 x^2 + a_4 x + a_6\quad (a_i\in R).
\end{equation}
\end{defn}
Such a curve comes with its point at infinity $p=[0:1:0]$, which defines a section $\sigma$ of $C\to S$. It also comes with a standard differential $\omega$, possibly with poles at the singular points: Writing the cubic equation as $w(x,y):=y^2 - x^3 + \cdots  =0$, one has $ \omega = dx/w_y  $ at points where $w_y\neq 0$, and $\omega = -dy/w_x$ at points where $w_x\neq 0$. 

\begin{lem}
Assume that $R$ is a \emph{normal} ring (i.e., $R$ is reduced and integrally closed in its total quotient ring). Then $\omega$ defines a section of the dualizing sheaf $\omega_{C/R}$. 
\end{lem}
\begin{pf}
Let $R[\underline{a}] = R[a_1,a_2,a_3,a_4,a_6]$---another normal ring. It will suffice to prove the assertion for the `universal Weierstrass curve' $p\colon \EuC \to \spec R [\underline{a}]$ defined by (\ref{Weierstrass equation}), since the formation of the dualizing sheaf is compatible with the specialization to particular values of the $a_i$.

The scheme $\EuC$ is normal: in the open set $U$ where (\ref{Weierstrass equation}) is valid, $a_6$ is a function of the other variables, so projection $U\to \spec R[x,y, a_1,a_2,a_3,a_4]$ is an isomorphism, and $R[x,y, a_1,a_2,a_3,a_4]$ is normal. Along the section at infinity $\sigma=[0:1:0]$, the fibers of $p$ are regular, and the base normal, so the total space is normal. The relative dualizing sheaf $\omega_{\EuC /R [\underline{a}]}$ is an invertible sheaf, since all its fibers are Gorenstein (being local complete intersections). The locus where the fibers of $p$ are singular is defined by $w_x = w_y = 0$. This locus has codimension 2 in $\EuC$: it maps to the codimension 1 locus $\{\Delta=0\} \subset \spec R[\underline{a}]$ defined by the vanishing of the discriminant, and it has codimension 1 in each fiber. Since $\omega$ is a section of $\omega_{\EuC/R [\underline{a}]}$  defined outside a codimension 2 subset of a normal scheme, it extends to a global section, by the algebraic counterpart to Hartogs's theorem.\footnote{One can take this to be the statement that an integrally closed subring $A$ of a field $K$ is the intersection of the valuation rings in $K$ which contain $A$ \cite[5.22]{AM}.}
\end{pf}
Thus an embedded Weierstrass curve functorially defines an abstract Weierstrass curve $(C,\sigma,\omega)$. By Riemann--Roch, every abstract Weierstrass curve is isomorphic to an embedded one. To specify the embedding into $\mathbb{P}^2$, one must give a basis of $ \h^0(\EuO_C(3\sigma))$ of the form $(1,x,y)$, where $1$ is the regular function with value 1, and $x\in  \h^0(\EuO_C(2\sigma))$. The denominator-free form of the argument is given at \cite[p. 68]{KM}, for instance. 

\subsubsection{Reparametrization group}\label{reparam}
The algebraic group $G\subset PGL(3)$ of elements which preserve Weierstrass form consists of matrices (up to scale) of the shape
\begin{equation} \label{def G}
\left [ \begin{array} {ccc} 
  u^3 	&  s 	& t \\
  0  	& u^2 	&   r \\
  0 	& 	0 		&   1
 \end{array}\right],\quad u \in \mathbb{G}_m. \end{equation}
 We shall call $G$ the \emph{reparametrization group} for embedded Weierstrass curves. It acts on embedded Weierstrass curves via the substitutions
\[  x=u^2 x'+ r,\quad y=u^3 y' + u^2s x' + t.   \]
The effects of a substitution  on the Weierstrass coefficients are listed in \cite{Del} or \cite{Sil}:
\begin{align}\label{subst}
ua_1'  & = a_1 + 2s  \\
u^2 a_2'  & = a_2 -s a_1 + 3r-s^2 \\
u^3 a_3'  & = a_3 + ra_1 +2t  \\
u^4 a_4'  & = a_4 - sa_3 + 2ra_2  -(t+rs) a_1 + 3r^2 -2st   \\
u^6 a_6'  & = a_6 +ra_4 + r^2a_2+ r^3 -ta_3 -t^2 -rta_1.
\end{align}
The unipotent subgroup $U \leq G$ of elements where $u=1$ is the subgroup which preserves the differential $\omega$. Thus if $g\in U$ then $g\colon C\to g(C)$ is an isomorphism of abstract Weierstrass curves.

\paragraph{The Lie algebra.} The Lie algebra $\frakg$ of $G$ is spanned by four vectors:
\[  \partial_s:=
\left [ \begin{array} {ccc} 
  0 		&   1 	&  0 \\
  0  	& 	0 	&   0 \\
  0 		& 	0 		&   0
  \end{array}\right],\;
  \partial_r:=
\left [ \begin{array} {ccc} 
  0 		&   0 	&  0 \\
  0  	& 	0 	&   1 \\
  0 		& 	0 		&   0
  \end{array}\right],\;
\partial_t:=
\left [ \begin{array} {ccc} 
  0 		&   0 	&  1 \\
  0  	& 	0 	&   0 \\
  0 		& 	0 		&   0
  \end{array}\right],\]
(these three span the Lie algebra $\mathfrak{u}$ of $U$) and 
\[
  \partial_u:=
\left [ \begin{array} {ccc} 
  3 		&   0 	&  0 \\
  0  	& 	2 	&   0 \\
  0 		& 	0 		&   0
  \end{array}\right].
  \] 
  The derivative of the $G$-action on 
 \begin{equation}\label{def W}
  W: = \spec\BbK[a_1,a_2,a_3,a_4,a_6] 
  \end{equation} 
  is an action of $\frakg$ on $W$ by a Lie algebra homomorphism
 \begin{equation}\label{Lie action} \rho \colon \frakg \to \vect(W),\end{equation} 
 which we think of as a map $\rho \colon \frakg \times W\to W$.  The partial derivative $ (\partial \rho /\partial w)|_{w=0}\colon \frakg\times W\to W$ makes $W$ a $\mathfrak{g}$-module. We can form a differential graded Lie algebra (DGLA) concentrated in degrees $0$ and $1$,
\begin{equation}\label{def L}  
\EuL = \{ \frakg \stackrel{d}{\to} W \} ,\quad d(\xi) =  \rho(\xi, 0),
\end{equation}
whose bracket combines the Lie bracket of $\frakg$ with the module structure of $W$. Thus $\EuL$ captures the truncation of $\rho$ where we only work in a first-order neighborhood of $0\in W$.

There are $\BbK^\times$-actions on $\frakg$ and on $W$, intertwined by $d$. The action on $W$ is given by $\tau \cdot a_j = \tau^{-j} a_j$; that on $\frakg$ by $\tau \cdot \partial_s = \tau^{-1} \partial_s$, $\tau \cdot \partial_r = \tau^{-2}\partial_r$, $\tau \cdot \partial_t = \tau^{-3}\partial_t $, $\tau\cdot \partial_u = \partial_u$. Thus $W$ and $\frakg$ are graded $\BbK$-modules.

Explicitly, taking $(\partial_s,\partial_r,\partial_t,\partial_u)$ as basis for $\frakg$, and 
$(a_1,a_2,a_3,a_4,a_6)$ as coordinates for $W$, one has 
\[  d  = 
\left[\begin{array}{cccc} 
2 & 0 &0 & 0\\
0 & 3 &0 & 0\\
0 &0 & 2 & 0 \\
0 &0 &0 &  0\\
0 &0 &0 &  0
\end{array} \right]\] 
and
\begin{align}\label{coker d} 
W^{\frakg} :=   \coker d  
& = \stackrel{a_1}{\frac{\BbK}{(2)}[1]  } \oplus \stackrel{a_2}{\frac{\BbK}{(3)} [2]}  \oplus \stackrel{a_3}{\frac{\BbK}{(2)} [3]}  \oplus\stackrel{a_4}{\BbK [4]} \oplus \stackrel{a_6}{\BbK [6]},  \\
 \label{ker d}
 \ker d & = 
\stackrel{\partial_u}{\BbK} \; \oplus \; \stackrel{\partial_s}{\frac{\BbK}{(2)}[1]} \; \oplus \; \stackrel{\partial_r}{\frac{\BbK}{(3)}[2]} \; \oplus \; \stackrel{\partial_t}{\frac{\BbK}{(2)}[3]}. 
\end{align}

\subsubsection{The cuspidal cubic}
The cuspidal Weierstrass curve 
\begin{equation}\label{def cusp}
\Ccusp=\{y^2-x^3=0\}\end{equation}
will play a special role in our story, stemming from the fact that the full subcategory of its derived category whose objects are the structure sheaf and the skyscraper at infinity is \emph{formal}. Let $X=\mathbb{P}^1$, and let $p \colon \spec\BbK \to X$ be the $\BbK$-point  $[0:1]$. Let $z $ denote the standard affine coordinate $\mathbb{A}^1 \to \mathbb{P}^1$,  $z\mapsto [z:1]$. One has the structure sheaf $\EuO_X$, and inside it the sheaf $\EuO^p \subset \EuO_X$ of functions $f$ such that $Df(p)=0$ (that is, in terms of the local coordinate $z$, functions $f(z) = f(0) + O(z^2)$). Let $X_{\mathsf{cusp}}$ denote the scheme $(X,\EuO^{p})$.
\begin{lem}\label{abstract cusp}
The abstract Weierstrass curve underlying $\Ccusp$ is isomorphic to 
\[ (X_{\mathsf{cusp}}, p , z^{-2}d z),\] 
\end{lem}
\begin{pf}
The normalization of $\Ccusp$ is a non-singular rational curve $\tilde{C}$ with a distinguished point $c$ which maps to the cusp under the normalization map $\nu\colon \tilde{C}\to \Ccusp$. We fix an isomorphism $X \to \tilde{C}$ mapping $p$ to $c$. The map $\nu$ is a homeomorphism in the Zariski topology, and so defines a scheme-theoretic isomorphism $(\tilde{C}, \nu^* \EuO_{\Ccusp})\to (\Ccusp,\EuO_{\Ccusp})$. One has $\EuO_{\Ccusp} \cong \EuO^p$: the local model near the cusp is the map of $\BbK$-algebras $\BbK[x,y]/(y^2-x^3) \to \BbK[z]$ given by $x\mapsto z^2$ and $y\mapsto z^3$, whose image is $\BbK.1 \oplus z^2\BbK[z]$. 

The $\EuO^p$-module of differentials $\Omega^1_{X_{\mathsf{cusp}}}$ is given by the submodule of $\Omega^1_{X}(2p)$ (meromorphic differentials on $X$ with a double pole at $p$) formed by the differentials with vanishing residue at $p$. In terms of the affine coordinate $z$ near $p$, the differential of a function $g(z)=a+ bz^2+ \dots $ is $d g = g' (z) dz$. The Weierstrass differential $\omega$ is given by $\omega=dx/(2y) = dy/(3x^2)$ (in characteristics 2 and 3 only one of these expressions makes sense). In terms of $z$,  one has $\omega =  z^{-2} dz$; this makes global sense because $\omega = - d(z^{-1})$.
\end{pf}

\section{Perfect complexes on Weierstrass curves}

\subsection{The two-object dg category associated with a Weierstrass curve}
In this subsection we explain how to pass from a Weierstrass curve $C\to \spec R$ to a two-object dg category $\EuB_C$ with standard cohomology. Consider a genus-one curve $C$ over a noetherian affine scheme $S$. It has a dg category $\VB(C)$, defined via an affine open covering, linear over the ring $\EuO_S$.
The dg category $\tw \VB(C)$ for an abstract Weierstrass curve $(C,\sigma,\omega)$ over $\spec R$ has extra structure in the form of a trace pairing $\tr$, as described in the introduction. It also has distinguished split-generators, namely, the structure sheaf and the skyscraper $\EuO_{C,\sigma}=\sigma_* \EuO_{\spec R}$ at $\sigma$ (more properly, its locally-free resolution $\EuO\to \EuO(\sigma)$): 
\begin{lem}\label{gen perf}
For a Weierstrass curve $C\to \spec R$ over a noetherian affine scheme, one has 
\[ \tw \VB(C) = \langle \EuO_C , \EuO_{C, \sigma} \rangle.  \]
\end{lem}
Here $\langle\cdot \rangle $ denotes the smallest dg subcategory of  $\tw \VB(C)$ closed under quasi-isomorphisms, shifts, mapping cones and passing to idempotents. 
\begin{pf}
We claim first that, if $\EuO(1)$ is a very ample line bundle on $C$, then $\tw \VB(C)$ is split-generated (i.e., generated under quasi-isomorphisms, shifts, mapping cones and passing to idempotents) by the twists $\{ \EuO(n) \}_{n<0}$. The argument is as in \cite[Lemma 5.4]{SeiQuartic}, which Seidel attributes to Kontsevich. Take  a locally free sheaf $V$ on $C$. By Serre's theorem that very ample implies ample \cite[Theorem II.5.17]{Har}, which is valid for noetherian projective schemes, one can find an epimorphism $\EuO(m)^{\oplus r} \to V$ for some $m \ll 0$. Iteratively, one can find for each $k$ a left resolution
\[   0\to V'\to  \EuO(m_k)^{\oplus r_k} \to \dots \to \EuO(m_1)^{\oplus r_1} \to V \to 0. \]
There results an exact triangle in $\der(C)$
\[    \{ \EuO(m_k)^{\oplus r_k} \to \dots \to \EuO(m_1)^{\oplus r_1} \} \to V \to  V'[k] \stackrel{+}{\to} . \]
Now, $\ext^k_R(V,V')=H^k(C,V^\vee\otimes V')$, and if we take $k> \dim C = 1+ \dim R$, this Ext-module must vanish. Consequently, the exact triangle splits and defines a quasi-isomorphism
\[   \{ \EuO(m_k)^{\oplus r_k} \to \dots \to \EuO(m_1)^{\oplus r_1} \} \to V \oplus  V'[k] . \]
Thus $V$ is a direct summand in the object on the LHS. Note also that $V'$ is a perfect complex, because it is the mapping cone of a map of perfect complexes; therefore, $V'$ is quasi-isomorphic to a strictly perfect complex. This proves that every locally free sheaf lies in the split-closure of the collection $\{ \EuO(n) \}_{n<0}$. It follows that the same is true of every object of $\tw \VB(C)$. 

Note next that $\langle \EuO_C , \EuO_{C, \sigma} \rangle$ includes $\EuO_C(n\sigma)$ for each $n \leq 0$, by a straightforward induction.  But $\EuO(3\sigma)$ is a very ample line bundle, so now the claim completes the proof.
\end{pf}

\paragraph{A two-object subcategory.}
Let $\EuB_{C}$ denote the full dg subcategory of $\tw \VB(C)$ with the two objects $ \EuO_C$ and $\EuO_{C, \sigma}$ and with the trace map $\tr_\omega$. It is defined up to quasi-isomorphisms acting trivially on cohomology. To be precise, we shall define $\EuB_C$ using the \v{C}ech complexes associated with an affine open covering $\EuU$. If we pick two coverings $\EuU_1$ and $\EuU_2$, we get dg categories
$\EuB_{\EuU_1}$ and $\EuB_{\EuU_2}$, and a zigzag of quasi-isomorphisms
\[  \EuB_{\EuU_1} \leftarrow \EuB_{\EuU_1\cup \EuU_2}\to \EuB_{\EuU_2}.  \]
The cohomology category $H^*\EuB_C$ is truly canonical---defined up to canonical isomorphism.

An isomorphism $g \colon C_1\to C_2$ of abstract Weierstrass curves is a homeomorphism $g$, together with a local isomorphism of sheaves of $S$-modules $g^\# \colon \EuO_{C_2}\to g_* \EuO_{C_1}$, respecting the sections and differentials. That means, first, that $\sigma_2 \colon S\to C_2$ is the composite $g\circ \sigma_1$; this implies a canonical isomorphism $\EuO_{C_2,\sigma_2}\to g_* \EuO_{C_1,\sigma_1}$. The isomorphism $g$ induces isomorphisms between the abelian categories of coherent sheaves on $C_1$ and $C_2$, preserving the objects $\EuO$ and $\EuO_\sigma$. This naturally extends to an isomorphism of dg categories $g_*\colon \EuB_{C_1}\to \EuB_{C_2}$, provided that we use an open covering $\EuU$ for $C_1$ and $g(\EuU)$ for $C_2$. Thus, if we have $g_{12}\colon C_1\to C_2$ and $g_{23}\colon C_2\to C_3$ with composite $g_{13}$, then the composite isomorphism $g_{23*}\circ g_{12*}\colon \EuB_{C_1}\to \EuB_{C_3}$ coincides with $g_{13*}$, provided again that we use the coverings $\EuU$, $g_{12}(\EuU)$ and $g_{13}(\EuU)$. If we do not, then we get instead the formal composite of chains of quasi-isomorphisms:
\[ \xymatrix{
&\EuB_{g_{12}(\EuU_1)\cup U_2} \ar[dr]\ar[d]  && \EuB_{g_{23}(\EuU_2) \cup \EuU_3}\ar[d]\ar[dr] \\
\EuB_{\EuU_1} \ar[r] & \EuB_{g_{12}(\EuU_1)} &\EuB_{\EuU_2} \ar[r] & \EuB_{g_{23}(\EuU_2)} & \EuB_{\EuU_3}
} \]

\paragraph{The cohomology category.}
We shall be interested in the map which assigns to each Weierstrass curve $(C,\omega,\sigma)$ a graded-linear cohomology category $A_C$ and a dg category-with-trace,
\[ (C,\omega,\sigma) \mapsto \EuB_C, \]
with an isomorphism $H^*\EuB_C\cong A_C$, defined up to quasi-isomorphisms acting trivially on $A_C$.

\begin{prop}
The category $A_C$ is independent of the abstract Weierstrass curve. Precisely: There is an $\EuO_S$-linear graded category $A$ with two objects $\EuO$ and $\EuO_\sigma$, equipped with a trace map $\tr$, such that the following holds: For any abstract Weierstrass curve $(f\colon C\to S,\omega,\sigma)$ the cohomology category $A_C=H^*(\EuB_C)$ is trace-preservingly isomorphic to $A$ in such a way that if $C_1\to C_2$ is any isomorphism of Weierstrass curves then the resulting map $A\to A$ is the identity.
\end{prop}

In other words, the category of Ext-modules between $\EuO_C$ and $\EuO_{C,\sigma}$ is independent of $(C,\sigma,\omega)$ as a graded $\BbK$-linear category with trace. 

To prove the proposition, we examine the structure of $H^*\EuB_C$. Writing $\EuO=\EuO_C$ and $\EuO_\sigma=\EuO_{C,\sigma}:=\sigma_* \EuO_S$, one has canonical isomorphisms
\begin{align*}
& \RHom_S(\EuO,\EuO) \cong \RR f_*(\EuO) \cong \EuO_S \oplus \RR^1f_*(\EuO),\\
& \RHom_S(\EuO_{\sigma},\EuO_\sigma) \cong \Lambda^* (\sigma^*TC) \cong \EuO_S \oplus 
\sigma^*\EuT_C
\end{align*}
Thus both endomorphism spaces are 2-dimensional and sit in degrees $0$ and $1$.
The trace isomorphisms
\begin{align*} 
& \tr_\omega \colon \RR^1f_*(\EuO) \to \EuO_S \\
& \tr_\omega \colon  \sigma^*\EuT_C\to  \EuO_S
\end{align*}
are, in the first case, the one that we have discussed (the composite of $\omega$ and the residue pairing) and in the second case the pullback by $\sigma$ of the composite 
\[ 	\EuT_C \to \EuT_C \otimes_{\EuO} \EuO
\xrightarrow{\id \otimes \omega} \EuT_C \otimes \omega_{C/S}
\xrightarrow{\mathsf{ev}_\sigma} \EuO,
\] 
where $\mathsf{ev}_\sigma$ is the map defined by evaluating $\omega$---viewed as a differential---on tangent vectors at $\sigma$. One has
\[
 \RHom_S(\EuO,\EuO_{\sigma}) =  \RR^0 \Hom_S(\EuO,\EuO_\sigma) \cong f_*(\EuO_\sigma)\cong \EuO_S.
\]
Finally, one has isomorphisms
\begin{align*}
\RHom_S(\EuO_{\sigma}, \EuO) & = \RR^1 \Hom_S(\EuO_{\sigma}, \EuO) \\
& \cong  \RR^1 \Hom_S(\EuO_{\sigma}, \EuO) \otimes_{\EuO_S} \EuO_S \\
& \cong  \RR^1 \Hom_S(\EuO_{\sigma}, \EuO) \otimes_{\EuO_S} \RR^0 \Hom_S(\EuO,\EuO_\sigma)\\
& \xrightarrow{\smile} \RR^1 \Hom_S(\EuO_\sigma ,\EuO_\sigma)\\
& \xrightarrow{\tr_\omega}  \EuO_S.
\end{align*}

We now describe the category $A$ demanded by the proposition above. Let $X = \EuO$ and $Y=\EuO_\sigma$. We have seen how to use $\omega$ to obtain algebra isomorphisms
\[ \End(Y) \cong \Lambda^*(\EuO_S[-1]) \cong \End(Y)  \]
such that the trace maps correspond to the identity map of $\EuO_S$. We also have exhibited isomorphisms $\Hom(X,Y) = \EuO_S$ and $\Hom(Y,X) = \EuO_S[-1]$. The composition maps are mostly dictated by the requirements of grading and unitality. The interesting ones are
\[ 	\Hom^1(Y,X)\otimes \Hom^0(X,Y) \to \Hom^1(Y,Y), \quad
  	\Hom^0(X,Y)\otimes \Hom^1(Y,X) \to \Hom^1(X,Y). \]
These are both given by the multiplication of functions
\[ \EuO_S \otimes  \EuO_S\to \EuO_S . \] 
The objects $X$ and $Y$ and their morphisms form a graded-linear CY category $(A,\tr)$, independent of $C$.

\paragraph{Differential graded structure.}
While the cohomology category $A=A_C$ is independent of $C$, the dg structure of $\EuB_C$ is fully sensitive to the curve $C$:

\begin{thm}[dg comparison theorem]\label{DG comparison}
Work over a field $\BbK$.
\begin{enumerate}
\item
Let $\EuB$ be a dg category with trace such that $H^*(\EuB) \cong A$. Then  there exist an abstract Weierstrass curve $(C,\omega,\sigma)$ and a trace-preserving $A_\infty$-quasi-isomorphism $\EuB \to \EuB_C$.
\item
If $(C,\omega,\sigma)$ and $(C',\omega',\sigma')$ give rise to quasi-isomorphic dg categories with trace, i.e., $\EuB_C$ is related to $\EuB_{C'}$ by a zig-zag of trace-preserving isomorphisms, then $(C,\sigma,\omega)\cong (C',\sigma',\omega')$.
\end{enumerate}
\end{thm}
We state this result now so as to indicate our aims. However, we will establish it as a corollary of a more detailed statement, Theorem \ref{Ainf comparison}, and it will in fact be the latter result which we use, not Theorem \ref{DG comparison}.

\begin{rmk}
The proofs will be given later, but we offer two hints. For the uniqueness clause, the point is that there is a construction which assigns to any such category $\EuB$ a sequence $\EuT_n$ of twisted complexes   in a uniform manner. When $\EuB=\EuB_C$, one has $\EuT_n\simeq \EuO_C(n\sigma)$. One further constructs multiplication maps $ \h^0(\EuT_n)\otimes  \h^0(\EuT_m)\to  \h^0(\EuT_{m+n})$. When $\EuB=\EuB_C$, these reproduce the multiplication $ \h^0(\EuO(m)) \otimes  \h^0(\EuO(n)) \to  \h^0(\EuO(m+n))$. Thus the coordinate ring of the affine curve $C^\circ$, the open complement of $\im\sigma$, is determined by $\EuB_C$.  The existence clause (1) is plausible because one has $ \h^0(\Lambda^2 T^*_C)=0$ and $H^2(\EuO_C)=0$. As a result, $\perf (C)$ has no Poisson deformations and no non-commutative deformations, and it is reasonable to expect all deformations of $ \perf (C)$ to be geometric. 
\end{rmk}

The cuspidal cubic $\Ccusp= \{ y^2-x^3=0 \}$ has the following special property, which already appeared in \cite{LP}:
\begin{lem}\label{formality}
The dga $\EuB_{\mathsf{cusp}}:=\EuB_{\Ccusp}$ is formal. 
\end{lem}
\begin{pf}
We may transfer the dg structure of $\EuB_{\mathsf{cusp}}$ to a minimal $A_\infty$-structure on $A=H^*\EuB_{\mathsf{cusp}}$. The transfer of dg structure will be described in detail in the proof of Lemma \ref{splittings}. The goal, then, is to prove that the $A_\infty$ structure maps $\mu^d$ vanish for $d>2$. 

$\Ccusp$ is the curve $\{Y^2Z = X^3\} \subset \mathbb{P}^2$. The multiplicative group $\mathbb{G}_m$ acts on $\Ccusp$ by $t\cdot (X,Y,Z) = (t^{-2}X,t^{-3}Y,Z)$, preserving the point $\sigma=[0:1:0]$, and therefore acts on $\EuB_{\mathsf{cusp}}$. The action of $\mathbb{G}_m$ on $\EuB_{\mathsf{cusp}}$  induces an action on the cohomology $A$, and the transfer will be set up equivariantly so that the resulting $A_\infty$-structure has the property that $\mu^d(t\cdot a_d,\dots,t\cdot a_1) = t\cdot \mu^d (a_d,\dots,a_1)$. 

A short computation leads to the following conclusion: the weight of the $\mathbb{G}_m$-action on a hom-space $\Hom^k_A(X,X')$ (where $X$ is $\EuO$ or $\EuO_\sigma$, ditto $X'$) is equal to the degree $k$. Now take $(X_0,\dots, X_d)$ a sequence of objects ($\EuO$ or $\EuO_\sigma$), and take $a_j\in \hom^{k_j}_{\EuB_{\mathsf{cusp}}}(X_{j-1},X_j)$. For the equation $\mu^d(t\cdot a_d,\dots,t\cdot a_1) = t\cdot \mu^d (a_d,\dots,a_1)$ to hold, one must have
$ k_1+\dots + k_d+2-d = k_1 + \dots + k_d, $ i.e., $d=2$.
\end{pf}

\subsection{\texorpdfstring{Stable vector bundles on $\EuT_0$}{Stable vector bundles on T0}}
When we come to prove Theorem \ref{mainth}, clause (iii), we will need to apply Thomason's theorem about Grothendieck groups \cite{Tho}, and for that we shall need to know $K_0(\EuT_0)$. We think of $\EuT_0$, the central fiber of the Tate curve, as the curve over $\spec \Z$ obtained from $\mathbb{P}^1$ by identifying $p=[1:0]$ and $q=[0:1]$. By definition, $K_0(\EuT_0)$ is the Grothendieck group of the abelian category of vector bundles (locally free sheaves of finite rank) on $\EuT_0$. It can also be thought of as $K_0(\perf \EuT_0)$, the Grothendieck group of the triangulated category of perfect complexes.\footnote{That is, the abelian group generated by the objects, with a relation $[B]=[A]+[C]$ for each distinguished triangle $A\to B \to C \to A[1]$.} The proof of the following lemma is more substantial than one might expect. As partial justification, we point out that $K_0(\EuT_0)$ is an \emph{absolute} invariant of the scheme $\EuT_0$---it is not defined `relative to $\spec\Z$'---and that $\EuT_0$ is \emph{2-dimensional} as a scheme.

\begin{lem}\label{K theory}
The map $(\rank,\det)\colon K_0(\EuT_0) \to \Z \oplus \mathsf{Pic}(\EuT_0)$ is an isomorphism.
Thus a vector bundle on $\EuT_0$ with trivial determinant is stably trivial. 
\end{lem}
\begin{pf}
Let $K_0(R)$ denote the Grothendieck group of finitely-generated projective modules over the commutative ring $R$. We also have the reduced group $\widehat{K}_0(R) = \ker(\rank\colon K_0(R)\to \Z)$ and the group of `stable endomorphisms' $K_1(R)$. For the following standard results in K-theory we refer to the text \cite{WeiK} (see in particular the `Fundamental theorem for $K_1$' (3.6)). We have $\widehat{K}_0(\Z)=0$ since $\Z$ is a PID. The units $R^\times$ are always a subgroup of $K_1(R)$, and  one has $K_1(\Z) = \Z^\times$.
Since $\Z$ is a regular ring, the inclusion-induced maps $\widehat{K}_0(\Z)\to \widehat{K}_0(\Z[t])$ and $K_1(\Z)\to K_1(\Z[t])$ are isomorphisms. One has a split injection $K_1(\Z[t]) \to K_1(\Z[t,t^{-1}])$, induced by the natural map $\Z[t]\to \Z[t,t^{-1}]$, whose cokernel is $K_0(\Z)$. Hence $K_1(\Z[t,t^{-1}])\cong \Z \oplus \Z^\times$. 

There is a group $K_0(\on{\Z[t]}{(t)})$ of complexes of f.g. projective $\Z[t]$-modules whose cohomology is bounded and supported on the ideal $(t)$; and an exact sequence
\[ K_1(\Z[t]) \xrightarrow{}  K_1(\Z[t,t^{-1}]) \xrightarrow{\partial} K_0(\on{\Z[t]}{(t)})  \xrightarrow{} \widehat{K}_0(\Z[t])   \]
\cite[(5.1)]{TT}. From our discussion, we see that the first map has cokernel $\Z$ and that $\widehat{K}_0(\Z[t])=0$, whence $K_0(\on{\Z[t]}{(t)}) \cong \Z$. 

Now consider the `node' $R=\Z[x,y]/(xy)$, and the normalization map $\nu \colon R\to\Z[t] \oplus \Z[t]$, namely, $\nu(f) = (f(x,0),f(0,y))$. One has $\wh{K}_0(\BbK[t] \oplus \BbK[t])=0$. The kernel of $\nu_*\colon \wh{K}_0(R)\to \wh{K}_0(\Z[t] \oplus \Z[t])$ is also zero, because a f.g. projective $R$-module $M$ is determined by $\nu_*M$---which is given by a pair of f.g. projective $\Z[t]$-modules $M_1$ and $M_2$---and `descent data', an isomorphism $\theta\colon M_1/tM_1\to M_2/tM_2$. By stabilizing, we may assume that $M_1$ and $M_2$ are free, and choosing appropriate bases we can make $\theta$ the identity matrix. Hence $\wh{K}_0(R)=0$.

We have (from \cite[(5.1)]{TT} again) a commutative diagram with exact rows
\begin{equation}\label{first K diagram}
\xymatrix{
&   K_1(R) \ar[r]\ar^{\nu_*}[d]  &K_1(R[x^{-1},y^{-1}]) \ar[r]\ar_{\cong}^{\nu_*}[d]  & K_0(\on{R}{(x,y)}) \ar[r]\ar^{\nu_*}[d]  & 0 \\
  0\ar[r]& \bigoplus^2{K_1(\Z[t])} \ar[r]  & \bigoplus^2{K_1(\Z[t,t^{-1}])}\ar[r] & \bigoplus^2{K_0(\on{\Z[t]}{(t)})}  \ar[r] &0 
}\end{equation}
The middle vertical arrow is an isomorphism, since it is induced by a ring isomorphism. Note that $K_1(R)$ contains the units $\Z^\times$. Chasing the diagram, we see that $ K_0(\on{R}{(x,y)}) \cong \Z^2 \oplus S$, where the summand $\Z^2$ mapped isomorphically by $\nu_*$ to $\bigoplus^2{K_0(\on{\Z[t]}{(t)})}$, and $S$ is either $\Z^\times$ or 0.  

We now switch from rings to schemes, referring to \cite{TT} for foundational matters. If $Z$ is a closed subscheme of $X$, $K_0(\on{X}{Z})$ is the Grothendieck group of the  triangulated category of perfect complexes on $\EuT_0$ whose cohomology sheaves are coherent and supported on $Z$. 

The map $\det \colon \wh{K}_0(\EuT_0)\to \Pic(\EuT_0)$ is surjective---apply it to a line bundle. We must prove its injectivity. Let $\nu \colon \mathbb{P}^1\to \EuT_0$ be the normalization map. Let $Z$ be the closure of the image of the nodal section of $\EuT_0$, and $j\colon U\to \EuT_0$ the inclusion of the open complement of $Z$.
There is a commutative diagram with exact rows
\begin{equation}\label{second K diagram}
\xymatrix{ 
K_1(U) \ar^{\partial}[r]\ar^{=}_{\nu^*}[d]  & K_0(\on{\EuT_0}{Z})  \ar^{\mu}[r] \ar_{\nu^*}[d] & \wh{K}_0(\EuT_0) \ar^{j^*}[r]\ar_{\nu^*}[d] & \wh{K}_0(U) \ar^{=}_{\nu^*}[d], \\
K_1(U) \ar^{\partial'}[r] & K_0(\on{\mathbb{P}^1}{\overline{\{ p,q \} } }) \ar^{\mu'}[r] & \wh{K}_0(\mathbb{P}^1) \ar[r] & \wh{K}_0(U) }. 
\end{equation} 
The groups $K_0(\on{X}{Z})$ have an excision property \cite[(3.19)]{TT}, which tells us for instance that restriction induces an isomorphism $K_0(\on{\EuT_0}{Z}) = K_0(\on{V}{Z})$, where $V$ is any open neighborhood of $Z$. Further, it tells us that $K_0(\on{V}{Z}) = K_0(\on{\wh{V}_Z}{\wh{Z}})$, where $\wh{V}_Z$ is the completion of $V$ along $Z$, and $\wh{Z}$ the image of $Z$. Similarly,
\[K_0(\on{\mathbb{P}^1}{\overline{\{ p,q \}}})= K_0(\on{\wh{\mathbb{P}^1}_{\overline{\{p,q} \}}}{ \{\wh{p},\wh{q}\}} ) = \bigoplus^2{K_0(\on{\Z\series{t}}{(t)})}=\bigoplus^2{K_0(\on{\Z[t]}{(t)})}.\] 
Now, $\wh{V}_Z$ is isomorphic to $\spec \Z\series{x,y}/(xy)$, and by naturality of normalization, the map $\nu \colon \wh{V}_Z \to \wh{\mathbb{P}^1}_{ \overline{\{p,q\}} }$ corresponds to the normalization map $\nu\colon \spec \Z\series{x,y}/(xy) \to \spec \bigoplus^2{\Z\series{t}}$, which is itself the completion of $\nu\colon \spec \Z[x,y]/(xy) \to \spec \bigoplus^2{\Z[t]}$. The upshot is that  the map $\nu^* \colon 
K_0(\on{\EuT_0}{Z})  \to K_0(\on{\mathbb{P}^1}{\overline{\{p,q\}}})$ in (\ref{second K diagram}) can be identified with the arrow $\nu_* \colon K_0(\on{R}{(x,y)})  \to \bigoplus^2{K_0(\on{\Z[t]}{(t)})} $ in (\ref{first K diagram}). Hence $\ker \left[K_0(\on{\EuT_0}{Z})  \to K_0(\on{\mathbb{P}^1}{\overline{\{p,q\}}})\right] = S$.

Now take $e \in \wh{K}_0(\EuT_0)$ with $\det e = 1 \in \Pic(\EuT_0)$. Then $\det \nu^*e =1 \in \Pic(\mathbb{P}^1)$, and $\nu^*e = 0 \in \wh{K}_0(\mathbb{P}^1)$ \cite[(4.1)]{TT}. Hence $j^*e=0\in \wh{K}_0(U)$. Thus $e$ is the image $\mu(f)$ of some $f \in K_0(\on{\EuT_0}{Z})$, and $\nu^*f$ maps to zero in $\wh{K}_0(\mathbb{P}^1)$. So $\nu^*f = \partial' g$ for some $g\in K_1(U)$. Then $f-\partial g \in S \subset K_0(\on{\EuT_0}{Z})$, and so $e \in \mu(S)$.

Let $\ell \in \Pic(\EuT_0)$ be the line-bundle obtained from the trivial line bundle on $\mathbb{P}^1$ by gluing the fibers over the $[0:1]$ and $[1:0]$ by the map $(-1)\in \Z^\times$. Let $x= [\ell]-1 \in \wh{K}_0(\EuT_0)$. Then $\det x \neq 1 \in \Pic(\EuT_0)$, so $x \neq 0$. We have $j^*x=0$, so $x= \mu(y)$, say.  Also $\nu^*x=0$; hence $\mu' ( \nu^* y) =0$, and so $\nu^* y = \partial' z$, say, and $y-\partial z\in S$. If $S$ were zero, we would have $x = \mu \circ \partial z =0$. Hence $S$ is not zero, which we have seen implies that $S=\Z^\times$, and therefore $\mu(S) = \{0,x\}$. So our class $e$ from the previous paragraph must be zero.  
\end{pf}

\section{Hochschild cohomology via algebraic geometry}\label{Hochschild section}
In the previous section we set up a two-object, graded linear category $A$, the cohomology $A_C=H^*\EuB_C$ of a dg category associated with an arbitrary Weierstrass curve. We can view (the sum of direct sum of the hom-spaces in) $A$ as a $\BbK$-algebra. In this section we compute its Hochschild cohomology $\HH^\bullet(A,A)$ as a bigraded algebra; or more precisely, the truncated version $\HH^\bullet(A,A)^{\leq 0}$ relevant to non-curved $A_\infty$ deformations \cite{SeiQuartic}. Later (Theorem \ref{bracket}) we shall complete the picture by determining the Gerstenhaber bracket. Our main result is as follows:
\begin{main}\label{HH for cusp curve}
Let $\BbK$ be any field. Introduce the commutative graded $\BbK$-algebra
 \[ T =  \BbK[x,y]/(y^2-x^3, 2y, -3x^2), \] 
concentrated in degree zero. Make it a bigraded algebra  by assigning the following internal degrees $s$ to the generators:
\[   s(x)=2;\quad s(y)=3. \]
Introduce also the free graded-commutative graded algebra 
\[ S^\bullet = \BbK[\beta,\gamma] = \BbK[\beta]\otimes \Lambda[\gamma] ,\quad  \deg\beta=2,\quad\deg \gamma=1, \] 
made a bigraded algebra by assigning internal degrees
\[ s(\beta) = -6; \quad 
s(\gamma)=\begin{cases}
0, & 6\neq 0\in \BbK\\
-2, & 3=0 ;\\
-3, & 2=0.
\end{cases} \]
Define the bigraded algebra $ Q^\bullet  = T\otimes S^\bullet$. Let $Q^{\bullet,\leq 0}$ denote the subalgebra of $Q^\bullet$ spanned by those classes with $s\leq 0$. Then there is a canonical isomorphism of bigraded algebras
\[   \HH^\bullet(A,A)^{\leq 0} \to  Q^{\bullet, \leq 0}. \]
\end{main}
(The bihomogeneous classes of $Q^\bullet$ with $s>0$---i.e., the bihomogeneous classes omitted in $Q^{\bullet, \leq 0}$---are spanned by $x$ and $x \gamma$ when $6\neq 0$; by $x$, $x^2$ and $x^2\gamma$ when $3=0$; and by $x$, $y$, $xy$ and $xy \gamma$ when $2=0$.)
\begin{rmk}
By definition, $\HH^{r+s}(A,A)^s = \ext^r_{(A,A)}(A,A[s])$, the $(r+s)$th derived homomorphism $A\to A[s]$ in the category of graded $(A,A)$-bimodules. Composition of bimodule-homomorphisms yields, on the derived level, an associative  product making $\HH^\bullet(A,A)$ a graded $\BbK$-algebra with an additional internal grading $s$. The product is graded-commutative with respect to the cohomological  grading $\bullet$ \cite{Ger}. The truncation $\HH^\bullet(A,A)^{\leq 0}$, in which $s$ is required to be non-positive, is a subalgebra. Theorem \ref{HH for cusp curve} does \emph{not} make any claims about the untruncated Hochschild algebra.
\end{rmk}

The rank of $\HH^r(A,A)^s$ was computed by the authors in \cite{LP} by a different method. Another approach has been found, in characteristic $0$, by Fisette \cite{Fis}, who also makes the link between $A_\infty$-structures and elliptic curves. Although the additive result is the only part that is essential, we choose to present here this more complete calculation, proved via algebraic geometry on a cuspidal cubic curve $\Ccusp$, by a method explained to us by Paul Seidel.

\begin{lem}
Let $\Ccusp$ be a cuspidal Weierstrass curve over the field $\BbK$. Then $\HH^\bullet(A,A)\cong \HH^\bullet(\Ccusp)$.
\end{lem}
This lemma is restated as part of Prop. \ref{HH for Weierstrass} and proved there. One can compute $\HH^\bullet(\Ccusp)$ using sheaf theory; the heart of the calculation is that of $\HH^\bullet(R,R)$, where $R=\BbK[x,y]/(y^2-x^3)$. The latter calculation can be done---as for any complete intersection singularity---by using a Koszul resolution to replace $R$ by a smooth affine dg manifold, for which a version of the Hochschild--Kostant--Rosenberg theorem is available.

One virtue of the method is that it adapts easily to yield a computation of the Hochschild cohomology of $A$ with a non-trivial $A_\infty$-structure arising from a Weierstrass curve. We shall carry out the computation for the case of a nodal curve, which we will be able to identify with the symplectic cohomology $SH^*$ for the punctured torus.

\subsection{Hochschild cohomology for varieties}\label{HH for var}
\emph{Note}: The ideas that we review in (\ref{HH for var}) have undergone a lengthy evolution, which we have not attempted to trace in detail. Our citations do not necessarily reflect priority.

\subsubsection{The global HKR isomorphism} Let $X$ be a quasi-projective scheme of dimension $d$ over a field $\BbK$. Its Hochschild cohomology $\HH^\bullet(X)$ is defined as
\[  \HH^\bullet(X):= \ext^\bullet_{X\times X}(\delta_*\EuO_X,\delta_*\EuO_X),  \] 
where $\delta\colon X\to X\times X$ is the diagonal map \cite{Swa}.  When $X$ is smooth, there is a natural morphism of complexes of sheaves
\[  (\mathbb{L}\delta^*)(\delta_*\EuO_X) \to \bigoplus_{r}\Omega^r_X[r]  \]
(trivial differential in the complex on the right), the Hochschild--Kostant--Rosenberg (HKR) map; see \cite{Yek} or the pr\'ecis in \cite{Cal}.  When $d!$ is invertible in $\BbK$---which means that $\ch(\BbK)$ is  either $0$ or $>d$---the HKR map is a quasi-isomorphism.\footnote{For a general field $\BbK$, this map may not be a quasi-isomorphism but it is nevertheless true (\cite[Ex. 9.1.3]{Wei},) that $\HH^\bullet(X,X) \cong \Lambda^\bullet T_X$ when $X=\mathbb{A}^n$.} From the HKR map and the adjunction $\RHom^\bullet_{X\times X}(\delta_*\EuO_X,\delta_*\EuO_X)\cong \RHom^\bullet_X((\mathbb{L}\delta^*)(\delta_*\EuO_X),\EuO_X)$, one obtains, under this assumption on $\BbK$, an isomorphism of graded $\BbK$-vector spaces
\[ \mathsf{HKR}^n \colon \HH^n(X) \to \bigoplus_{p+q=n}{\h^p({\Lambda^q T_X})}.  \]
\paragraph{The Hodge spectral sequence.} For $X$ quasi-projective over $\BbK$, there is a sheaf of graded algebras $\mathcal{HH}^*$ on $X$, the sheafification of a natural presheaf whose sections are $\Gamma(U) =\HH^*(U)$. There is by \cite{Swa} a local-to-global (or `Hodge') spectral sequence $E^{**}_*$ converging to $\HH^*(X)$, with
\begin{equation} \label{Hodge}
E_2^{pq} \cong \h^p(X, \mathcal{HH}^q).
\end{equation}
When $X$ is smooth of dimension $d$, and $d!$ is invertible in $\BbK$, one has $\mathcal{HH}^q \cong \Lambda^q T_X$ by the HKR isomorphism. Comparing dimensions of $\HH^n(X)$ and $E_2^n$, one sees that the spectral sequence degenerates at $E_2$. One therefore has
\[ E_\infty^{pq} = E_2^{pq} \cong  \h^p(X, \Lambda^qT_X).  \]
As a sheaf of graded algebras, $\mathcal{HH}$ has cohomology spaces which form a bigraded algebra  $\h^*(X, \mathcal{HH}^\bullet)$. Like any local-to-global spectral sequence computing self-Exts, the Hodge spectral sequence is multiplicative. From this one sees:
\begin{lem}\label{edge}
For $X$ a quasi-projective $\BbK$-scheme, the edge-map $\HH^\bullet(X) \to E_2^{0,\bullet} =  \h^0(X; \mathcal{HH}^\bullet)$ is an algebra homomorphism. 
\end{lem}

\subsubsection{Derived categories} Let $\mathsf{QC}(X)$ denote a dg enhancement for the unbounded derived category of quasi-coherent sheaves on the $\BbK$-scheme $X$ (e.g. \cite{LO,Toe}). Thus one has an equivalence of $H^0(\mathsf{QC}(X))$ with the unbounded derived category.  There is a natural isomorphism of rings \cite{Toe}
\[   \HH^*(X) \to \HH^*(\mathsf{QC}(X)), \]
where by the latter we mean the Hochschild cohomology of the dg category $\mathsf{QC}$ (i.e., self-Ext of the identity functor). Let $\widetilde{\perf}(X)$ be the full dg subcategory of $\mathsf{QC}(X)$ of objects which map to perfect complexes in the unbounded derived category. Since $\widetilde{\perf}(X)$ is a full subcategory of $\mathsf{QC}(X)$, there is a restriction (ring) map
\[  \HH^*(\mathsf{QC}(X)) \to \HH^*(\widetilde{\perf}(X)),   \]
and this too is an isomorphism when $X$ is quasi-projective, because $\mathsf{QC}(X)$ is the \emph{ind-completion}
of $\widetilde{\perf}(X)$ \cite{BFN}. 
When $X$ is quasi-projective, we also have $\widetilde{\perf}(X)\simeq \tw \VB(X)$. Putting these facts together, we obtain the following:
\begin{lem}\label{perf computes Hochschild}
For any quasi-projective $\BbK$-scheme $X$, one has a canonical algebra-isomorphism 
\[  \HH^\bullet(\tw \VB(X)) \cong \HH^\bullet(X). \]
\end{lem}
If $T$ is a split-generator for $\tw\VB (X)$, the Morita invariance of Hochschild cohomology (e.g. \cite{Toe}) implies that the restriction map
\[ \HH^\bullet(\tw\VB (X)) \to \HH^\bullet(\End(T))  \]
is an isomorphism. Hence we have $\HH^\bullet(\End(T))\cong \HH^\bullet(X)$, and consequently
\begin{prop}\label{HH for Weierstrass}
Take an abstract Weierstrass curve $(C,\sigma,\omega)$ and let $T=\EuO\oplus \EuO_\sigma$. Then there is a canonical algebra-isomorphism
\[ \HH^\bullet(C)\to \HH^\bullet( \End(T)).  \]
In particular, by Lemma \ref{formality},
\[ \HH^\bullet(\Ccusp) \cong \HH^\bullet(A,A).  \]
\end{prop}
Note that the latter isomorphism respects bigradings. The internal grading of $\HH^\bullet(A,A)$ comes from the grading of $A$.  The internal grading of $\HH^\bullet(\Ccusp)$ arises because of the $\BbK^\times$-action on $\Ccusp$ discussed in the proof of (\ref{formality}), which gives rise to an action on $\tw\VB (\Ccusp)$, hence on Hochschild cohomology.  They agree because the action on $\Ccusp$ induces the grading of $\End(T)$ (cf. the proof of (\ref{formality})).

\subsubsection{Curves}
We now specialize to quasi-projective curves $C$ over $\BbK$. HKR reads as follows:
\begin{lem}\label{HKR curves}
For a non-singular curve $C$ over an arbitrary field $\BbK$, one has an HKR isomorphism of $\BbK$-modules
\[ \mathsf{HKR}^n \colon \HH^n(C) \to \h^n(\EuO)\oplus \h^{n-1}(TC).  \]
\end{lem}
For more general curves, we have the following
\begin{lem}\label{curve ses}
For any quasi-projective curve $C$, one has short exact sequences
\[ 0 \to \h^1(C, \mathcal{HH}^{n-1}) \to  \HH^n(C)  \to \h^0(C, \mathcal{HH}^n) \to 0.\] 
\end{lem}
\begin{pf}
Since curves have cohomological dimension 1, the Hodge spectral sequence (\ref{Hodge}) is supported at $E_2$ in two adjacent columns $p \in \{ 0,1\}$, and therefore degenerates. The edge-maps then give rise to these short exact sequences.
\end{pf}
One has $\mathcal{HH}^0 = \EuO_C$, while $\mathcal{HH}^1\cong \deriv(\EuO_C)$, the sheaf of derivations of $\EuO_C$, also known as the tangent sheaf $\EuT_C$. 

To go further, we suppose that $C$ is given with a point $s\in C(\BbK)$, and that $C\setminus \{s\}$ is non-singular. We take an affine open  cover $\EuU=\{U, V\}$ of $C$ such that $s\not \in V$.  We compute sheaf cohomology $\h^*(C,\mathcal{HH}^\bullet)$ using the \v{C}ech complex $ \check{C}^*=\check{C}^*(\EuU, \mathcal{HH}^q)$. To validate such a computation, we must check that $\h^i(\mathcal{HH}^q|_Y)=0$ when $Y=U$, $V$ or $U\cap V$ and $i>0$. As coherent sheaves, $\mathcal{HH}^0$ and $\mathcal{HH}^1$ have no higher cohomology on affine open subsets. When $q>1$, $\mathcal{HH}^q$ is supported at $s$, by (\ref{HKR curves}). It therefore has no higher cohomology. We now proceed with the computation. 

Since $V$ is non-singular, (\ref{HKR curves}) implies that $\HH^0(V)\cong\Gamma(V,\mathcal{O}_V)$, $\HH^1(V)\cong \Gamma(V, TV)$ and $\HH^q(V)=0$ for $q>1$. The same applies over $U\cap V$. So the \v{C}ech complex is
\[ 0 \to \HH^q(U) \oplus \HH^q(V) \to \HH^q(U\cap V) \to 0 \]
supported in degrees 0 and 1. Hence:

\begin{lem}\label{sheaf coh of HH}
For $q>1$, one has
\[ \h^0(C,\mathcal{HH}^q) \cong \HH^q(U) , \quad  
   \h^1(C,\mathcal{HH}^q)=0. \]
Therefore, for $n>2$,
\[  \HH^n(C)  \xrightarrow{\cong} \h^0(C, \mathcal{HH}^n) \cong \HH^n(U) = 
\HH^n(\Gamma(\EuO_U),\Gamma(\EuO_U)).  \]
\end{lem}
In low degrees, we have an isomorphism $\HH^0(C) \cong \h^0(\EuO_C)$ and two short exact sequences
\begin{align}\label{two ses}
& 0\to \h^1(C, \EuO_C) \to \HH^1(C) \to  \h^0(C,\EuT_C)  \to 0\\
& 0\to \h^1(C,\EuT_C)\to \HH^2(C)\to \HH^2(U)\to 0 .
 \end{align}
In summary,
\begin{leftbar}
\emph{$\HH^\bullet(C)$ coincides with $\HH^\bullet(\Gamma(\EuO_U),\Gamma(\EuO_U))$ with corrections in degrees $\leq 2$ from the cohomology of the functions and of the vector fields on $C$.}
\end{leftbar}

\subsubsection{Plane curve singularities}
If  $R=\BbK[x,y]/(f)$ is a plane curve singularity, over a field $\BbK$, one can compute $\HH^\bullet(R,R)$ via Koszul resolutions and a version of HKR. We have been informed that such calculations go back to Quillen \cite{Qui}, but an explicit recipe, valid for complete intersections, is explained by Kontsevich in \cite{KonApp}. The result is that
\begin{equation}\label{kontsevich}  
\HH^\bullet(R,R) \cong H^\bullet(D, d_D), 
\end{equation}
where $(D, d_D)$ is a certain dga, namely, the supercommutative $\BbK$-algebra
\[ R \otimes \BbK[\beta, x^*, y^*], \quad \deg \beta=2,\quad \deg x^*=\deg y^* = 1 \]
and $d_D(R)=0$, $d_D(\beta)=0$, $d_D x^* = f_x \beta$, $d_D y^* = f_y \beta$. Thus 
\[  D_{2n} = R \beta^n \oplus R \beta^{n-1}x^*y^*, \quad D_{2n+1} = R \beta^n x^* \oplus R \beta ^n y^*.\]

With the isomorphism (\ref{kontsevich}) understood, one immediately reads off the even Hochschild cohomology of $R$.
Let $T= R/(f_x, f_y)$. In the literature $T$ is often called the \emph{Tjurina algebra} of this isolated hypersurface singularity; it parameterizes a miniversal deformation of the singularity; see for instance \cite[Theorem 14.1]{HarDef} for the case of curve singularities. 
\begin{lem}\label{tjurina}
Assume that $f_x$ and $f_y$ are not both zero.  Then the map
\[ T \to \HH^{2n}(R,R) , \quad [\alpha] \mapsto [\alpha \beta^n]  \]
is an isomorphism for each $n>0$. Hence $\bigoplus_{n>0}{\HH^{2n}(R,R)}\cong T\otimes \beta \BbK[\beta] $. 
\end{lem}
Let $M$ be the $R$-submodule of $R\oplus R$ of pairs $(\alpha_1,\alpha_2)$ with $\alpha_1 f_x + \alpha_2 f_y =0 \in R$.
It has a submodule $N$ generated by $(f_y, -f_x)$. There is a skew-symmetric pairing
\[   \omega\colon M\otimes M\to T, \quad (\alpha_1,\alpha_2; \gamma_1,\gamma_2)\mapsto
[\alpha_1 \gamma_2-\alpha_2 \gamma_1]  \]
such that $\omega(N\otimes M)=0$. When $f$ is irreducible and $df\neq 0$, one has $\omega=0$. Indeed, $R$ is then an integral domain, and the matrix 
\[ \left[ \begin{array}{cc}  \alpha_1 & \alpha_2\\ \gamma_1 & \gamma_2\end{array}\right] , \]
over $R$, has a non-trivial kernel and therefore vanishing determinant.

Observe also that the $R$-module $M/N$ is actually a $T$-module.
\begin{lem}
One has a surjective map
\[  M \to  \HH^{2n+1}(R,R),\quad (\alpha_1,\alpha_2) \mapsto \beta^n [\alpha_1 x^* + \alpha_2 y^* ] \]
The kernel is $0$ when $n=0$ and is $N$ when $n>0$.
\end{lem}
The proof is an easy check, in light of (\ref{kontsevich}). Note that $\HH^1(R,R)=\deriv R$, the derivations of $R$; the isomorphism $M = H^1(D,d_D) \to \deriv(R)$ induced by (\ref{kontsevich}) is the map $(\alpha_1,\alpha_2)\mapsto \alpha_1\partial_x + \alpha_2 \partial_y$. We deduce:

\begin{prop}\label{HH for R}
There is a canonical map of algebras
\[ \HH^\bullet(R,R) \to  Q^\bullet, \]
which is an isomorphism in degrees $\geq 2$.  Here 
\[ Q^\bullet = \BbK[\beta] \otimes \left(T \oplus \frac{M}{N}[-1]\right),\quad
\deg T=0,\quad \deg\beta=2.  \]
The product in $Q^\bullet$ combines the algebra structure of $T$, the left and right $T$-module structures of $M/N$, and the skew pairing $\omega$:
\[ (\beta^n \otimes (t + m) ) \cdot (\beta^{n'}\otimes  (t'+m'))
= \beta^{n+n'} \otimes (t\cdot t' + t\cdot m' + m \cdot t') + \beta^{n+n' +1} \otimes \omega(m,m').  \]
\end{prop}

The following observation, whose proof is immediate, is helpful in computing $M/N$:
\begin{lem}\label{koszul}
$M/N$ is the middle homology $H_1( f_x,f_y; R)$ of the Koszul complex $K(f_x,f_y)$ , i.e., the chain complex
\[  
R \xrightarrow{\left[\begin{array}{c}f_y \\ -f_x\end{array}\right]} R^2 \xrightarrow{\left[\begin{array}{cc} f_x & f_y \end{array}\right]}  R . \]
\end{lem}
$K(f_x,f_y)$ (notation from Serre \cite[ch. IV]{Ser}) is the tensor product $K(f_x)\otimes_R K(f_y)$ of Koszul complexes for the elements $f_x$ and $f_y$. Here $K(a)$ denotes the complex  $ 0\to R \stackrel{a}{\to} R \to 0 $ where the map---multiplication by $a$---maps degree 1 to degree 0.

\subsection{Hochschild cohomology for the cuspidal cubic}
\subsubsection{Hochschild cohomology for the affine curve}
We take $ f(x,y) = y^2-x^3.$ The Tjurina algebra $T$ (\ref{tjurina}) is as follows:
\[ T \cong  \begin{cases}
\BbK[x]/(x^2) & \text{if }6\neq 0,\\
\BbK[x]/(x^3) &\text{if } 3=0,\\
\BbK[x,y]/ (x^2,y^2) & \text{if }2=0.
\end{cases}\]
(The isomorphism takes $x\in T$ to $x$, and takes $y \in T$ to $0$ in the first two cases and to $y$ in the third.)
\begin{lem}\label{M mod N}
There are isomorphisms of $R$-modules 
\[ M/N  \cong  
\begin{cases}
\Ann_{R/f_y}(f_x) =(x) 		& \text{if }6 \neq 0,\\
\Ann_{R/f_y}(f_x) = R/f_y	& \text{if }3=0,\\
\Ann_{R/f_x}(f_y) = R/f_x 	& \text{if }2=0
\end{cases}
\]
and 
\[ T\to M/N,\quad t \mapsto  
\begin{cases}
 [2xt,3yt] 	&  	\text{if }6 \neq 0,\\
 [t,0] 		&  	\text{if }3 = 0,\\
 [0, t]  		&	\text{if }2 =0.
\end{cases} \]
\end{lem}
\begin{pf}
By Lemma \ref{koszul}, $M/N$ is equal as an $R$-module to $H_1(f_x,f_y;R)$, and hence isomorphic to $H_1(K(f_x;R)\otimes K(f_y;R))$. 

If $2\neq 0 \in \BbK$ then the element $f_y=2y$ is not a zero-divisor in $R$. 
Therefore the complex $K(f_y)$ has homology only in degree zero, and $H_0(K(f_y))=R/(f_y)$.  Projection $R \to R/(f_y)$ defines a quasi-isomorphism $K(f_y) \to R/(f_y)$. Hence $K(f_x,f_y; R)\simeq K(f_x;R)\otimes R/f_y$, and the latter complex is 
\[ 0\to R /(f_y) \stackrel{f_x}{\to} R/(f_y)\to 0,  \] 
where the differential maps degree 1 to degree 0. So $H_1( f_x,f_y; R)$ is isomorphic as an $R$-module (hence also as a $T$-module) to the annihilator of $f_x$ in $R/f_y$. Explicitly, the annihilator is $(x)\subset \BbK[x,y]/(x^3,y)$ if $3\neq 0$ and $\BbK[x,y]/(x^3,y)$ if $3 =0$. If $2=0$ then the element $f_x=-3x^2$ is not a zero-divisor in $R$. We can then run the same argument with the roles of $f_x$ and $f_y$ interchanged, with a similar outcome. Here $R/f_x$ is $\BbK[x,y]/(x^2,y^2)$, and $f_y=0$. 

At this stage, we can easily see that $\dim_\BbK T = \dim_\BbK (M/N)$. The given map $T\to M/N$ is injective, by another easy check, hence an isomorphism.
\end{pf}
Since $f$ is irreducible and $df\neq 0$, we have the
\begin{lem}\label{skew}
The skew pairing $\omega \colon (M/N) \times (M/N) \to T$ is zero.
\end{lem}

Collating results (Prop. \ref{HH for R} and Lemmas \ref{M mod N}, \ref{skew})  we obtain
\begin{lem}\label{HH for cusp singularity}
Consider the free graded-commutative algebra $Q^\bullet = T\otimes \BbK[\beta,\gamma]$, where $\deg \beta=2$ and $\deg \gamma=1$.  One then has a map of graded $\BbK$-algebras 
\[ \HH^\bullet(R,R)\to Q^\bullet \]
which is an isomorphism except in degrees 0 and 1.
\end{lem}
\begin{pf}
In light of Prop. \ref{HH for R}, we need only see that $T\otimes \BbK[\beta,\gamma]$ agrees with the algebra named $Q^\bullet$ there. In view of the isomorphism $T\cong M/N$, and the vanishing of $\omega$, the latter algebra is
\[  \BbK[\beta]\otimes (T\oplus T[-1]) \cong \BbK[\beta,\gamma]\otimes T,\quad \deg \gamma=1.  \]
\end{pf}
We can be more explicit. We know that $\HH^\bullet(R,R)\cong H^\bullet(D)$. The map 
$H^\bullet(D)\to Q^\bullet$ is $R$-linear.  It maps $\beta$ to $\beta$. It maps a certain class $[ax^* +b y^*] \in H^1(D)$ to $\gamma$.  The coefficients are given by $(a,b)=(2x,3y)$ when $6 \neq 0$; $(a,b) = (1,0)$ when $3=0$; and $(a,b)=(0,1)$ when $2=0$. 

\subsubsection{Global calculation}
\begin{lem}\label{derivations}
The Lie algebra $\h^0(\Ccusp,\EuT_{\Ccusp})$ contains linearly independent vector fields $v_0$, $v_1$ which restrict to $U$ as the vector fields
\[  v_0 |_U = 2x\partial_x+ 3y\partial_y,  \quad v_1|_U =  2y \partial_x + 3x^2 \partial_y.  \]
When $6\neq 0$, $v_0$ and $v_1$ span. When $3=0$, $\h^0(\Ccusp,\EuT_{\Ccusp})$ is 3-dimensional, spanned by $v_0$ and $v_1$ together with
\[ v_{-2} = -\partial_x. \]
When $2=0$, it is 4-dimensional, spanned by $v_0$, $v_1$ and
\[ v_{-1} = x \partial_y, \quad v_{-3} = \partial_y. \]
\end{lem}
(We will see presently that the subscript $s$ of $v_s$ represents the internal degree.)
\begin{pf}
As shown to us by Seidel, the calculation of $\h^0$ is straightforward when one takes the `abstract' view of $\Ccusp$ (Lemma \ref{abstract cusp}), as $\mathbb{P}^1$ with the marked point $c$, and non-standard structure sheaf $\EuO^c$. Thus $\EuT_{\Ccusp}= \deriv(\EuO^c)$ is the sheaf of meromorphic vector  fields $\theta$ on $\mathbb{P}^1$, with a pole only at $c$, which preserve $\EuO^c$. The latter condition forces $\theta$ to be regular at $c$ except when $6=0$. When $\theta$ is regular at $c$, it must vanish there, except in characteristic 2. One has  $\h^0(T_{\mathbb{P}^1}) = \mathfrak{pgl}_2$ (the Lie algebra of the automorphism group of $\mathbb{P}^1$), spanned by $v_{-1}= \partial_z$, $v_0= z \partial_z$ and $v_1= z^2 \partial_z$. In characteristic 2, all three lie in  $\h^0(\deriv(\EuO^c))$, while when $2\neq 0$, only $v_0$ and $v_1$ do (they span the Borel subalgebra $\mathfrak{b}\subset \mathfrak{pgl}_2$ of  upper triangular $2\times 2$ matrices modulo scalars).

When $3=0$, there is a additional derivation $v_{-2}= z^{-1}\partial_z$; when $2=0$, there is again one additional derivation,  $v_{-3} = z^{-2}\partial_z$.

To interpret these vector fields as derivations of $R$, we observe that the normalization map $\nu\colon \mathbb{P}^1\to \Ccusp$ corresponds to the map of rings $R=\BbK[x,y]/(y^2-x^3)\to \BbK[z]$, $x\mapsto z^2$, $y\mapsto z^3$. Using this we compute that the $v_i$ restrict to $U$ in the way stated above.
\end{pf}
\begin{lem} \label{H1 vanishes}
One has $\h^1(\Ccusp, \EuT_{\Ccusp})=0$.
\end{lem}

\begin{pf}
Consider $\mathbb{P}^1$ with the cusp modifying its structure sheaf at $c=[0:1]$. Let $b=[1:0]$. We use the affine coordinate $[z:1]\mapsto z$. Take a derivation $\theta$ of $\EuO_{\mathbb{P}^1}$ (i.e., a vector field) over $\mathbb{P}^1\setminus \{b,c \}$. We must show that it is  the difference $u-v$ of vector fields $u$ on $U=\mathbb{P}^1\setminus \{b\}$ and $v$ on $V=\mathbb{P}^1\setminus \{c\}$. We can extend $\theta$ to a meromorphic vector field on $\mathbb{P}^1$. We proceed by induction on the order $d$ of the pole of $\theta$ at $b$. If $\theta$ is regular at $b$ (i.e., $d\leq 0$) then we take $u=0$ and $v= - \theta$. For the inductive step, say $\theta\sim az^d$ at $b$, where $a\neq 0$ and $d>0$. Let $u=az^d$. Then $u$ defines a derivation of $\EuO^c$ near $c$, because $d>0$. Moreover, $\theta': =u -\theta$ has a pole of order $<d$ at $b$, so by induction we can write $\theta' = u'-v'$ for $u'$ on $U$ and $v'$ on $V$; then $\theta=(u-u') +v' $ and we are done.
\end{pf}

At this stage it becomes useful to bring internal gradings into play. The Hochschild cohomology $\HH^\bullet(\Ccusp)$ carries an internal grading $s$, arising from the $\BbK^\times$-action on $\Ccusp$. Equivalently, this is the internal grading of $\HH^\bullet(A,A)$ arising from the grading of $A$. There is also an action on $\mathcal{HH}^\bullet$, hence an internal grading on $\h^0(\Ccusp,\mathcal{HH}^\bullet)$.

Under the $\BbK^\times$-action on $\Ccusp$, the functions $X$, $Y$ and $Z$ on $\Ccusp$ have respective weights $w(X)=2$, $w(Y)=3$ and $w(Z)=0$. Hence on the affine part, the functions $x=X/Z$ and $y=Y/Z$ have weights $w(x)=2$ and $w(y)=3$. The dga $(D^*,d_D)$ inherits a $\BbK^\times$-action, i.e., a grading $s$, in which the weights of $x$ and $y$ are $s(x)=2$ and $s(y)=3$. The variables $x^*$ and $y^*$ have weights $ s(x^*)=-2$ and $s(y^*)=-3$. It then follows that $s(\beta)=-6$. 

In Lemma \ref{derivations}, one should understand $\partial_x$ to have weight $s=-2$ and $\partial_y$ weight $s=-3$. Then the vector field $v_k$ has weight $k$.

Introduce the internal grading on $Q^\bullet$ as in the statement of Theorem \ref{HH for cusp curve}. It is is set up so that the homomorphism $\HH^\bullet(R,R)\to Q^\bullet$ respects it. Let $Q^{\bullet,\leq 0}$ be the subalgebra of $Q^\bullet$ where the internal grading is non-positive. Similarly, define $\h^0(\Ccusp,\mathcal{HH}^\bullet)^{\leq 0}$.

\begin{lem}\label{HH to H0}
The map $\HH^{\bullet,\leq 0}(\Ccusp)\to \h^0(\Ccusp, \mathcal{HH}^{\bullet})^{\leq 0}$ is an isomorphism.
\end{lem}
\begin{pf}
The kernel of the surjective map $\HH^\bullet(\Ccusp)\to \h^0(\Ccusp, \mathcal{HH}^\bullet)$ is spanned by $\h^1(\EuT_{\Ccusp})$ in degree 2---this vanishes by Lemma \ref{H1 vanishes}---and $\h^1(\EuO)$ in degree $1$. The action of $\BbK^{\times}$ has weight $1$ on $\h^1(\EuO)$. Hence the restricted map
$\HH^\bullet(\Ccusp)^{\leq 0}\to \h^0(\Ccusp, \mathcal{HH}^\bullet)^{\leq 0}$ is injective, and so an isomorphism. 
\end{pf}

\begin{lem}
One has an isomorphism of graded algebras $ \h^0(\Ccusp,\mathcal{HH}^\bullet)^{\leq 0} \to Q^{\bullet,\leq 0}. $
\end{lem}
\begin{pf}
There is a restriction map $r\colon \h^0(\Ccusp,\mathcal{HH}^\bullet)\to \HH^\bullet(U)=\HH^\bullet(R,R)$. By Lemmas \ref{sheaf coh of HH} and \ref{H1 vanishes}, $r$ is an isomorphism in degrees $\geq 2$. 
In degree 1, it is the restriction map $\h^0(\Ccusp,\EuT_{\Ccusp}) \to \h^0(U,\EuT_U)$, and therefore again injective. The composite
\[  \h^0(\Ccusp,\mathcal{HH}^\bullet) \xrightarrow{r}\HH^\bullet(U)=\HH^\bullet(R,R) \to Q^\bullet \]
of $r$ with the map from Lemma \ref{HH for cusp singularity}, restricted so as to map between the non-positively graded subalgebras, is the sought-for map. It is certainly an isomorphism in degrees $\bullet > 1$.   It is also an isomorphism in degree $0$, where both sides reduce to $\BbK$. 

The non-positively graded part of $\h^0(\Ccusp,\deriv(\EuO_{\Ccusp})^{\leq 0})$ has dimension $\dim T -1$, 
by Lemma \ref{derivations}: it is spanned by $v_0$, together with $v_{-2}$ in characteristic 3, and $v_{-1}$, $v_{-3}$ in characteristic 2. One also has $Q^1=T$, and $\dim (Q^1)^{\leq 0} = \dim T-1$.

We have identified $v_j$ as derivations of $R$, hence as elements of $H^1(D,d_D)=M$. One has $Q^1\cong M/N$, in such a way that the map $H^1(D,d_D)\to Q^1$ corresponds to the quotient map $\h^1(D,d_D)=M\to M/N$. To show that $\h^0(\Ccusp,\deriv(\EuO_{\Ccusp})^{\leq 0})\to (Q^1)^{\leq 0}$ is injective, one need only show that the relevant elements $v_j$ are linearly independent in $M/N$. But as elements of $M$, one has $v_0=[2x, 3y]$; $v_{-2} = [-1,0]$; $v_{-1}=[0,x]$; and $v_{-3}=[0,1]$. Linear independence in $M/N$ is easily seen.
\end{pf}

\begin{pf}[Proof of Theorem \ref{HH for cusp curve}]
The isomorphism we want is the composite of the isomorphisms described in the last two lemmas:
\[ \HH^\bullet(\Ccusp)^{\leq 0}\to \h^0(\Ccusp,\mathcal{HH}^{\bullet})^{\leq 0}\to Q^{\bullet, \leq 0}. \]
\end{pf}

Before leaving $\Ccusp$, it will be worthwhile to spell out our findings about $\HH^1(\Ccusp)^{\leq 0}$ in a clean form. They are as follows:
\begin{prop}\label{HH1}
\begin{enumerate}
\item
The canonical map 
\[ \HH^1(\Ccusp)^{\leq 0}  \to \h^0(\Ccusp, \EuT_{\Ccusp})^{\leq 0}  \]
is an isomorphism. 
\item
Let $\ker d \subset W$ be as at (\ref{ker d}); it is the Lie subalgebra of $\mathfrak{pgl}_3(\BbK)$ of those vector fields on $\mathbb{P}^2(\BbK)$ which preserve $\Ccusp$. Thus there is a canonical map 
\[ \ker d \to  \h^0(\Ccusp, \EuT_{\Ccusp})^{\leq 0},\] 
natural in $\BbK$. The latter map is an isomorphism.
\end{enumerate}
\end{prop}
\begin{pf}
Only the second clause has not already been proved. As a vector field on $\mathbb{P}^2$, we have $\partial_u = (\partial x/\partial u) \partial_x+ (\partial y/\partial u) \partial_y = -2 x\partial_x - 3y\partial_y$. Similarly, $\partial_s = -x\partial_y$; $\partial_r = - \partial_x$; and $\partial_t = -\partial_y$. Now, $\partial_u$ is tangent to $\Ccusp$, and the restriction of $\partial_u$ to $\Ccusp$ is the vector field $v_0$. When $2=0$, $\partial_s$ and $\partial_t$ are also tangent to $\Ccusp$, and they restrict to $\Ccusp$ as the respective vector fields $v_{-1}$ and $v_{-3}$. When $3=0$, $\partial_r$ is tangent to $\Ccusp$ and restricts to $\Ccusp$ as $v_{-2}$. Thus, by Lemma \ref{derivations}, the map $\ker d\to  \h^0(\Ccusp, \EuT_{\Ccusp})^{\leq 0}$ is an isomorphism in every case.\end{pf}

\subsection{Hochschild cohomology of the nodal cubic curve}

\emph{Note: this subsection is not used elsewhere in the paper.} 

Let $\EuT_0 \to \spec \Z$ be the central fiber of the Tate curve, defined by $w(x,y)=0$ where $w(x,y) = y^2+xy - x^3$: the proposed mirror to the punctured 2-torus $T_0$. The Hochschild cohomology $\HH^\bullet(\EuT_0)$ is also the Hochschild cohomology of $\tw \VB(\EuT_0)$. Since Hochschild cohomology is invariant under $A_\infty$ quasi-equivalences, Theorem \ref{mainth} says that $\HH^\bullet(\EuT_0)\cong\HH^\bullet(\EuF(T_0)^{\exact})$. So, taking that theorem for granted for now, one can regard this subsection as a computation of $\HH^\bullet(\EuF(T_0)^{\exact})$. In general, there is a natural map to the Hochschild cohomology of the exact Fukaya category from the \emph{symplectic cohomology} \cite{SeiBias} of the manifolds, and in certain cases \cite{SeiDef} this is expected to be an isomorphism. The graded ring we compute here is indeed isomorphic to $SH^\bullet(T_0)$, though we do not check that the map is an isomorphism.

\begin{thm}\label{HH nodal}
Over any field $\BbK$, there is an isomorphism of graded algebras
\[ \HH^\bullet(\EuT_0) \to 
\frac{\BbK[\beta,\gamma_1,\gamma_2]}{\left(\gamma_1\gamma_2,\beta( \gamma_1 - \gamma_2)\right)}  \]
where
\[ \deg \gamma_1 = \deg\gamma_2=1, \quad \deg \beta = 2,\]
and $\BbK[\beta,\gamma_1,\gamma_2]$ is a free supercommutative algebra.
\end{thm}
The computation follows similar lines to that for the cuspidal cubic, and we shall be terse. Taking $R= \BbK[x,y]/(w)$, one readily checks the following:
\begin{lem}\label{Tjurina nodal}
The Tjurina algebra $T$ is reduced to $\BbK$. Its module $M/N$ is 1-dimensional, and therefore the skew pairing $\omega\colon M/N\otimes M/N\to \BbK$ is zero. Hence there is a surjective map of algebras
\[ \HH^\bullet(R,R) \to \BbK[\gamma,\beta],\quad \deg \beta=2,\quad \deg \gamma=1, \] 
with kernel concentrated in degree 1. 
\end{lem}

The normalization of $\EuT_0$ is isomorphic to $\mathbb{P}^1$; the normalization map $\nu\colon \EuT_0 \to \mathbb{P}^1$ carries two points, $p$ and $q$, say, to the node. The pullback $\nu^* \EuO_{\EuT_0}$ is contained in $\EuO_{\mathbb{P}^1}$ as the sheaf of functions $f$ such that $f(p)=f(q)$. Similarly, 
$\nu^* \deriv(\EuO_{\EuT_0})$ is contained in $\deriv(\EuO_{\mathbb{P}^1})$ as the sheaf of vector fields $v$ with $v(p)=v(q)=0$. Hence 
\[ h^0(\EuT_0,\deriv(\EuO_{\EuT_0})) = h^0(\mathbb{P}^1,\deriv( \nu^*\EuO_{\EuT_0} )) = 1.  \]
By Riemann--Roch for $\mathbb{P}^1$, 
\[ h^1(\mathbb{P}^1,\deriv( \nu^*\EuO_{\EuT_0} )) = 0.  \]
In \v{C}ech terms, this means that if we cover $\mathbb{P}^1$ by $U=\mathbb{P}^1\setminus \{p,q\}$ and $V=\mathbb{P}^1\setminus \{r\}$, then every vector field on $U\cap V$ can be expressed as the difference $u-v$ of vector fields on $U$ and $V$ with descend to $\nu(U)$ and $\nu(V)$ respectively. Hence 
\[  \h^1(\EuT_0, \deriv(\EuO_{\EuT_0}))=0.  \]
This vanishing result gives us the following lemma:
\begin{lem}\label{surj}
The map $\HH^{\bullet}(\EuT_0) \to \HH^0(\EuT_0, \mathcal{HH}^\bullet)$ is surjective with kernel $\h^1(\EuO)[-1]$.
\end{lem}
We know that $\h^0(\EuT_0, \mathcal{HH}^\bullet) = \HH^\bullet(R,R)$ in degrees $\bullet \geq 2$. We have
\[  \h^0(\EuT_0, \mathcal{HH}^1) = \h^0(\EuT_0, \deriv(\EuO)) \cong \BbK.   \]
\begin{lem}
The composite of the maps
\[  \h^0(\EuT_0, \mathcal{HH}^1) \to \HH^1(R,R) \to  \BbK \{\gamma\} , \]
(the restriction map and the surjection from Lemma \ref{surj}) is an isomorphism.
\end{lem}
\begin{pf}
We know that $\h^0(\EuT_0, \mathcal{HH}^1)$, the space of global vector fields, has dimension 1, so we need only show that the composite map is non-zero. The composite 
$\h^0(\EuT_0, \mathcal{HH}^\bullet) \to \HH^\bullet(R,R) \to  \BbK [\gamma,\beta]$ is a map of algebras. Moreover,  by Lemma \ref{Tjurina nodal}, $\h^0(\EuT_0, \mathcal{HH}^\bullet)$ is generated in degrees 1 and 2. If the map on $\h^0(\EuT_0, \mathcal{HH}^1)$ were zero, the same would be true of the map $\h^0(\EuT_0, \mathcal{HH}^3)\to \BbK {\beta\gamma}$; yet we know that the latter map is an isomorphism. 
\end{pf}

\begin{pf}[Proof of Theorem \ref{HH nodal}]
The previous lemma implies that
\[ \h^0(\EuT_0,\mathcal{HH}^\bullet) \to  \BbK[\beta,\gamma] \]
is an isomorphism. 
Now consider the surjective map 
\[ \HH^\bullet(\EuT_0) \to \h^0(\EuT_0,\mathcal{HH}^\bullet). \]
Its kernel is $\h^1(\EuT_0,\deriv(\EuO))[-2] \oplus \h^1(\EuT_0,\EuO)[-1]$, and the first summand vanishes. Thus we have a surjection of algebras 
\[\HH^\bullet(\EuT_0) \to \BbK[\beta,\gamma]  \]
with kernel $\h^1(\EuO)[-1]$. Since this map is an isomorphism in degrees $\geq 2$, we write $\beta$ to denote a class in $\HH^2(\EuT_0)$, etc.

Let $\eta$ be a generator for the 1-dimensional image of $\h^1(\EuO)$ in $\HH^1(\EuT_0)$. We claim that $\beta \eta =0$.  To see this, observe that the $E_\infty$-page of the Hodge spectral sequence is an algebra isomorphic to $\gr \HH^\bullet (C)$, the associated graded algebra for the filtration giving rise to the spectral sequence. The $E_\infty$-page here is supported along the $0$th row and the $0$th column (in non-negative total degrees). Hence $\gr^i \HH ^\bullet(\EuT_0) = \HH^i(\EuT_0)$ for $i\neq 1$. Thus, to show that $\beta \eta =0$, it suffices to show that it vanishes in $\h^1(\EuT_0,\mathcal{HH}^2)$, which is obvious since this module vanishes.

Now let $\gamma_1$ be any lift of $\gamma$ to $\HH^1(\EuT_0)$, and let $\gamma_2=\gamma_1+\eta$. Then 
$\{ \gamma_1,\gamma_2\}$ is a basis of $\HH^1(\EuT_0)$. Moreover, $\gamma_1 \beta = \gamma_2\beta$, and $\gamma_1\gamma_2$ is zero since it maps to $\gamma^2 =0 \in \BbK[\beta,\gamma]$. The result now follows.
\end{pf}

\section{\texorpdfstring{Weierstrass curves versus $A_\infty$-structures}{Weierstrass curves versus A-infinity structures}}

In this section we shall prove our `dg comparison theorem' \ref{DG comparison}, and refinements of it. 
\subsection{Cochain models and splittings}
Our plan is to reformulate Theorem \ref{DG comparison} in terms of minimal $A_\infty$-structures on the fixed algebra $A$, and prove it in sharper form in that language. To do so, we need homological perturbation theory. 
\begin{defn}
Let $(C^\bullet,\delta)$ be a cochain complex over a commutative ring $R$, with an action of a group $\Gamma$ by automorphisms. A \emph{splitting} for $C^\bullet$ is an internal direct sum decomposition 
\[ C^k =  \underbrace{\EuH^k \oplus \im\delta^{k-1}}_ {\ker \delta^k} \oplus  \EuI^k \] 
for each $k$, with $\Gamma$-invariant summands. Equivalently, it is a $\Gamma$-equivariant linear map $s\in \hom^{-1}(C,C)$ such that $\delta s\delta = \delta$
(for given given the direct sum decomposition we put $s|_{\EuH \oplus \EuI}=0$ and $s|_{\im \delta} = \delta^{-1}\colon \im \delta\to \EuI$, while given such an $s$ we have $C^\bullet =  \underbrace{ \left[ \ker \delta \cap \ker (\delta\circ s)\right] \oplus \im\delta}_ {\ker \delta} \oplus \im (s\circ \delta)$).
\end{defn}
The set of splittings will be denoted by $\spl_\Gamma(\EuA)$. When $\EuA$ is a $\Gamma$-equivariant $A_\infty$-algebra, with cohomology $A$, a splitting for $\EuA$ as a cochain complex gives rise, via the homological perturbation lemma \cite{SeiBook}, to a canonical, $\Gamma$-equivariant $A_\infty$-structure $\mu^\bullet_\EuA$ on $A$, together with equivariant $A_\infty$-homomorphisms $i\colon A \to \EuA$ and $p \colon \EuA\to A$ such that $i$ and $p$ induce the identity map $\id_A$ on cohomology; equivariance follows from the naturality of the construction. 

If one merely splits the cocycles, writing $\ker \delta = \im\delta\oplus \EuH$, but does not complement the cocycles, \emph{the conclusion is the same except that one does not get the map $p$} \cite{Kad}. One can find a splitting of the cocycles (not necessarily equivariant) whenever $A$ is a projective module.

\subsubsection{\texorpdfstring{A \v{C}ech model for Weierstrass curves}{A Cech model for Weierstrass curves}} Let $\EuO_W=\BbK[a_1,\dots,a_6]$, and and let $\EuC$ be the universal embedded Weierstrass curve $y^2z + a_1 xyz + a_3 y z^2= x^3 + a_2 x^2 z+ a_4xz^2 + a_6z^3$ over $\EuO_W$.  We are interested in cochain models $\EuB_{\EuC}$ for the endomorphism algebra $\ext^\bullet(T,T)$, where $T = \EuO \oplus \EuO_\sigma$ over the curve $\EuC$. A technical irritation is that $\EuO_\sigma$, though a perfect complex, is not locally free. However, there is an autoequivalence $\tau$ of $\tw \VB(\EuC)$---a twist along the spherical object $\EuO$ \cite{ST}---such that $\tau(\EuO)=\EuO$ and $\tau(\EuO_\sigma)=\EuO(-\sigma)$. There is an autoequivalence $\tau'$ which is inverse to $\tau$ in that $\tau\circ \tau'$ and $\tau'\circ \tau$ are both naturally isomorphic to the identity functor. There is an induced isomorphism $\ext^\bullet(T,T) \owns [h] \mapsto [\tau'\circ h\circ \tau] \in \ext^\bullet(T',T')$, where $T' = \EuO \oplus \EuO(-\sigma)$. 
We shall set up a cochain complex $\EuB_{\EuC}'$ which computes $\ext^\bullet(T',T')$, and use $\tau$ and $\tau'$ to obtain from it a complex $\EuB_{\EuC}$ which computes $\ext^\bullet(T,T)$. We observe that a splitting for the former complex will transfer to one for the latter. 

Our dg model for $\EuB_{\EuC}'$ will be the \v{C}ech complex 
\begin{equation} 
\EuB_{\EuC} ' : = \left(\check{C}^*(\EuU; \SheafEnd(T')),\delta\right), \quad \EuU = \{U,V\}, 
\end{equation}
where $U:=\{ z \neq 0\}$ is the complement of the point at infinity $\sigma=[0:1:0]$ and $V= \{ y \neq 0\}$. The multiplicative group $\BbK^{\times}$ acts on $C$, covering its an action on $W$. The action preserves each of the two sets in the covering $\EuU$; hence $\BbK^{\times}$ acts on $\EuB_{\EuC}'$ by automorphisms---a \emph{strict} action \cite[(10b)]{SeiBook}.

\begin{lem}\label{splittings}
The set $\spl_{\BbK^{\times}}(\EuB_C')$ of $\BbK^{\times}$-equivariant retractions is non-empty; there exists a distinguished splitting $r$.
\end{lem}

\begin{pf}
The sheaf  $\SheafEnd(T')$ has four line-bundle summands (`matrix entries'): $\EuO$ (twice), $\EuO(\sigma)$, and $\EuO(-\sigma)$. The complex $\EuB_\EuC'$  is the direct sum of the \v{C}ech complexes for these line bundles, so it suffices to handle them separately. Evidently, we must put $\EuI^1=0$ and $\EuH^0 = \ker \delta^0$. Our tasks are to identify $\EuH^1$ and $\EuI^0$ for the line bundles $\EuO$, $\EuO(\sigma)$ and $\EuO(-\sigma)$.
 
{\bf 1.} We consider endomorphisms of $\EuO$ or of $\EuO(-\sigma)$. We have  $\SheafEnd(\EuO)=\SheafEnd(\EuO(-\sigma))=\EuO$. To describe $\coker\delta=H^1(\EuO)$, take a function $g=\gamma_{UV} \in \Gamma(U\cap V,\EuO)$. If it is regular at $\sigma$ then $g = \delta(0, - \gamma_{UV})$.  If $g$ has a pole of order $d \geq 2$ at $\sigma$ then we can find a function $\zeta_U$ on $U$ such that $g-\zeta_U$ has a pole of order $< d$ at $\sigma$. On the other hand, if $g_{UV}$ has a simple pole at $\sigma$, with $\sigma(s) = \frac{1}{s}+O(1)$ with respect to a local uniformizer $s$ at $\sigma$, then $H^1(\EuO)= [g]\cdot \BbK[W]$. (The restriction of $g$ to the locus $\BbK[W][\Delta^{-1}]$ where the discriminant $\Delta$ is non-zero cannot be a coboundary, since there is no degree 1 rational map defined on an elliptic curve. It is of no consequence to our argument whether $g$ is a coboundary when we restrict to $\Delta=0$.)

We take $g= \gamma_{UV} = \frac{X^2}{YZ}$, which is indeed regular on $U\cap V= \{ YZ\neq 0\}$. One has a local uniformizer $s=X/Z$ at $\sigma$, and $C$ is cut out formally as the graph of $t\in \BbK[W] \series{s}$ with $t=s^3+ O(s^4)$; and $g= s^2/t = s^{-1} + O(1)$, as required. Thus, we put $\EuH^1 = \BbK[W]\cdot g$, and we then have $ \check{C}^1(\EuU; \EuO) =  \im \delta \oplus \EuH^1$. (\emph{Note:} $g$ restricts to the cuspidal fibres as a \v{C}ech coboundary, but that does not affect our argument.)

Turning to the degree 0 part, $\ker \delta^0$ consists of pairs $(c,c)$ with $c$ constant. A complement $\EuI^0$ is given by
\[  \EuI^0  =  \EuO_U \oplus I  \subset \EuO_U \oplus \EuO_V, \quad I= \{ \gamma_V\in \EuO_V :\gamma_V(\sigma) =0 \} . \]

{\bf 2.} We consider $\EuH om(\EuO,\EuO(\sigma)) = \EuO(\sigma)$.
One has  $ \check{C}^1(\EuU,\EuO(\sigma)) = \im\delta^0,$ so $\EuH^1 = 0$; and 
\[ \ker \delta^0 = \{ (c,c)\in \EuO(\sigma)_U \oplus \EuO(\sigma)_V:  c\text{ constant} \}.  \]
Thus
\[ \check{C}^0(\EuU,\EuO(\sigma)) = \ker\delta \oplus \EuI^0, \quad \EuI^0  = \EuO(\sigma)_U \oplus I,  \]
where $I$ consists of sections over $V$ which vanish at $\sigma$. 

{\bf 3.}
We consider $\EuH om(\EuO(\sigma),\EuO)= \EuO(-\sigma)$. There are no global sections, so we put $ \EuI^0 = \check{C}^0(\EuU,\EuO(-\sigma)).$
We have 
\[   \check{C}^1(\EuU,\EuO(-\sigma)) = \im\delta^0 \oplus \EuH^1,\]
where $ \EuH^1=\BbK[W] \cdot g$, for similar reasons to those explained in case 1. 
\end{pf}
We now come to a key point in the construction:

\begin{leftbar}
Take $r$ from Lemma \ref{splittings}. By applying the homological lemma to ($\EuB_\EuC', r)$, we obtain a minimal $A_\infty$ structure $\EuA_W = (A\otimes W,\mu^\bullet_W)$, linear over $\BbK[W]$. 
\end{leftbar}
Said another way, we obtain a family of minimal $A_\infty$-structures $\EuA_w=(A,\mu^\bullet_w)$ parametrized by $w\in W$, whose structure coefficients depend polynomially on $w$. 
The $\BbK^\times$-equivariance of $r$ implies equivariance of $\mu^\bullet_W$; precisely, for $w\in W_k$ (the $k$th graded part), we have
\begin{equation} 
\mu^d_{\epsilon^k w} = \epsilon^{d-2} \mu^d_w. 
\end{equation} 

The dg comparison theorem was stated in terms of abstract Weierstrass curves, but the $A_\infty$-version will be formulated using embedded Weierstrass curves. First we set up the relevant categories:
\begin{enumerate}
\item
The groupoid $\EuW$ whose objects are embedded Weierstrass curves, thought of as elements $w\in W$. The group of reparametrizations $G$, from (\ref{reparam}), acts on $\ob \EuW$; we set
\[\mathsf{mor}_{\EuW}(w_1,w_2) = \{ g \in G: g(w_1) =w_2\}.\] 

\item
The groupoid $\EuM$ of minimal $A_\infty$-structures $\EuA$ on $A$. Let $\hoch^{r+s}(A,A)^s$ denote the part of the $r$th Hochschild cochain group in cohomological degree $r+s$ and internal degree $s$,
\[ \hoch^{r+s}(A,A)^s =\Hom^s_{\BbK[W]}(A^{\otimes r},A).   \]
and write $\hoch^k(A,A)^{\leq 0}=\prod_{s\leq 0}{\hoch^k(A,A)^s}$: these are the truncated Hochschild cochains, and they govern deformations of $A_\infty$-structures \cite{SeiQuartic}. Let $\EuG$ denote the group of `rescalings and gauge transformations', namely, the group of $u\in \hoch^1(A,A)^{\leq 0}$ whose leading term $u^1 \in \Hom(A,A)$ is  $r\cdot \id$ for some $r\in \BbK^\times$. Then $\EuG$ acts on the objects of $\EuM$, $(u, \EuA)\mapsto u_*\EuA$. We think of $u$ as an $A_\infty$-functor $\EuA\to u_*\EuA$ acting trivially on objects. The morphisms in $\EuM$ are given by $\mathsf{mor}_{\mathbf{\EuM}}(\EuA,\EuA') = \{ u\in \EuG:  u_*\EuA = \EuA'\}/\sim$, where $\sim$ is the equivalence relation which identifies homotopic functors. `Homotopy' has the following meaning: say $u_1$ and $u_2$ are $A_\infty$-functors $\EuA\to \EuA'$ with the same action on objects. Their difference $D= u_1-u_2$ then defines a natural transformation, i.e., a morphism $D = (D^0,D^1,D^2,\dots) \in \hom^1_\EuQ(u_1,u_2)$, where $\EuQ$ is the $A_\infty$-category of non-unital functors $\mathsf{funct}(\EuA,\EuA')$, satisfying $\mu^1_\EuQ D=0$. One puts $D^0=0$ and $D^k = u_1^k-u_2^k$ for $k>0$. A \emph{homotopy} from $u_1$ to $u_2$ is a $T\in \hom^0_\EuQ(u_1,u_2)$, with $T^0=0$, such that $D=\delta_\EuQ T$. Functors are called \emph{homotopic} if a homotopy exists. 
\end{enumerate}

\begin{main}[$A_\infty$ comparison theorem]\label{Ainf comparison}
Let $\BbK$ be a normal commutative ring. The passage from Weierstrass curves over $\BbK$ to $\BbK$-linear minimal $A_\infty$-structures on $A$ defines a functor $\EuF \colon \EuW \to \EuM$. Precisely:
\begin{enumerate}
\item 
Each embedded Weierstrass curve $C=C_w$ gives rise to a minimal $A_\infty$-structure $\EuF(C) = \EuA_w = (A,\mu^\bullet_w)$, by applying the homological perturbation lemma to $\EuB_C$, defined via the open cover $\EuU$  and the splitting $r$. These $A_\infty$-structures depend polynomially on $w\in W$, and are $\BbK^\times$-equivariant, meaning that $\mu^d_{t\cdot w}= t^{d-2} \mu^d_w$. In particular, the cuspidal cubic $C_0$ gives rise to the minimal $A_\infty$-structure $\EuA_0$ with no higher products: $\mu^d_0 = 0$ for $d\neq 2$.

\item Each pair $(g,w)\in G\times W$, where $G$ is the group of projective transformations acting on Weierstrass curves, gives rise to an element $\EuF(g,w) = u^g(w) \in \EuG$ such that $u^g(w) _* \EuA_w = \EuA_{g(w)}$. Moroever, $u^g(t\cdot w)=t^{d-1}u^g(w)$ for $t\in \BbK^\times$.

\item One has $u^1(w)=\id$. For any $(g_1,g_2)\in G\times G$, there exists a homotopy
\[    u^{g_2g_1}(w) \simeq u^{g_2}(g_1w) \circ u^{g_1}(w)  \]
depending algebraically on $w$.
\end{enumerate}
The functor $\EuF$ commutes with base-change (i.e., the operation $\cdot \otimes_{\BbK}\BbK'$ when $\BbK\to \BbK'$ is a ring homomorphism. If $\BbK$ is either an integrally closed noetherian domain of characteristic zero, or a field then $\EuF$ is an equivalence of categories. 
\end{main}

The equivalence clause says that the following three properties hold:
\begin{itemize} 
\item
\emph{Essential surjectivity:} every minimal $A_\infty$-structure on $A$ is isomorphic in $\EuM$ to one of the form $\EuA_w$.

\item
\emph{Faithfulness}: If $g_1(w)=g_2(w)$, and if $u^{g_1}(w)$, $u^{g_2}(w) \in \ob \mathsf{funct}\, (\EuA_{w},\EuA_{g_1(w)}) $ are homotopic functors, then $g_1=g_2$.

\item
\emph{Fullness:} if $u\in \EuG$ and $u_*\EuA_{w_1} = \EuA_{w_2}$ then $u \in \ob \mathsf{funct}\,(\EuA_{w_1},\EuA_{w_2})$ is homotopic to $u^g(w_1)$ for some $g\in G$ such that $g(w_1)=w_2$.
\end{itemize}
The normality condition is there because that is the condition under which we know that the Weierstrass differential $\omega$ on the universal Weierstrass curve $C$ defines a section of the dualizing sheaf.

\begin{pf}[Proof of Theorem \ref{DG comparison} assuming Theorem \ref{Ainf comparison}.] Over a field $\BbK$, every dg structure $\EuB$ can be transferred to a quasi-isomorphic $A_\infty$-structure $\EuA$ on the cohomology $A\cong H^*(\EuB_C)$. By Theorem \ref{Ainf comparison}, $\EuA$ is gauge-equivalent to $\EuA_w$ for some $w$, which in turn is quasi-isomorphic to the dg category $\EuB_w$. Hence $\EuB\simeq \EuB_w$ in the $A_\infty$-sense. Furthermore, if $\EuB_{C_1} \simeq \EuB_{C_2}$, realize $C_i$ as an embedded Weierstrass curve $C_{w_i}$; so $\EuB_{w_1}\simeq \EuB_{w_2}$, and hence $\EuA_{w_1}\simeq \EuA_{w_2}$. Pick a gauge-equivalence $u \colon \EuA_{w_1}\to \EuA_{w_2}$, and then use the theorem to replace $u$ by a homotopic gauge-transformation $u(g)$, where $g(w_1)=w_2$. Hence $C_1\cong C_2$.
\end{pf}

\paragraph{Coherence.} One can ask whether the homotopies $H(g_1,g_2)$ from $u^{g_2g_1}$ to $u^{g_2}(g_1\cdot) \circ u^{g_1}$ can be chosen coherently in the sense of \cite[(10b)]{SeiBook}. Let $\tilde{\EuM}$ be the category of minimal $A_\infty$-structures on $A$, over $\BbK$, in which morphisms $\EuA\to \EuA'$ are rescaled gauge transformations $u\in \EuG$ such that $u_*\EuA=\EuA'$. Let $M= \aut(\id_{\tilde{\EuM}})$; it is an abelian group under composition (by an Eckmann--Hilton argument), and a $G(\BbK)$-module. The obstruction to coherence is a group cohomology class $o\in \h^2(G(\BbK); M)$. If one works over a field $\BbK$ and asks that the homotopies to be continuous in $g_1$ and $g_2$ the obstruction lies in the continuous group cohomology $\h^2_{\mathsf{cont}}(G(\BbK); M)$; if wants them algebraic, it lies in an algebraic version of group cohomology. 

We do not pursue these obstructions in detail, but content ourselves with an easy case. Work over $\C$, and let $G_{\mathsf{an}}$ denote $G(\C)$ with its analytic (not Zariski) topology. We ask that the homotopies be continuous on $G_{\mathsf{an}} \times G_{\mathsf{an}}\times \C^5$. If we restrict $g_1$ and $g_2$ to the uni-triangular normal subgroup $U_{\mathsf{an}}\subset G_{\mathsf{an}}$ then the obstruction lies in $\h^2_{\mathsf{cont}}(U_{\mathsf{an}}; M)$. But $U_{\mathsf{an}}$ is contractible, hence $\h^2_{\mathsf{cont}}(U_{\mathsf{an}};M)= H^2(BU_{\mathsf{an}}; M)=H^2(\{\mathsf{pt}.\}, M) =0$. Thus we \emph{can} make our homotopies coherent for $U_{\mathsf{an}}$. 

The Lyndon--Hochschild--Serre spectral sequence for $U_{\mathsf{an}}\subset G_{\mathsf{an}}$ is concentrated in a single row, which tells us that $\h^2_{\mathsf{cont}}(G_{\mathsf{an}};M) \cong \h^2_{\mathsf{cont}}(\C^\times, M^{U_{\mathsf{an}}})\cong M^{U_{\mathsf{an}}}$. So $M^{U_{\mathsf{an}}}$ is where to find the obstruction to extending the coherent homotopies from $U_{\mathsf{an}}$ to $G_{\mathsf{an}}$.

\subsection{\texorpdfstring{The functor $\EuF$ on hom-spaces} {The functor F on hom-spaces}}
The proof of the equivalence clause of Theorem \ref{Ainf comparison} will be given in the next section; for now, we shall set up the functor.

We have already set up the functor on objects---namely, we have constructed $\EuA_W$. Now take $g\in G$ and $w\in W$, and notice that they define an isomorphism $g_* \colon C_w\to C_{g(w)}$. We have a diagram of $A_\infty$ quasi-isomorphisms
\[ \xymatrix{ 							
										&												&				\EuB_{g(w), g(\EuU)\cup \EuU}\ar[dl]_{f}\ar[dr]^{f'} 	& 	\\
 \EuB_{w,\EuU} \ar^{g_*}[r] \ar_{p_w}[d]	& \EuB_{g(w),g(\EuU)}	\ar@{.>}@<-1ex>_s[ur]	& 	  		&  	\EuB_{g(w),\EuU}\ar_{p_{g(w)}}[d] 						\\
 \EuA_{w}\ar@<-1ex>_{i_w}[u]			&												&			&	\EuA_{g(w)}\ar@<-1ex>_{i_{g(w)}}[u]
}  \]
Here $f$ and $f'$ are maps of dga; they forget one of the three open sets in the covering $g(\EuU)\cup \EuU = \{ U, V,g(V)\}$. The dotted arrow marked $s$ is an $A_\infty$-homomorphism which is inverse to $f$, up to homotopy; $s$ is still to be constructed. The $A_\infty$-maps $i_w$ and $p_w$ are mutual inverses up to homotopy; they are associated with the splitting $r$, via the homological perturbation lemma \cite[Remark 1.13]{SeiBook}. 

Once $s$ has been constructed, we shall define $\EuF(g,w)$ as the $A_\infty$-composite $p_{g(w)}\circ f'\circ s\circ g_*  \circ i_w$. For composition of $A_\infty$-functors, see \cite[(1e)]{SeiBook}. 
To obtain $s$, we consider the following picture (cf. \cite[Cor. 1.14]{SeiBook}):
\[  \xymatrix{
	  	\EuB_{g(w), g(\EuU)} \ar_{g^*p_w}[d]							& \EuB_{g(w),g(\EuU) \cup \EuU} \ar^{f}[l] \\
	 	g^*\EuA_w\ar^{v^{-1}}[r]											& \tilde{\EuA}_{g(w)} \ar@<1ex>[l]^{v}  \ar_{\tilde{i}_{g(w)}} [u]
}\] 
The left vertical arrow is the pullback by $g$ of the $A_\infty$-morphism $p_w\colon \EuB_{w,\EuU}\to \EuA$, produced by means of the splitting $r$ and the homological perturbation lemma.  
\begin{leftbar}
\emph{Claim:}  the $\BbK[W]$-cochain complex $\EuB_{g(w),g(\EuU) \cup \EuU}$ admits a $\BbK^\times$-invariant splitting for its cocycles.
\end{leftbar}
The proof of the claim will be given below. The splitting gives rise to a minimal $A_\infty$-structure $\tilde{\EuA}_{g(w)}$ on its cohomology $A$, and the $A_\infty$-quasi-isomorphism $\tilde{i}_{g(w)}$. The composite $(g^*p_w) \circ f \circ \tilde{i}_{g(w)}$ is an $A_\infty$-morphism inducing the identity map $A\to A$ on cohomology. Thus it is a gauge transformation $v \in \hoch^1(A,A)^{\leq 0}$. As such, it has a strict inverse $v^{-1}$. We put  $s  =  \tilde{i}_{g(w)}\circ v^{-1}\circ (g^*p_w)$. Then $s\circ f$ and $f\circ s$ are homotopic to identity maps. (To prove this, use the fact that composition with a fixed functor preserves homotopy and observe that therefore, since $v=(g^*p_w) \circ f \circ \tilde{i}_{g(w)}$, one has $(g^*i_w)\circ v \circ \tilde{p}^*_w  \simeq f$.) 

We now consider the existence of homotopies $\EuF(h,g(w))\circ\EuF(g,w)\simeq \EuF(hg,w)$. For this, contemplate the diagram
\begin{equation}  \label{big diagram}
\xymatrix{ 							
&&\EuB_{hg(w), hg(\EuU)\cup h(\EuU)\cup \EuU}\ar@[blue][d]\ar@[cyan][dr] 
\\
&\EuB_{g(w), g(\EuU)\cup \EuU}\ar@[blue][d]\ar@[red][dr] \ar@<-1ex>_{h_*}[r] 
&\EuB_{hg(w), hg(\EuU)\cup h(\EuU)}\ar@[blue][l]  
&\EuB_{hg(w), h(\EuU)\cup \EuU}\ar@[magenta][dr]\ar[d]	
\\
\EuB_{w,\EuU} \ar@<-1ex>_{g_*}@[red][r] \ar_{p_w}[d]	 \ar@{.>}@/^3pc/[uurr] 
&  \EuB_{g(w),g(\EuU)}\ar@[blue][l]	\ar@{.>}@<-1ex>@[red][u]	
& \EuB_{g(w),\EuU} \ar_{p_{g(w)}}@[red][d] \ar@<-1ex>@[magenta]_{h_*}[r]
&  \EuB_{hg(w),h(\EuU)}\ar[l]	\ar@{.>}@<-1ex>@[magenta][u]	
& 	    	\EuB_{hg(w),\EuU}\ar_{p_{hg(w)}}@[magenta][d] 	 
\\
\EuA_{w} \ar@[red]@<-1ex>_{i_w}[u]	
&&\EuA_{g(w)}\ar@<-1ex>@[magenta]_{i_{g(w)}}[u] 
&&\EuA_{hg(w)}\ar@<-1ex>_{i_{hg(w)}}[u]
}  \end{equation}
In the lower part of the diagram we see the juxtaposition of the maps that go into the definitions of $\EuF(g,w)$ (red arrows) and $\EuF(h, g(w))$ (magenta arrows). In the top row is the \v{C}ech complex associated with a 4-set open cover. The long, curved  arrow pointing to it is a homotopy-inverse to the four-step composite formed by the blue arrows; it must be constructed. 

The arrows marked $g_*$ or $h_*$ have strict inverses $g^{-1}_*$ and $h^{-1}_*$. Now consider the arrow pointing down and right from the top of the diagram, shown in cyan. It is homotopic to the composite of five arrows going the other way round the hexagonal region.  Moreover, the composite of the four blue arrows in (\ref{big diagram}) is equal to the forgetful map  $\EuB_{hg(w), hg(\EuU)\cup h(\EuU)\cup \EuU}\to  \EuB_{hg(w),hg(\EuU)}$ followed by $(hg)_*^{-1}$.  As a result, we see that $\EuF(h, g(w)) \circ \EuF(g, w)$ is homotopic to the map $\EuA_{w}\to \EuA_{hg(w)}$ which factors through the curving and cyan arrows:
\begin{equation} \label{small diagram}
\xymatrix{ 							
&& \EuB_{hg(w), hg(\EuU)\cup h(\EuU)\cup \EuU}\ar[dl]\ar@[cyan][dr] \\
 \EuB_{w,\EuU} \ar^{(hg)_*}[r]\ar@{.>}@/^2pc/[urr]
 & \EuB_{hg(w),hg(\EuU)}\ar@<-1ex>@{.>}[ur] 	 
 && 	\EuB_{hg(w),\EuU}\ar_{p_{hg(w)}}[d]   \\
 \EuA_{w}\ar@<-1ex>_{i_w}[u]	
 &&& \EuA_{hg(w)}
 }\end{equation}
Therefore it will suffice to construct a homotopy-inverse to the latter forgetful map (this homotopy-inverse is indicated by the straight dotted arrow in (\ref{small diagram})). For this, it is sufficient to prove the 
\begin{leftbar}
\emph{Claim:}  the $\BbK[W]$-cochain complex $\EuB_{hg(w), hg(\EuU)\cup h(\EuU)\cup \EuU}$ admits a $\BbK^\times$-invariant splitting for its cocycles.
\end{leftbar}
 
The $A_\infty$-homomorphism $\EuA_w\to \EuA_{hg(w)}$ indicated by (\ref{small diagram}) is homotopic to $\EuF(hg,w)$. That is because the solid arrows make the following diagram commutative; hence the dotted homotopy-inverses form a homotopy-commutative diagram:
\[ \xymatrix{ 							
& \EuB_{hg(w), hg(\EuU)\cup h(\EuU)\cup \EuU}\ar[ddl]\ar[ddr]\ar[d] \\
& \EuB_{hg(w), hg(\EuU) \cup \EuU}\ar[dl]\ar[dr] \\
    \EuB_{hg(w),hg(\EuU)}\ar@<-.5ex>@{.>}[uur] \ar@<-.5ex>@{.>}[ur] 	 && 	\EuB_{hg(w),\EuU} 
 }\]
Hence $\EuF(hg,w)\simeq \EuF(h,g(w))\circ \EuF(g,w)$.  

It remains to prove the two claims highlighted above. Recall that we constructed a splitting for $\EuB_{w,\EuU}$ by constructing one for $\EuB'_{w,\EuU}$, meaning that we considered $\SheafEnd(\EuO\oplus \EuO(-\sigma))$ instead of $\SheafEnd(\EuO\oplus \EuO_\sigma)$. We shall do the same here. And as before, it suffices to consider the line bundles $\EuO$, $\EuO(\sigma)$ and $\EuO(-\sigma)$ which form the matrix entries for $\SheafEnd(\EuO\oplus \EuO(-\sigma))$. These sheaves have cohomology in any given degree which is either zero (in which case the splitting of the cocycles is trivial) or is a free module of rank 1, which means that one can complement $\im \delta$ in $\ker\delta$ by choosing any representative for the generator of cohomology---and the resulting splitting will be $\BbK^\times$-invariant. This establishes the two claims, and thereby completes the construction of the functor $\EuF$.

\subsection{Comparison of deformation theories}
Introduce the shifted Hochschild cochain complex $\EuD^\bullet$, given by
\[ \EuD^\bullet = (\hoch^\bullet(A,A)^{\leq 1}) [1]. \]
A minimal $A_\infty$-structure on the algebra $A$ is a sequence of maps $\mu^d\in \Hom_{\BbK}^{2-d}(A^{\otimes d},A)$ for $d\geq 2$ such that $\mu^2$ is the multiplication for $A$.  We can view the structure as a truncated  Hochschild cochain $\mu^\bullet\in \EuD^1$. Thus the functor $\EuF$ is defined, on objects, by a map $\mu^\bullet \colon W\to \EuD^1$. Pick a $\BbK^*$-equivariant splitting of the cocyles for $\EuB_\EuC$, where $\EuC$ is the universal Weierstrass curve over $\BbK[W]$; then we obtain, for each $w\in W$, a minimal $A_\infty$-structure $\alpha(w):=\mu^\bullet_w$ on $A$.
\begin{lem}\label{homogeneous}
Pick $i>2$ and  $(a_1,\dots,a_i)\in A$. The map $p\colon W_d \to \BbK$ given by $w\mapsto \mu^i_w(a_1,\dots,a_i)$ is zero except when $i-2$ is a multiple of $d$. It is linear if $i=d+2$, and in general is homogeneous of degree $(i-2)/d$.
\end{lem}
\begin{pf} 
The map $p$ is a polynomial function, equivariant under $\BbK^\times$, meaning that $p(\epsilon^d w) = \epsilon^{i-2} p(w)$, and is therefore homogeneous, of degree $(i-2)/d$.
\end{pf}
We now consider the derivative at $0$ of the map $w\mapsto \EuA_w$,
\[  \lambda := D \mu^\bullet |_{w=0} \colon W \to \EuD^1 . \]
Take $w\in W_d$. By definition, $\lambda_d(w)$ is a sequence $(x^2,x^3,x^4,\dots)$ where $x^i \in \Hom^{2-i}_\BbK(A^{\otimes i},A)$. By Lemma \ref{homogeneous}, $x^i=0$ except when $i=d+2$. Thus, we consider $\lambda|_{W_d}$ as a map 
\[ \lambda_d\colon W_d\to (\EuD^1)_{-d}.  \]
It is the map which assigns to $w$ its \emph{primary deformation cocycle} (see e.g. \cite{SeiQuartic}). Because $\mu^i_w=0$ for $i<d+2$, we have $\delta\circ \lambda_d=0$.

The effect of  $\EuF$ on morphisms is encoded in a map
\[ u\colon G\times W\to \EuD^0,  \]
which has a partial derivative
\[ \left. \frac{\partial u}{\partial g} \right|_{g=1} \colon \mathfrak{g} \times W\to \EuD^0.   \]

We define 
\[\kappa^0\colon \mathfrak{g}\to \EuD^0,\quad \xi \mapsto  \left. \frac{\partial u}{\partial g} \right|_{g=1}(\xi,0). \]
When $\xi \in \frakg_d$, the only non-vanishing component of $\kappa^0(\xi)$ lies in $\EuD^0$.

\begin{lem} 
The vertical maps $\kappa^j$ in the diagram 
\[\xymatrix{
& 0 \ar[r] &\frakg \ar^{d}[r]\ar_{\kappa^0}[d] & W \ar[r]\ar^{\kappa^1=\lambda}[d] & 0\ar[d]^{\kappa^2}\\
0 \ar[r] &\EuD^{-1} \ar^{\delta}[r] &\EuD^0 \ar^{\delta}[r]&\EuD^1 \ar^{\delta}[r] &\EuD^2
}\]
define a map $\kappa^\bullet$ of cochain complexes.
\end{lem}

\begin{pf}
We have seen that $\delta\circ \kappa^1=0$. To prove that $\delta\circ \kappa^0= \kappa^1 \circ d$, observe that the action of gauge transformations on $A_\infty$-structures is defined through a map $a\colon \EuD^0 \times  \EuD^1 \to \EuD^1$. We have $a(u^g(0),\mu^\bullet_0) = \mu^\bullet_{g(0)}$; differentiating this relation with respect to $g$, and setting $g=1$, we obtain the sought equation. 
\end{pf}

\begin{thm}\label{kappa iso}
The maps $\kappa^0$ and $\kappa^1$ induce isomorphisms 
\begin{align}
& [\kappa^0]\colon \ker d \to \Hoch^1(A,A)^{\leq 0}\\
& [\kappa^1]\colon  \coker d \to \Hoch^2(A,A)^{\leq 0}
\end{align}
when the base ring $R$ is 
\begin{enumerate}
\item [(i)] a field $\BbK$; or
\item [(ii)] the ring of integers $\Z$; or more generally
\item [(iii)] an integral domain of characteristic zero.
\end{enumerate}
\end{thm}
\begin{pf}
(i) We begin with $[\kappa^1]$. We claim that $\ker [\kappa^1]=0$. Indeed, if $w$ lies in in this kernel, let $C_{wt}$ be the Weierstrass curve $C$ over $\BbK[t]/t^2$ with parameters $a_i = t w_i$. Thus $C_{wt}$ specializes to $\Ccusp$ at $t=0$. The resulting minimal $A_\infty$-algebra $\EuA_{wt}$ is then formal over $\BbK[t]/t^2$; indeed, its class $[\kappa^1 (wt)] \in t \HH^2(A,A)$ exactly measures non-triviality of $\EuA_{wt}$ as an $A_\infty$-deformation of $A$. Since the quasi-isomorphism class of $\EuA_{C_{wt}}$ determines the curve, one finds that $C_{wt}\cong \Ccusp\times_{\spec \BbK} \spec \BbK[t]/t^2$. It follows that $wt$ defines a trivial first-order deformation of $\Ccusp$. Hence $w \in \im d$.

We quote from \cite{LP} or section \ref{Hochschild section} the result that, as graded $\BbK$-modules, we have
\[ \HH^2(A,A)^{\leq 0} =  \stackrel{s=-1}{\BbK/(2)} \oplus  \stackrel{-2}{\BbK/(3)}\oplus  \stackrel{-3}{\BbK/(2)} \oplus \stackrel{-4}{\BbK} \oplus \stackrel{- 6}{\BbK}.\]
Thus $\coker\delta$ is abstractly isomorphic to $\HH^2(A,A)^{\leq 0}$ as a graded vector space, and hence the injection $[\kappa^1]$ is an isomorphism.

Now consider $[\kappa^0] \colon \ker d \to \HH^1(A,A)^{\leq 0}$. Recall that $W$ is a $\mathfrak{g}$-module. On the level of cohomology, $W_\frakg= \coker d$ is a $\ker d$-module. Moreover, this module is easily checked to be faithful: that is, the action homomorphism $\ker d \to \End(W_\frakg)$ is injective. Moreover, $\Hoch^2$ is a $\Hoch^1$-module, and the map $[\kappa^\bullet]$ respects the actions on cohomology:
\[  [\kappa^0 \xi ] \cdot [\kappa^1 w] = [\kappa^1(\xi\cdot w)].  \] 
Given $\xi \in \ker d$, pick a $w\in W_\frakg $ such that $\xi\cdot w \neq 0 \in W_\frakg$. We then have $[\kappa^0(\xi)]\neq 0$. Hence $[\kappa^0]$ is injective. Both domain and codomain of $[\kappa^0]$ are isomorphic to 
\[ \stackrel{s=0}{\BbK} \oplus  \stackrel{-1}{\BbK/(2)} \oplus \stackrel{-2}{\BbK/(3)} \oplus \stackrel{-3}{\BbK/(2)};
 \]
therefore $[\kappa^1]$ is an isomorphism.

(ii) It will be helpful to note at the outset that for any ring $R$ one has $CC^\bullet_R(A\otimes R, A\otimes R) = CC^\bullet(A,A)\otimes R$, that $\mathfrak{g}_R = \mathfrak{g}_\Z \otimes R$, and that these canonical isomorphisms are compatible with the construction of the map $\kappa^{\bullet}$. Let $\BbK=\Z/(p)$ be a residue field of $\Z$. Form the $\Z$-cochain complex $K^\bullet = \cone \kappa^\bullet$, a complex of free abelian groups, and note, using (i), that $H^1(K^\bullet \otimes \BbK)=H^2(K^\bullet \otimes \BbK)=0$ while $\Hoch^0_{\BbK}(A\otimes \BbK,A\otimes \BbK)^{\leq 0} \to H^0(K^\bullet \otimes \BbK)$ is an isomorphism (so the latter vector space is 1-dimensional). By universal coefficients \cite[3.6.2]{Wei}, $H^j(K^\bullet  \otimes \BbK)$ has a direct summand $H^j(K^\bullet ) \otimes \BbK$. Consequently $H^1(K^\bullet ) \otimes \BbK= H^2(K^\bullet )\otimes\BbK=0$. Since $H^j(K^\bullet)$ is a finitely-generated abelian group, one deduces  $H^1(K^\bullet ) = H^2(K^\bullet )= 0$. One then has $\mathsf{Tor}_1^{\Z}(H^1(K^\bullet), \BbK)=0$, and so by universal coefficients again, $H^0(K^\bullet)\otimes \BbK \cong H^0(K^\bullet \otimes\BbK) \cong \BbK$. Hence $H^0(K^\bullet)$ is free of rank 1. Over $\Z$, one has $\ker \delta \cong \Z$, and a look at the exact sequence of the mapping cone then tells us that $\HH^0(A,A)^{\leq 0} \to H^0(K^\bullet)$ is an isomorphism. The same exact sequence then tells us that $[\kappa^0]$ and $[\kappa^1]$ are isomorphisms. 

(iii) Take the $\Z$-cochain complex $K^\bullet$ as before. Universal coefficients now tells us that 
\begin{align*}
H^0(K^\bullet \otimes_\Z R) &\cong H^0(K^\bullet)\otimes_\Z  R \oplus \mathsf{Tor}^1_\Z(H^1(K^\bullet), R) \cong R, \\
H^1(K^\bullet \otimes_\Z R) &\cong H^1(K^\bullet)\otimes_\Z  R \oplus \mathsf{Tor}^1_\Z(H^2(K^\bullet), R) = 0, \\
H^2(K^\bullet \otimes_\Z R) &\cong H^2(K^\bullet)\otimes_\Z  R \oplus \mathsf{Tor}^1_\Z(H^3(K^\bullet), R) = \mathsf{Tor}^1_\Z(H^3(K^\bullet),R).\\
\end{align*}
The ring $R$ is torsion-free as an abelian group, so $\mathsf{Tor}^1_\Z(H^3(K^\bullet),R)=\mathsf{Tor}^1_\Z(R,H^3(K^\bullet))=0$. The exact sequence of the mapping cone then tells us that $[\kappa^1]$ is an isomorphism. Over $R$, one has $\ker[\delta\colon \mathfrak{g}\to W]\cong R$, again because $R$ is torsion-free. Part of the exact sequence of the mapping cone reads
\[   0\to \Hoch^0(A_R,A_R)^{\leq 0}\to R \to R \xrightarrow{[\kappa^0]} \Hoch^1(A_R,A_R)^{\leq 0}\to 0, \]
and since $\Hoch^0(A_R,A_R)^{\leq 0}$ is non-zero, and $R$ an integral domain, the map $R\to R$ must be zero. Hence $[\kappa^0]\colon R\to \Hoch^1(A_R,A_R)^{\leq 0}$ is also an isomorphism. 
\end{pf}

\begin{rmk}
Out of caution, work over a field $\BbK$ in this remark. The map $[\kappa^1]$ has a straightforward deformation-theoretic meaning (a first-order deformation of Weierstrass curves gives a first-order deformation of $A_\infty$-structures on $A$). The map $[\kappa^0]$ may then be characterized as the unique map that makes the following diagram commute:
\[ \xymatrix{
\ker d \ar[r]\ar[d]_{[\kappa^0]}  & \End W_\frakg \ar[d]^{[\kappa^1]} \\
\Hoch^1(A,A)^{\leq 0} \ar[r] & \End \Hoch^2(A,A)^{\leq 0}.
} \]
The horizontal arrows are the module actions. This leads to a derived-categorical construction of $[\kappa^0]$ as the composite of the canonical isomorphisms of graded $\BbK$-modules
\begin{align*}
 \ker d &  \xrightarrow{\cong} \h^0(\Ccusp, \EuT)^{\leq 0} && \text{Prop. \ref{HH for Weierstrass}(2)} \\
&  \xrightarrow{\cong} \HH^1(\Ccusp)^{\leq 0} && \text{Lemma \ref{HH to H0}} \\
&  \xrightarrow{\cong} \HH^1(\tw \VB(\Ccusp))^{\leq 0} && \text{Prop. \ref{perf computes Hochschild}} \\
&  \xrightarrow{\cong} \HH^1(\EuB'_{\mathsf{cusp}},\EuB'_{\mathsf{cusp}})^{\leq 0} && \text{Prop. \ref{HH for Weierstrass}} \\
&  \xrightarrow{\cong} \HH^1(A,A)^{\leq 0} && \text{Prop. \ref{formality}}.
\end{align*}
\end{rmk}

\begin{rmk}
Before embarking on the proof of Theorem \ref{Ainf comparison}, we say a word about the methodology. We are very close here to the framework for deformation theory which uses differential graded Lie algebras (DGLA).  For instance, the DGLA $\EuK^\bullet = \hoch^{\bullet-1}(A,A)^{\leq 0}$ determines, in characteristic zero, a deformation functor which assigns to a local artinian $\BbK$-algebra $A$ the solutions in $\EuK^1\otimes m_A$ to the Maurer--Cartan equation $\delta \mu^\bullet + \half[\mu^\bullet, \mu^\bullet] =0$, modulo the group $\exp(\EuK^0\otimes m_A)$ (see e.g. \cite{Man}). A standard approach would be to show that $\kappa^\bullet$ is a map of DGLA, and conclude, given its effect on cohomology, that the deformation theories controlled by $\frakg\oplus W$ and $\EuK^\bullet$ coincide \cite{Man}. 

The route we have actually taken is a variant of this standard approach. Everything here works over arbitrary fields; there is no need for characteristic zero.  A minimal $A_\infty$-structure on $A$ is a Hochschild cochain $\mu^\bullet \in \hoch^2(A,A)^{< 0}$ which satisfies $\delta \mu^\bullet + \mu^\bullet \circ \mu^\bullet=0$. The Gerstenhaber square $\mu^\bullet \circ \mu^\bullet$ agrees with $\half [\mu^\bullet , \mu^\bullet]$ when $2$ is invertible.  There are no artinian rings in the picture, but the length filtration of Hochschild cochains serves as a substitute.  The map $\kappa^\bullet$ is not quite a Lie algebra homomorphism; it might be possible to promote it to an $L_\infty$-homomorphism, but we have chosen to use the functor $\EuF$ more directly. 
\end{rmk}

\begin{pf}[Proof of Theorem \ref{Ainf comparison}]
We have set up the functor $w\mapsto \EuA_w$ over regular rings $R$, and must prove that it is an equivalence in the stated sense. 

\emph{Essential surjectivity:} If two minimal $A_\infty$-structures on $A$, with composition maps $\mathbf{m}^k$ and $\mathbf{n}^k$, agree for $k\leq 8$, then the two structures are gauge-equivalent. Indeed, one proves inductively that if $\mathbf{m}^k=\mathbf{n}^k$ for $k<d$ then one can find a gauge transformation $u$ such that $(u_*\mathbf{m})^k=\mathbf{m}^k$ for $k<d$ and $(u_*\mathbf{m})^d=\mathbf{n}^d$. To do so, one notes that $\mathbf{m}^d-\mathbf{n}^d$ defines a class in $\HH^2(A,A)^{2-d}$. This Hochschild module  is zero for $d>8$: when $R$ is a field, this holds by Theorem \ref{HH for cusp curve}; when $R=\Z$, it then follows by the universal coefficients; when $R$ is an integral domain of characteristic zero, it then follows by universal coefficients from the $\Z$ case (cf. the proof of Theorem \ref{kappa iso}). One uses a trivialization of $\mathbf{m}^d-\mathbf{n}^d$ to define the gauge transformation (cf. \cite{SeiQuartic}).

Now let $\EuA = (A, \mu^\bullet)$ be a minimal $A_\infty$-structure on $A$, over $\Z$. Our goal is to show that $\EuA$ is gauge-equivalent to $\EuA_w$ for some $w\in W$. We shall repeatedly apply gauge transformations to $\EuA$, without notating them. If $\EuA$ is formal then $\EuA\simeq \EuA_0$. If it is not formal then we apply a gauge transformation so as to arrange that $\mu^k=0$ for $2<k<d$ but $[\mu^d] \neq 0 \in \HH^2(A,A)^{2-d}$. By Theorem \ref{kappa iso}, one has $[\mu^d]= [\kappa^1(w)]$ for a unique $w \in W_{d-2}$; we may assume, by applying another gauge transformation, that in fact $\mu^d= \kappa^1(w)$. By Lemma \ref{homogeneous}, $\mu^k_w =0$ for $k<d$. Thus $\mu^k=\mu^k_w$ for $k\leq d$. The difference $\mu^{d+1}-\mu^{d+1}_w$ is a cocycle, as one checks using the $A_\infty$-relations. If it is exact, one can adjust $\mu^d$ by a gauge transformation which leaves $\mu^k$ untouched for $k\leq d$ such that $\mu^{d+1}$ equals $\mu^{d+1}_w$, whereupon $\mu^{d+2}-\mu_w^{d+2}$ is a cocycle, and we can repeat the process. What we find is that either $\mu^\bullet$ is gauge-equivalent to $\mu^\bullet_w$, or else there is a $d' \leq 8-d$ such that, after applying a  gauge transformation to $\mu^\bullet$, one has  $\mu^{d+k}=\mu_w^{d+k}$ for $k=0,\dots, d'-1$ and $[\mu^{d+d'}-\mu_w^{d+d'}]\neq 0$. In the latter case, write $\mu^{d+d'}-\mu_w^{d+d'}= \kappa^1(w')$, and consider $\mu^\bullet_{w+w'}$. The differences $\mu^k-\mu^k_{w+w'}$ can be killed by gauge transformations for $k \leq d+d'$. We continue in the same fashion; the process stops once has made $\mu^k$ agree with $\mu^k_{w+w'+\dots}$ for $k\leq 8$. 

\emph{Faithfulness:} 
We shall consider the case of automorphisms: say $g(w)=w$, and that $u^g(w)\simeq \id$ ($\simeq$ means `is homotopic to'); we shall show that then $g=1$. From this, the general case will follow: say $g_1(w)=g_2(w)$ and $u^{g_1}(w)\simeq u^{g_2}(w)$. Let $g= g_1^{-1}g_2$. Then $g(w) = w$ and  $u^g(w)\simeq (u^{g_1})^{-1}\circ u^{g_2}\simeq \id$.

Recall that the Lie algebra $\frakg$ is a graded $R$-module: $\frakg = \frakg_0 \oplus \frakg_{-1} \oplus \frakg_{-2} \oplus \frakg_{-3}$. There is an induced filtration by Lie subalgebras
\[  \frakg = \frakg_{\leq 0} \supset \frakg_{\leq -1} \supset  \frakg_{\leq - 2}  \supset \frakg_{\leq -3}\supset \frakg_{\leq -4} = 0,
\quad \frakg_{\leq j} = \bigoplus_{k\leq j} {\frakg_k}.  \]
One easily finds algebraic subgroups $G \supset G_{-1} \supset G_{-2} \supset G_{-3} \supset G_{-4}=\{1\}$, with $\mathsf{Lie}\,G_{j}=\frakg_{\leq j}$; for instance, $G_{-1}$ is the unipotent subgroup $U$.

We also filter the group $\EuG$, so that $\EuG = \EuG_{0} \supset \EuG_{-1} \supset \dots $, as follows: take $v=(v^1, v^2,v^3, \dots ) \in \EuG$ (so $v^j\in \hoch^1(A,A)^{1-j}$, and $v^1 = c\, \id$ with $c$ a unit). If $v^1=\id$, say $v\in \EuG_{-1}$; if in addition, $v^k=0$ for $1<k \leq d-1$, say $v\in \EuG_{1-d}$.  If $g\in G_{1-d}$ then $u^g(w) \in \EuG_{1-d}$.
If in addition $d>2$ then we have $[u_d] \neq 0 \in \HH^1(A,A)^{1-d}$.  Hence we have the following observation:

\begin{leftbar}
if  $g(w)=w$, if $u=u^g(w)\in \EuG_{1-d}$ for some $d>2$, and if $[u^d] \neq 0$, then $g\in G_{-d}$.
\end{leftbar}

Now let $u^g(w)=(u^1, u^2, u^3,\dots)$, where $u^j\in \hoch^1(A,A)^{1-j}$, with $g(w)=w$ and $\id \simeq u^g(w)$. Since $u^g(w)\simeq \id$, $u^1$ induces the identity on $A=H^*\EuA_w$. Hence $g$ lies in the unipotent subgroup $U=G_{-1} \subset G$. Hence $u^1=\id$, i.e., $u^g\in \EuG_{-1}$. Thus 
$\delta u^2 = 0$, but the fact that $u^g(w)\simeq \id$ implies that $u^2$ is a coboundary. By the highlighted observation, we deduce that $g\in G_{-2}$. We argue similarly that $u^2$ is a coboundary, hence that $g\in G_{-3}$, and finally that $u^4$ is a coboundary, hence that $g\in G_{-4} = \{1\}$.

\emph{Fullness}. Next consider the assertion that given $v\in \EuG$ and given $w_1,w_2$ such that $v_* \EuA_{w_1} = \EuA_{w_2}$, there is some $g \in G$ such that $g(w_1)=w_2$ and $v \simeq u^g(w_1)$. 

We first note a non-emptiness statement: if $v_* \EuA_{w_1} = \EuA_{w_2}$ where $v\in \EuG$, then there is some $h \in G$ such that $h(w_1)=w_2$. This is true because the formal diffeomorphism-type of the $A_\infty$-structure $\EuA_w$ determines the affine coordinate ring, and hence the curve. That is: 
\[ \mathsf{mor}_{\EuM}(\EuA_{w_1},\EuA_{w_2})\neq \emptyset \Rightarrow \mathsf{mor}_{\EuW}(w_1,w_2)\neq \emptyset.\] 
Likewise, if $v_* \EuA_{w_1} = \EuA_{w_2}$ where $v\in \EuG_{-1}$, then there is some $h \in U=G_{-1}$ such that $h(w_1)=w_2$. Indeed, if we know $\EuA_w$ up to gauge-equivalence then we can reconstruct not only the affine coordinate ring of the curve, but also the Weierstrass differential on the curve.

Hence, since $\EuW$ and $\EuM$ are groupoids, and $\EuF$ functorial, it will be enough to prove fullness when $w_1=w_2$. We wish to show that the map
\[ \EuF_{w,w} \colon  \mathsf{mor}_{\EuW}(w,w) \to  \mathsf{mor}_{\EuM}(\EuA_w,\EuA_w), \] 
which we already know to be injective, is also surjective.  Moreover, by a rescaling argument, we see that it suffices to prove this under the assumption that $v\in \EuG_{-1}$. 

Take some $v \in \EuG_{-1}$ with $v_*\EuA_w = \EuA_w$. If $v$ is homotopic to the identity, we are done. If not then there is a $d>1$ such that $v$ is homotopic to some $x\in \EuG_{-d}$, $x_*\EuA_w=\EuA_w$, where $[x^d]\neq 0 \in \HH^1(A,A)^{1-d}$. By Theorem \ref{kappa iso}, $[x^d]=[\kappa^0(\xi_1)]$ for some $\xi_1 \in \frakg_d$. We then have $g_1:= 1+\xi_1 \epsilon \in U(\BbK[\epsilon]/\epsilon^2)$, and $g_1([w]) = [w] \in W \otimes_{\BbK} \BbK[\epsilon]/\epsilon^2$. It is possible to lift $g_1$ to $g_2 = 1 + \xi_1\epsilon + \xi_2 \epsilon^2 \in U(\BbK[\epsilon]/\epsilon^3)$, with $\xi_2 \in \frakg_{-d}$. We wish to do so in such a way that  $g_2([w]) = [w]\in W \otimes_{\BbK} \BbK[\epsilon]/\epsilon^2$. The obstruction is the class of $g_2(w)-w$ in $T_w W / \im \rho_w= \coker \rho_w$. Now, $[\kappa^1]$ maps $W$ to $\HH^2(A,A)^{\leq 0}$, and carries $\im \rho_w$ to $[\EuA_w, \HH^1(A,A)^{\leq 0}]$. Because $x^d$ extends to a gauge transformation which preserves $\EuA_w$, $[\kappa^1(g_2(w)-w)] \in [w, \HH^1(A,A)^{\leq 0}]$. Hence $g_2(w)-w\in \im \rho_w$. Inductively, we extend $g_1$ to $1+\xi_1 \epsilon + \dots + \xi_3 \epsilon^3$ (mod $\epsilon^4$) with $\xi_j \in \frakg_{\leq j}$. We can then find a homomorphism $\theta\colon \BbK\to \frakg$ which has this series as its 3-jet, and we put $g=\theta(1)$. One has $g(w)=w$. 

We next ask whether $(u^g)^{-1} \circ v$ is homotopic to $\id \colon \EuA_w \to \EuA_w$. The obstruction we encounter now lies in $\HH^1(A,A)^{1-d'}$ with $d'>d$. Repeating the argument, we eventually obtain an  element $g\in U$ with $u^g \simeq v$.
\end{pf}

\subsection{\texorpdfstring{The Gerstenhaber bracket on $\HH^\bullet(A,A)^{\leq 0}$}{The Gerstenhaber bracket on HH(A,A)}}

\emph{Note: This subsection will not be used elsewhere, except in the variant method (\ref{variants}) of our proof-by-elimination that $\EuT_0$ is mirror to $T_0$.}

The truncated Hochschild cohomology $H^\bullet = \HH^\bullet(A,A)^{\leq 0}$, in addition to its graded algebra structure, carries the Gerstenhaber bracket $[\, , \,]$. These two operations make $H^\bullet$ a Gerstenhaber algebra over $\BbK$: the product is graded commutative, the bracket makes $H^\bullet[1]$ a Lie superalgebra, and the bracket is a bi-derivation for the product.  Moreover, the internal grading $s$ is additive under the bracket. Since the brackets can readily be computed from our results, we record it here. We use the notation of Theorem \ref{HH for cusp curve}.

\begin{thm}\label{bracket}
$Q^{\bullet,\leq 0}$ carries a Gerstenhaber bracket $[\cdot,\cdot]$, respecting the internal grading $s$ and making the (canonical) isomorphism $\HH^\bullet(A,A)^{\leq 0} \to  Q^{\bullet, \leq 0}$ from Theorem \ref{HH for cusp curve} a map of Gerstenhaber algebras. The Gerstenhaber bracket on $Q^{\bullet,\leq 0}$ is as follows. The Lie algebra $L=(Q^1)^{\leq 0} = (T\otimes \gamma)^{\leq 0}$ is given by
 \begin{align*}
&\BbK \gamma, &&  &&\text{if }6\neq 0\\
&\BbK \gamma \oplus  \BbK x\gamma, &&  [\gamma,x\gamma]= - \gamma, &&\text{if }3= 0\\
&\BbK \gamma \oplus  \BbK x\gamma \oplus \BbK y \gamma , 
&& [\gamma,x\gamma]=0, \;  [\gamma,y\gamma] =\gamma, \; [x\gamma,y\gamma]=x\gamma, &&\text{if } 2=0.
\end{align*}

One has $ \left[Q^2,Q^2 \right]=0. $
The brackets $[L, Q^2]$ are given by
\[\begin{array}{l l  l l l l}
& [\gamma,   x\beta] = -4 x\beta,  & [\gamma,   \beta] = -6 \beta , & & & \text{if }6\neq 0;\vspace{5mm} \\

& [x\gamma,  x^2\beta] = x^2\beta, & [x\gamma,  x\beta] =  -x\beta, & [x\gamma, \beta] =0 , & \\
& [\gamma ,   x^2\beta] = x\beta, &  [\gamma, x\beta] = -\beta, & [\gamma, \beta] = 0  & & \text{if }3= 0 ; \vspace{5mm}\\

& [y\gamma, xy\beta] = xy\beta,& [y\gamma, y\beta] = y\beta, & [y\gamma, x\beta]=0, & [y\gamma, \beta] =0,\\
& [x\gamma, xy\beta] =0, & [x\gamma, y\beta] = x\beta,&  [x\gamma, x\beta] = 0,& [x\gamma,\beta]=0, \\
& [\gamma, xy\beta] = x\beta, & [\gamma, y\beta] = \beta,& [\gamma, x\beta]= 0,& [\gamma,\beta]=0, &\text{if }2=0.
\end{array}\]
$Q^{\bullet,\leq 0}$ is generated as a unital $\BbK$-algebra by $L$ and $Q^2$. The remaining brackets are therefore determined by the Leibniz rule.
\end{thm}
\begin{pf}
The structure of $L$ was already computed in Prop. \ref{HH1}; the first assertion here is essentially a restatement. A Hochschild 2-cocycle $c$ for $A$ determines a first-order deformation of $A$ as an $A_\infty$-algebra, and that the Gerstenhaber square $[c\circ c] \in \HH^3(A,A)$ is the obstruction to lifting it to a second-order deformation. By Theorem \ref{kappa iso}, all first-order deformations of $A$ come (via $\kappa^1$) from deformations of $\Ccusp$ as a Weierstrass curve. These Weierstrass deformations lift to second order; hence the algebraic deformation of $A$ also lifts. This shows that $c\circ c=0$. The bracket on $\HH^2$ is given by $[a,b] = a\circ b + b\circ a = (a+b)\circ (a+b) - a\circ a - b\circ b$, so $[\HH^2,\HH^2]=0$.  The adjoint action of the Lie algebra $\HH^1(A,A)^{\leq 0}$ on $\HH^2$ is the natural action of infinitesimal ($A_\infty$) automorphisms of $A$ on first order deformations. But the Lie algebra $\h^0(\EuT_{\Ccusp})^{\leq 0}$ acts via Lie derivatives on Weierstrass deformations of $\Ccusp$. Namely, take $\xi \in \ker d = \h^0 (\EuT_{\Ccusp})^{\leq 0}$. We have a vector field $\rho(\xi)$ on $W$. Take $w\in W$, and regard it as a translation-invariant vector field on $W$; then the adjoint action of $\xi$ on $w$ is given by the Lie bracket of vector fields: $w\mapsto [\rho(\xi),w](0)$. One computes these brackets using the formulae (\ref{subst}). The maps $[\kappa^\bullet]$ from (\ref{kappa iso}) intertwine this Lie derivative with the adjoint action of $\HH^1$ on $\HH^2$.  The statement about generation is clear from the definition of $Q^{\bullet,\leq 0}$. 
\end{pf}
\begin{rmk}
Once the Lie algebra $L$ has been computed, another approach to obtaining the rest of the brackets is to compute the BV operator $\Delta\colon Q^{\bullet, s}\to Q^{\bullet-1, s}$. Along with the product, that determines the brackets, and it is sharply constrained by the bidegrees (which, for instance, force $\Delta \beta=0$). We have used this method as an independent check of the above formulae for $[L,L]$ and $[L,Q^2]$---in particular, as a check on the signs.
\end{rmk}

\part{Symplectic geometry and the mirror map}

\section{Fukaya categories}
\subsection{The exact Fukaya category}

\paragraph{$T_0$ as a Liouville manifold.} Let $T$ be a 2-dimensional torus marked with a point $z$ and equipped with a symplectic form $\omega$. In the complement $T_0=T\setminus \{z\}$, $\omega$ is exact; say $\omega|_{T_0}=d\theta$. 

A small punctured neighbourhood of $z$ is symplectomorphic to the negative end in the symplectization of the contact manifold $(c_z,\theta|_{c_z})$, where $c_z$ is a small loop encircling $z$. The disc $D_z$ bounded by $c_z$ has the property that there is a diffeomorphism $D_z \setminus \{z\} \to c_z \times [0,\infty)$, which is the identity on the boundary, under which $\theta$ pulls back to $\theta|_{c_z}e^{r}$. $T_0$ is a \emph{Liouville manifold} (see e.g. \cite{SeiBias}): an exact symplectic manifold $(M,\theta)$ whose Liouville vector field $\lambda$ is complete, with a given compact hypersurface $H$ enclosing a compact domain $D$ in $M$, such that $\lambda$ is nowhere vanishing on $M \setminus \interior{D}$. The Liouville structure of $T_0$ includes the selection of the curve $c_z$, but is independent of this choice, as well as that of the form $\theta$, up to \emph{Liouville isomorphism.} A Liouville  isomorphism $f\colon (M_0,\theta_0)\to (M_1,\theta_1)$ is a diffeomorphism such that $f^*\theta_1-\theta_0=dK$ for a compactly supported function $K\in C^\infty_{\mathsf{c}}(M)$.

\paragraph{Grading.} We specify a \emph{grading} of $T_0$ as a symplectic manifold, that is, a trivialization of the square of the canonical line bundle; in two dimensions, that amounts to an unoriented line field $\ell \subset T(T_0)$. Anticipating our discussion of $\EuF(T,z)$, we shall always choose $\ell$ extends over $z$ to a line field on $T$. Such line fields form a torsor for $C^\infty(T,\R P^1)$; one has $\pi_0 C^\infty(T,\R P^1) = H^1(T;\Z)$. 

\paragraph{Branes.} The objects of $\EuF(T_0)^{\exact}$, the exact Fukaya category of the Liouville manifold $T_0$, as defined in \cite[chapter 2]{SeiBook}, are `exact Lagrangian branes' $L^\#$: 
\begin{leftbar}
embedded closed curves $L \subset T_0$ such that $\int_L \theta =0$, equipped with spin structures and gradings. 
\end{leftbar}

A grading for $L$ is a homotopy-class of paths from $\ell|_{L}$ to $TL$ inside $T(T_0)|_L$. If $L$ and $L'$ are graded curves then a transverse intersection point $y\in L\cap L'$ has a degree $i(y) = \lfloor \alpha/\pi \rfloor +1$, where $\alpha$ is the net rotation of the path from $T_y L $ to 
$\ell_y$ to $T_y L'$. When $\ell$ is \emph{oriented}---pointing along $\beta$, say---the grading for a curve $L$ induces an orientation for $L$.  The sign $(-1)^{i(y)}$ is then the intersection sign $[L'] \cdot [L]$ (note the order!), regardless of which orientation for $\ell$ was selected; see Figure \ref{HMScorrespondence} for relevant examples.
There are four inequivalent spin structures on each curve $L$, defined by the two orientations and the two double coverings. 

The morphism-spaces in $\EuF(T_0)^\exact$ are $\Z$-linear cochain complexes;  $\hom_{\EuF(T_0)^\exact}(L^\#,L'^\#)$ is a cochain complex $CF(\phi(L^\#),L'^\#)=\Z^{\phi(L) \cap L'}$ computing the Floer cohomology $HF(L^\#,L'^\#)$; here $\phi(L)$ is the image of $L$ under the time-1 map $\phi$ of a specified Liouville isotopy. The degrees $i(y)$ define the grading of the complex; the spin structures determine the sign-contributions of the holomorphic bigons defining the differential. Floer cohomology is invariant under Liouville isotopies $\phi_t$ (i.e., 1-parameter families of Liouville automorphisms) in that $HF(\phi_t (L^\#),  L'^\#) \cong HF(L^\#, L^\#)$. As a result, Liouville-isotopic objects are quasi-isomorphic in $\EuF(T_0)^{\exact}$. 

One has $\hom_{\EuF(T_0)^\exact}(L^\#,L^\#) \simeq \Z[-1] \oplus  \Z$: this is just a formula for the Morse cochain complex for $L$. If $L^\star$ denotes the brane sharing the same oriented, graded curve as $L^\#$ but with the other double covering, then one has $\hom_{\EuF(T_0)^{\exact}}(L^\#,L^\star) \simeq \{ \Z \xrightarrow{2} \Z[-1] \}$, which is the Morse complex for $L$ with the local system with fiber $\Z$ and holonomy $-1$.

The composition map $\mu^2 \colon \hom(L_1,L_2) \otimes \hom(L_0,L_1) \to \hom(L_0,L_2)$ can be understood purely combinatorially \cite[(13b)]{SeiBook}, provided that $L_0$, $L_1$ and $L_2$ are in general position. The coefficient for $y_0$ in $\mu^2(y_2,y_1)$ is a signed count of immersed triangles bounding $L_0$, $L_1$ and $L_2$ in cyclic order, with convex corners at $y_1$, $y_2$ and $y_0$. We shall discuss the sign when we discuss the relative Fukaya category. When $(L_0, L_1,L_2)$ are not in general position, one moves $L_1$ and $L_2$ by exact isotopies so as to make them so. 

Beyond the differential $\mu^1$ and the composition $\mu^2$, the higher $A_\infty$-structure maps $\mu^d$ are defined through inhomogeneous pseudo-holomorphic polygons for a complex structure $j$ on $T_0$ compatible with the orientation and the Liouville structure at infinity. We refer to \cite{SeiBook} for the foundations; since we shall not make any direct calculations with the higher structure maps, we need not say more here.

\paragraph{Describing the objects.} Oriented simple closed curves in $T_0$ do not realize all free homotopy classes of loops, but just one free homotopy class per non-zero homology class. Indeed, suppose that $\gamma_0$ and $\gamma_1$ are two such curves representing the same class in $H_1(T_0)$. It follows from the `bigon criterion' \cite{FM} that they can be disjoined by isotopies in $T_0$. Assuming that they are disjoint, they divide $T$ into two annuli. Only one of the annuli contains $z$, and hence an isotopy from $\gamma_0$ to $\gamma_1$ can be realized in $T_0$. If, moreover, $\gamma_0$ and $\gamma_1$ are isotopic simple closed curves which are both exact then they cobound an immersed annulus of area zero, and hence are Liouville-isotopic. A non-exact simple closed curve representing a non-trivial homology class can be isotoped to an exact one (just take an isotopy of an appropriate flux). 

\paragraph{Generation.} Pick a basis $(\alpha,\beta)$ for $H_1(T_0;\Z)$ with $\alpha\cdot \beta = 1$. Let $L_0^\#$ and $L_{\infty}^\#$ be objects of $\EuF(T_0)^{\exact}$ of respective slopes $\alpha$ and $\beta$. In Figure \ref{HMScorrespondence}, the line field $\ell$ is chosen to be parallel to $\beta$. With that choice, we grade $L_0$ by the trivial homotopy from $\ell|_{L_0}$ to $TL_0$, and $L_\infty$ by the homotopy depicted, so that $i(x)=0$ where $x$ is the generator for $CF(L_0,L_\infty)$. We also orient these curves in the respective directions of $\alpha$ and $\beta$ (so,if $\ell$ is oriented in the direction of $\beta$, the orientations are those obtained from the orientation of $\ell$ and the grading of the curve). It is permissible to choose a different line field $\ell'$, so long as it is defined on all of $T$, twisting say, $a$ times along $L_\infty$ and $b$ times along $L_0$; but in that case we put the same number of extra twists into the gradings of $L_0$ and $L_\infty$, so that one still has $i(x)=0$. We take the spin-structures on $L_0$ and $L_{\infty}$ to be the \emph{non-trivial} double covers. A convenient way to keep track of double coverings of curves $L$ is to mark a point $\star_L \in L$, and declare the double cover to be trivial on $L\setminus \{ \star_L\}$ and to exchange the sheets over $\star_L$. Such stars appear in Figure \ref{HMScorrespondence}.

Let $\EuA$ be the full $A_\infty$ sub-category of $\EuF(T_0)^{\exact}$ with objects $(L_0^\#, L_{\infty}^\#)$. Every object of $\EuF(T_0)^{\exact}$ whose double covering is non-trivial is quasi-isomorphic, in $\tw \EuF(T_0)^{\exact}$, to a twisted complex in $\EuA$. This is because such objects are given, up to possible reversal of orientation and shift in degree, by iterated Dehn twists $(\tau_{L_\infty}^\#)^n(L_0^\#)$, $n\in \Z$, and the effect of the Dehn twist $\tau_{L_\infty^\#}$ on the Fukaya category is the twist functor along the spherical object $L_\infty^\#$, provided that the double covering of $L_\infty^\#$ it is non-trivial. For a systematic account, including the double-cover condition, see \cite[Theorem 17.16]{SeiBook}; a low-tech account of the case at hand is given in \cite{LP}.  

Curves with trivial spin structures are not quasi-isomorphic in $\tw \EuF(T_0)^\exact$ to object of $\tw \EuA$; in fact, we shall see later that these objects do not represent classes in the subgroup of $K_0(\tw \EuA)$ of $K_0(\tw \EuF(T_0)^{\exact})$. What is true, however, is that for any object $X$ of $\EuF(T_0)^{\exact}$, the direct sum $X\oplus X[2]$ is quasi-isomorphic, in $\tw \EuF(T_0)^{\exact}$, to a twisted complex in $\EuA$. One can write down an explicit twisted complex representing $X\oplus X[2]$ using \cite[Cor. 5.8, Theorem 17.16, and formula 19.4]{SeiBook}, bearing in mind that there is an exact isotopy $(\tau_{L_0} \tau_{L_\infty})^6 \simeq \tau_{c_z}$ where $\tau_\gamma$ is a symplectic Dehn twist along $\gamma$. The isotopy takes place in a compact domain containing the surface bounded by $c_z$, and is trivial near the boundary. The relation arises from the monodromy of an anticanonical Lefschetz pencil on $\C P^2$ with 8 of the 9 base-points blown up. Consequently, \emph{$\EuA$ split-generates $\EuF(T_0)^{\exact}$}, i.e., $\twsplit \EuA \to \twsplit \EuF(T_0)^{\exact}$ is a quasi-equivalence.

\subsection{The wrapped category}
The \emph{wrapped} Fukaya $A_\infty$-category $\EuW(M)$ of a Liouville manifold $(M,\theta)$ is set up in \cite{ASei}. It is a $\Z$-linear $A_\infty$-category containing $\EuF(M)^\exact$ as a full subcategory. Its objects are again certain Lagrangian branes $L^\#$. Precisely, $L$ must be a properly embedded, eventually conical Lagrangian submanifold, exact in the strong sense that $\theta|_L = dK$ for some $K\in C^\infty_{\mathsf{c}}(L)$. The brane structure consists of a spin structure and a grading on $L$, and both can be taken to be eventually translation-invariant on the conical end. The morphism spaces $\hom_{\EuW(M)}(L^\#,L'^\#)$ are cochain complexes $CW^*(L^\#,L'^\#)$ computing wrapped Floer cohomology $HW^*(L^\#,L'^\#)$, which is Lagrangian Floer cohomology $HF^*(\phi(L^\#), L'^\#)$ for a Hamiltonian diffeomorphism $\phi$ which `accerelates' on the conical end. That is, $\phi$ acts on the conical end $N\times \R_+$ as $\phi(x,r)= \phi^{\mathsf{Reeb}}_{f(r)}(x)$, where $\{\phi^{\mathsf{Reeb}}_t\}$ is the time $t$ Reeb flow on the contact cross-section $N$ and $f(r)$ is a function which increases rapidly in a precise sense. As a consequence $\phi$ `wraps' $L$ around the end many times, typically producing an infinity of intersections with $L'$.

An object in $\EuW(T_0)$ is either an object of $\EuF(T_0)$, or an eventually-straight oriented arc $A$, with a spin-structure, of which there is only one isomorphism class per orientation, and a grading. A grading is, again, a way to rotate from $\ell|_\alpha$ to $TA$ inside $T(T_0)|_A$. 

The oriented arc $A$ has an initial segment $\{a_{\mathsf{in}}\} \times [r,\infty)$ and  a final segment $\{a_{\mathsf{out}}\}\times [r,\infty) $, where $a_{\mathsf{in}}$ and $a_{\mathsf{out}}$ are points on the circle $c_z$. Exact arcs with fixed $a_{\mathsf{in}}$ and $a_{\mathsf{out}}$ are Liouville-isotopic if and only if they have the same slope, i.e., represent the same primitive class in $H_1(T,  \{z\})$. Exact arcs with the same slope but different endpoints are not Liouville-isotopic, but are nonetheless quasi-isomorphic in $\EuW(T_0)$.

If $L^\#$ is an object of $\EuW(T_0)$ such that $L$ is closed, and $\Lambda^\#$ any object of $\EuW(T_0)$ such that $[L] \not \in \{ [\Lambda], -[\Lambda]\}$ in  $H_1(T,\{z\})$, one can apply a Liouville isotopy to $L_0$ so as to reduce the number of intersections to the minimal intersection number $i(L,\Lambda)$. There are then no bigons, so $HW^*(L^\#,\Lambda^\#)= \Z^{i(L,\Lambda)}$.
 
If $A$ is an arc of the same slope as the closed curve $L$, one can apply a Liouville isotopy to $A$ that disjoins it from $L$, whence $HW^*(A^\#,L^\#)=0$. It then follows that $HW^*(L^\#,A^\#)=HW^{1-*}(A^\#,L^\#)^\vee = 0$: this is an instance of Floer-theoretic Poincar\'e duality, which applies only when one of the objects is a closed Lagrangian.  One can show that $HW^*(A^\#,A^\#)$ is the non-commutative (tensor) algebra $T(u,v)$ on generators of degree 1 modulo the two-sided ideal $(u^2,v^2)$, but we shall not need to use this assertion.

\paragraph{Generation.} Let $A$ be any non-compact exact arc, equipped with a brane structure so as to make it an object of $\EuW(T_0)$. Let $L_0$ be a simple closed curve of the same slope as $A$, and $L_\infty$ a curve which intersects $L_\infty$ transversely at a single point. Equip these three curves with brane structures, with the double-coverings of the two closed curves both non-trivial. Then $\{ A^\#, L_0^\#, L^\#_\infty\}$ split-generates $\EuW(T_0)$. Indeed, we have already seen that 
$\{L_0^\#, L^\#_\infty\}$ split-generates $\EuF(T_0)^{\exact}$, so we need only consider the arcs. Any arc can be obtained from $A$ by a sequence of Dehn twists along $L_0$ and $L_\infty$, and so by an easy adaptation of \cite[Theorem 17.16]{SeiBook}, or by a much more elementary argument which applies to the surface case, can be represented as a twisted complex in $A^\#$, $L_0^\#$ and $L_\infty^\#$.

\subsection{The relative Fukaya category}
The relative Fukaya category $\EuF(T,z)$ (cf. \cite{SeiDef, SeiQuartic, She}) has the same objects as $\EuF(T_0)^{\mathsf{ex}}$; and 
\[\hom_{\EuF(T,z)}(L^\#,L'^\#)=\hom_{\EuF(T_0)^{\mathsf{ex}}}(L^\#,L'^\#)\otimes_\Z \Z\series{q}, \] 
i.e., $\hom_{\EuF(T,z)}(L^\#,L'^\#)$ is a free $\Z\series{q}$-module on the intersections $L\cap L'$. 
The line field $\ell$ on $T$ defines a grading, and hence makes the hom-spaces graded modules.

The $A_\infty$-structure $\{\mu^d\}$ is defined through inhomogeneous pseudo-holomorphic polygons, now in $T$. Such polygons $u$ count with a weight $\varepsilon(u) q^{u\cdot z}$, where $\varepsilon(u)$ is a sign, defined just as in $\EuF(T_0)^{\mathsf{ex}}$, and $u\cdot z$ is the intersection number with $z$. A formula for $\varepsilon(u)$ is given in  \cite{SeiGen2} (see also \cite{LP}). We shall give here only a special case in which all corners of the polygon (i.e., $y_{k+1} \in \hom(L_k,L_{k+1})$ for $k=0,\dots, d-1$ and $y_0 \in \hom(L_d,L_0)$) have \emph{even index} $i(y_k)$. In that case, $\varepsilon = (-1)^s$, where where $s$ is the number of stars on the boundary (recall that the stars designate non-trivial monodromies for the double covers). 

The same perturbations that define the $A_\infty$-structure in $\EuF(T_0)^{\mathsf{ex}}$ succeed in defining an $A_\infty$-structure here too. The proof uses automatic regularity for holomorphic maps to surfaces \cite[(13a)]{SeiBook}, and is otherwise unchanged from the proof in the the exact case.  The resulting $A_\infty$-structure is an invariant of $(T,\omega, z;\theta)$. Moreover, up to quasi-equivalence it is independent of $\theta$; the proofs of these assertions are straightforward adaptations of their analogues for $T_0$ given in \cite{SeiBook}.

\paragraph{Generation.} The objects $L_0^\#$ and $L_\infty^\#$ split-generate $\EuF(T,z)$. The proof is essentially the same as for $\EuF(T_0)$. One can see curves with non-trivial double coverings as explicit twisted complexes in $L_0^\#$ and $L_\infty^\#$, cf. \cite{LP}. For arbitrary objects $X$, the sum $X\oplus X[2]$ is again a twisted complex in $L_0^\#$ and $L_\infty^\#$; the proof uses the main results of \cite{SeiBook} just as before, via the relation $(\tau_{L_0}\circ \tau_{L_\infty})^6 \simeq \id [2]$ in the graded symplectic mapping class group. It is significant here that $c_1(T)=0$, and more particularly that $\ell$ extends over $T$, since Seidel's argument depends on the presence of absolute gradings; with that point noted, the argument applies to $(T,z)$ as it does to $T_0$. 

\subsection{The closed-open string map}
Besides the `open string' invariants $\EuF(M)^\exact$ and $\EuW(M)$, Liouville manifolds $M$ have a `closed string' invariant, the symplectic cohomology algebra $SH^\bullet(M)$ (see \cite{SeiBias} for an exposition and foundational references).  In the first, place $SH^\bullet(M)$ is a graded-commutative graded ring; the grading depends on a choice of grading for $M$ as a symplectic manifold. It also comes with a ring map $v \colon H^\bullet(M)\to SH^\bullet(M)$, which in our grading convention is a map of graded rings; this pins down our normalization for the grading of $SH^\bullet(M)$, which in some other accounts (such as \cite{SeiBias}) differs from ours by $\dim_\C M$. As a simple algebraic variant, we can work with $SH^\bullet(M;\BbK)$, an algebra over the commutative ring $\BbK$.

\begin{lem}\label{symp coh}
For any commutative ring $\BbK$, there is an isomorphism of graded $\BbK$-modules
\[ \theta\colon  \BbK[\beta,\gamma_1,\gamma_2]/(\gamma_1\gamma_2,\beta(\gamma_1-\gamma_2))\to SH^\bullet(T_0;\BbK), \] 
canonical after a choice of basis of $H_1(T_0)$, where $\deg \gamma_1= \deg\gamma_2 = 1$ and $\deg \beta = 2$. 
Moreover, $\theta(\gamma_1)\cdot \theta(\gamma_2)=0$.
\end{lem}
When $\BbK$ is a field, the algebra on the right is isomorphic to $\HH^\bullet(\EuF_\BbK(\EuT_0)^\exact)$ by Theorem \ref{HH nodal}. We will see presently that $\theta$ respects products.
\begin{pf}
We begin with a comment about grading. The grading of $SH(T_0;\BbK)$ depends, \emph{a priori}, on the choice of line field $\ell$. If $\ell' = h(\ell)$ is another choice, obtained from $\ell$ by a map $g\colon T_0\to \R P^1$ representing a class $c =g^*(o) \in H^1(T_0;\Z)$, where $H^1(\R P^1;\Z) = \Z \, o$, the degrees of generators, which are 1-periodic Hamiltonian orbits $y$, change according to the formula $|y|_{\ell'}=|y|_\ell \pm 2\langle c, [y]\rangle$ (we do not bother with the sign). In the cochain complex for $SH(T_0)$ we shall describe, all orbits are null-homologous, and hence the grading of $SH(T_0)$ is independent of $\ell$. 

The graded $\BbK$-module $SH(T_0;\BbK)$ is described  in \cite[ex. 3.3]{SeiBias}. One uses an autonomous Hamiltonian which is a perfect Morse function $h$, accelerating appropriately on the cylindrical end. One has a natural map of algebras $v\colon H^\bullet(T_0;\BbK)\to SH^\bullet(T_0;\BbK)$ which is an isomorphism onto $SH^{\leq 1}(T_0;\BbK)$; this is the contribution of the minimum $m$ and two saddle-points $s_1$ and $s_2$ of $h$. We choose $h$ so that $([s_1],[s_2])$ is the chosen basis for the Morse homology $H_1(T_0)$. Define $\theta(1)=1$, $\theta(\gamma_1)=s_1$ and $\theta(\gamma_2)=s_2$. 

For each $q\geq 1$, there is a Reeb orbit $o_q$ which winds $q$ times around the puncture; $H^*(o_q;\BbK)$ contributes classes $c_{2q}\in SH^{2q}(T_0;\BbK)$ and $c_{2q+1} \gamma\in SH^{2q+1}(T_0;\BbK)$ which span those respective modules. We define $\theta(\beta^q)=c_{2q}$ and $\theta(\gamma_1\beta^q) = c_{2q+1}$.
We have $[s_1]\cdot [s_2]=0$, because $o_1$ is not contractible. 
\end{pf}

The exact Fukaya category of a Liouville manifold $M$ is tied to its symplectic cohomology via the `closed-open string map' \cite{SeiDef} to its Hochschild cohomology
\begin{equation}\label{CO} 
\EuC\EuO \colon  SH^\bullet(M;\BbK) \to  \Hoch^\bullet \left(\EuF(M)^{\exact}_\BbK\right). \end{equation}
For the details of the construction of $\EuC\EuO$ we refer to S. Ganatra's doctoral thesis \cite{Gan}. 

\begin{lem}[see \cite{Gan}]
$\EuC\EuO $ is a homomorphism of $\BbK$-algebras: it intertwines the pair-of-pants product on $SH^\bullet$ with the cup product on $\Hoch^\bullet$.
\end{lem}

\section{The punctured torus}
This part of the paper pinpoints the Weierstrass curve $C_{\mathsf{mirror}} \to \spec \Z\series{q}$ such that the minimal $A_\infty$-structure $\EuA_{\mathsf{mirror}}$ it induces on the algebra $A$ is gauge-equivalent to the $A_\infty$-structure $\EuA_{\mathsf{symp}}$ obtained from the Fukaya category $\EuF(T,z)$. Theorem \ref{Ainf comparison} assures us that $C_{\mathsf{mirror}}$ exists, and is unique as an abstract Weierstrass curve. Our aim is to show that $C_{\mathsf{mirror}}\cong \EuT$.  This will be accomplished in the final section of the paper by an argument involving $\theta$-functions.  In this section we offer two alternative proofs that  $C_{\mathsf{mirror}}|_{q=0} \cong \EuT_0$. One is by eliminating all possibilities other than $\EuT_0$; the other is by a calculation of `Seidel's mirror map'. We also prove our mirror-symmetry theorem for the wrapped category.

\subsection{\texorpdfstring{First proof that $C_{\mathsf{mirror}}|_{q=0} =  \EuT_0$: by elimination}{First identification of the mirror at q=0: by elimination}}
Lemma \ref{Tate at q=0} characterized the Weierstrass curve $\EuT_0\to \spec \Z$ as having a section which is a node at $p$, for any $p\in \spec \Z$. Symplectic topology now enters the picture:
\begin{prop}\label{elimination}
Take $p\in \spec \Z$, and $C\to \spec \mathbb{F}_p$ a Weierstrass curve (here $\mathbb{F}_0=\Q$). Suppose that the exact category $\twsplit \EuF(T_0)^{\mathsf{ex}}$, taken with coefficients in $\mathbb{F}_p$, is quasi-equivalent to $\tw\VB C$. Then $C$ is nodal.
\end{prop}

\begin{pf}
Let $\BbK=\mathbb{F}_p$.  One has Viterbo's map $v\colon H^\bullet(T_0;\BbK)\to  SH^\bullet(T_0;\BbK)$, a graded algebra homomorphism \cite{Vit, SeiBias}, and restriction homomorphisms
\[  \Hoch^\bullet(\EuF(T_0)^{\exact}) \to \Hoch^\bullet\left(\hom_{\EuF(T_0)^{\exact}}(L^\#,L^\#)\right),  \]
(everything over $\BbK$), one for each Lagrangian brane $L^\#\in \ob \EuF(T_0)^{\exact}$. The Hochschild cochain complex $\hoch^\bullet$ for an $A_\infty$-algebra \cite{SeiBook} has a filtration $F^r \hoch^\bullet$ by the length of cochains (so $\hoch^\bullet=F^0\hoch^\bullet\supset F^1 \hoch^\bullet \supset \dots$), and so for each $L$ one has a quotient map
\[  \Hoch^\bullet (\hom(L^\#,L^\#))\to  H(\hoch^\bullet/F^1\hoch^\bullet) = H^\bullet (\hom(L^\#,L^\#))= HF(L^\#,L^\#).  \]
Now, $HF(L^\#,L^\#) \cong H^\bullet(L;\BbK)$, ordinary cohomology,  by a canonical isomorphism \cite[(8c)]{SeiBook}. It follows easily from the definitions, plus the gluing theorem for Hamiltonian Floer theory---of which \cite{Sch} has a meticulous account---that the composite of the maps
\[ H^\bullet(T_0;\BbK) \xrightarrow{v} SH^\bullet(T_0;\BbK)\xrightarrow{\EuC\EuO} \HH^\bullet(\EuF(T_0)^{\exact}) \to \HH^\bullet(\hom(L^\#,L^\#))\to H^\bullet(L;\BbK)  \] 
is the classical restriction map $H^\bullet(T_0;\BbK)\to H^\bullet(L;\BbK)$---a \emph{surjective} map. Taking the sum of these composite maps for the objects $L_0^\#$ and $L_\infty^\#$ produces an \emph{isomorphism}
\begin{equation}\label{sum of composites}
H^1(T_0;\BbK)\to H^1(L_0;\BbK)\oplus H^1(L_\infty;\BbK). 
\end{equation}
We note one more feature of the composite $\EuC\EuO\circ v$, which is that it maps $H^1(T_0;\BbK)$ to $\Hoch^1(\EuF(T_0)^{\exact})^{\leq 0}$, the part spanned by cocycles in $F^1\hoch^1$. Indeed, the length-zero component of $\EuC\EuO\circ v$ returns, for each Lagrangian brane $L$, a count of index-zero pseudo-holomorphic discs attached to $L$, with one marked boundary point, of which there are none by exactness (constant discs have index $-1$). 

On the algebro-geometric side, suppose that $C=C_{\mathsf{mirror}}\to \spec \BbK$ is a Weierstrass curve mirror to $T_0$ over $\BbK$. By Lemma \ref{perf computes Hochschild}, one has $\Hoch^\bullet(\tw \VB C)\cong \Hoch^\bullet(C)$. Hochschild cohomology is invariant under passing to $\twsplit$, by a form of Morita invariance \cite{Kel,Toe}. Hence, under the hypotheses of the proposition, one has $\Hoch^\bullet \left(\EuF(T_0)^{\mathsf{ex}}\right)\cong  \Hoch^\bullet(C)$. Consequently, using (\ref{sum of composites}) we obtain a map of $\BbK$-algebras $ SH^\bullet(T_0;\BbK) \to \Hoch^\bullet(C)$ such that the composite
\[ H^1(T_0;\BbK) \to SH^1(T_0;\BbK) \to \Hoch^1(C)^{\leq 0} \to \ext^1(\EuO_C,\EuO_C)  \oplus \ext^1(\EuO_{C,\sigma},\EuO_{C,\sigma}) \]
is an isomorphism of $\BbK$-modules. We assert that such a homomorphism exists only if $C$ is nodal. To prove this, we must eliminate the cuspidal and smooth cases.

We claim that if $C$ were smooth, one would have $\HH^\bullet(C)\cong \Lambda^\bullet[\alpha_1,\alpha_2]$, where $\deg \alpha_1 = 1= \deg \alpha_2$. Additively, this follows the degeneration of the Hodge spectral sequence (\ref{Hodge}). The spectral sequence is multiplicative, and so the $\BbK$-algebra $\bigoplus_{p,q}{\h^p(\Lambda ^q \EuT_C)}$ is isomorphic to the associated graded algebra $E^\infty = \gr \Hoch^\bullet(C)$. We have $\alpha_1 \cdot \alpha_2\neq 0 \in \Hoch^2(C)$, since this is even true in the associated graded algebra; this establishes the claim. 
Any homomorphism of graded unital $\BbK$-algebras
\[ \theta \colon SH^\bullet(T_0) \to  \BbK[\alpha_1,\alpha_2] \]
obeys, in the notation of Lemma \ref{symp coh}, $\theta(\gamma_1)\theta(\gamma_2)=0$, and therefore $\theta$ fails to surject onto the 1-dimensional part 
$\{ \alpha_1,\alpha_2\}$. The composite of $\theta$ with the maps to $\ext^\bullet(\EuO_C,\EuO_C)$ and $\ext^\bullet(\EuO_{C,\sigma},\EuO_{C,\sigma})$ cannot then both be surjective.

We can rule out the possibility that $C$ is cuspidal, i.e., that $C\cong \Ccusp$, by noting that $\ext^1(\EuO_{\Ccusp},\EuO_{\Ccusp})\cong H^1(\EuO_{\Ccusp}) = \BbK\cdot \omega$. This module transforms under the $\BbK^\times$-action on $\Ccusp$, and it has weight $+1$. The restriction map $\Hoch^1(\Ccusp)\to \ext^1(\EuO_{\Ccusp},\EuO_{\Ccusp})$ respects the weight of the $\BbK^\times$-action---that is, it is a map of graded vector spaces---and hence its restriction to $\Hoch^1(\Ccusp)^{\leq 0}$ is zero. 
\end{pf}
\subsubsection{Variants}\label{variants}
There are other methods for ruling out $\Ccusp$:
\begin{itemize}
\item
When $6\neq 0$, the map that $\Hoch^1(\Ccusp)^{\leq 0}$ is 1-dimensional, and therefore the composite map $H^\bullet(T_0;\BbK) \to \ext^1(\EuO,\EuO)  \oplus \ext^1(\EuO_\sigma,\EuO_\sigma)$ unavoidably has a kernel. 
\item
In \cite{LP}, we used Abouzaid's model \cite{AboPlumb} for the Fukaya category of a plumbing to describe the structure maps of $\EuA$ and thereby prove non-formality, assuming $6\neq 0$.
\item
When $6=0$, the Gerstenhaber  bracket on $\HH^1(\Ccusp)^{\leq 0}$ is non-zero by Theorem \ref{bracket}.  For any Liouville domain $M$, the bracket is zero on the image of Viterbo's map $v\colon H^\bullet(M)\to SH^\bullet(M)$. Indeed, $v$ is a ring homomorphism, and it is easy to see that $\Delta \circ v =0$, where $\Delta$ is the BV operator on symplectic cohomology \cite{SeiBias}, and the bracket is, in accordance with the rules of BV algebras, 
\[  [x,y]=(-1)^{|x|}\Delta(x\cdot y)- x\cdot\Delta y - (-1)^{|x|}( \Delta x) \cdot y ,\] 
which implies that $[ v(a),v(b) ] =0$. Moreover, \emph{$\EuC\EuO$ preserves Gerstenhaber brackets}: this was first stated by Seidel \cite{SeiDef}, but there is no published proof. In characteristic 2 we have verified it for ourselves using standard gluing methods. Over other fields the signs are tricky, so we regard it as conjectural. However, taking this assertion for granted, the composite map $H^1(T_0;\BbK)\to \Hoch^\bullet(\Ccusp)^{\leq 0}$  preserves brackets and so cannot be injective.
\end{itemize}
\begin{pf}[Proof of Theorem \ref{mainth} clause (iii).]
We want to construct an $A_\infty$-functor
\[ \psi \colon \EuF(T_0)\to  \tw\VB(\EuT|_{q=0}). \]
We already have an isomorphism $H\psi\colon A = H^*\EuA \to H^*\EuB_C$, valid for any Weierstrass curve $C$ defined over $\Z$. Theorem \ref{Ainf comparison} implies that there is a unique $C$ for which $H\psi$ lifts to a quasi-isomorphism $ \psi\colon \EuA \to \EuB_C.$ For that particular $C$, $\psi$ extends naturally to a quasi-isomorphism
\[ \tw^\pi \EuA \to \tw^\pi \EuB_C. \]
Since the inclusion maps $\tw^\pi \EuA \to \tw^\pi \EuF(T_0)^{\exact}$ and $\tw^\pi \EuB_C \to \tw^\pi \VB(C)$ are quasi-equivalences, one obtains a quasi-isomorphism
\[  \tw^\pi \EuF(T_0)^{\exact} \to \tw^\pi \VB(C). \]
Composing this with a quasi-inverse to the inclusion $\tw\VB (C)\to \tw^\pi \VB(C)$, one obtains a quasi-isomorphism
\[  \tw^\pi \EuF(T_0)^{\exact} \to \tw \VB(C)  \]
whose restriction to $\EuF(T_0)$ is the functor we want.

Our task, then, is to identify the mirror Weierstrass curve $C$. In light of Prop. \ref{elimination}, it must be a curve which is nodal over $\Q$. \emph{A priori}, the node is only a $\overline{\Q}$-point. However, the normalization of $C$, defined over $\Q$, has two points (over $\overline{\Q}$) which map to the node, and since the normalization is a rational curve, these points are actually defined over $\Q$. Hence the same is true of the node. By clearing denominators, we obtain integer coordinates for the node. It then defines a section of $C\to \spec \Z$, which by the proposition must map to a node over $\mathbb{F}_p$ for every prime $p$. By Lemma \ref{Tate at q=0}, $C$ is therefore equivalent to $\EuT|_{q=0}$. 

We have proved a weakened form of Theorem \ref{mainth} clause (iii): we have shown that $\tw^\pi \EuF(T_0)^{\exact} \to  \tw\VB(C) $ is a quasi-equivalence, while the theorem claims that $\tw \EuF(T_0)^{\exact} \to  \tw\VB(C)$ is already an equivalence. We formulate this step as a separate statement, Prop. \ref{ex Fuk split closed}, whose proof completes that of the theorem.
\end{pf}

\begin{prop}\label{ex Fuk split closed}
$\tw\psi\colon \tw \EuF(T_0)^{\exact}\to \perf(C)$ is a quasi-equivalence.
\end{prop}

The triangulated $A_\infty$-category $\tw \VB(\EuT_0)$ is split-closed. Let $\EuI \subset \tw\VB(\EuT_0)$ denote the image of $\tw \EuF(T_0)^{\exact}$ under $\tw \psi$. We must show that the inclusion $\EuI \to \tw\VB(\EuT_0)$ is a quasi-equivalence. By Thomason's theorem \cite{Tho}, it is sufficient to prove equality of Grothendieck groups: $K_0(\EuI)=K_0( \tw \VB(\EuT_0))$.\footnote{By $K_0$ of a triangulated $A_\infty$-category $\EuC$, we mean $K_0(H^*\EuC)$, the Grothendieck group of the classical triangulated category $H^*\EuC$---in this case, $K_0(\perf \EuT_0)$.} 

\begin{lem}
Consider the map $s\colon \ob \EuF(T_0)^{\exact}\to \Z/2$ which maps an exact Lagrangian brane $L$ to $0$ if and only if the double covering $\tilde{L}\to L$ (part of the brane structure) is trivial. This map descends to a homomorphism $s\colon K_0(\tw \EuF(T_0)^{\exact})\to \Z/2$.
\end{lem}
\begin{pf} This is an adaptation of a result of Abouzaid \cite[Prop. 6.1]{AboSurfaces}. One views a spin-structure, lifting a given orientation, as a local system for the group $\{ \pm 1\} \subset U(1)$. The switch from closed, higher-genus surfaces to $T_0$ is irrelevant.
\end{pf}

\begin{pf}[Completion of the proof of \ref{ex Fuk split closed}.]
We must show that 
\[ K_0(\psi_0)\colon K_0(\tw \EuF(T)^{\exact})\to K_0(\tw \VB \EuT_0)\] 
is onto. We have $K_0(\tw \VB \EuT_0)=K_0(\EuT_0)$, and by Lemma \ref{K theory}, $(\rank, \det)\colon K_0(\EuT_0) \to \Z\oplus \Pic(\EuT_0)$ is an isomorphism.  One also has an isomorphism 
$(\deg,\rho)\colon \Pic(\EuT_0)\to \Z\oplus \Z^\times$, where $\deg$ is the degree and $\rho$ describes the `descent data' under normalization, as at the end of the proof of that lemma. The image of $\psi_0$ contains $\EuO$ (rank 1, degree 0) and $\EuO(\sigma)$ (rank 1, degree 1). Hence $(\rank,\deg)\circ K_0(\psi_0) 
\colon K_0(\tw \EuF(T)^{\exact})) \to \Z^2$ is onto. Now take some $L^\# \in \ob \EuF(T_0)^{\exact}$, and let $L^\star$ be the same object with the other double covering. Then the class $[L^\#]-[L^\star]$ is 2-torsion in $K_0(\tw \EuF(T)^{\exact})$, by the last lemma and the fact that the change of covering can be accomplished by an involution of $\tw \EuF(T_0)^{\exact}$ (namely, tensoring spin-structures on Lagrangians by the restrictions of some real line-bundle $\ell\to T_0$). Thus $[\psi_0(L)]-[\psi_0(L')]$ is again 2-torsion in $K_0(\tw \VB \EuT_0)$. There is a unique 2-torsion class in $K_0(\tw \VB \EuT_0)$, detected by $\rho$. Therefore $\rho\circ K_0(\psi_0)$ is surjective, and hence $K_0(\psi_0)$ is surjective.
\end{pf}

\begin{rmk}
There is an alternative to the argument just given which does not appeal to Abouzaid's analysis, but instead observes that the objects $L_0^\#$ and $L_0^\star$ (which differ only in their double coverings) map under $\psi$ to perfect complexes whoses $K_0$-classes differ by the generator of $\Z/2 \in K_0(\EuT)$. For this, we regard $L^\star_0$ as $L_0^\#$ with a local system with fiber $\Z$ and holonomy $-1$. We have $\psi(L_0^\#) = \EuO$, and it follows from an easy adaptation of Lemma \ref{where local systems go} below that $\psi(L_0^\star) = \EuO(\sigma-\sigma')$, where $\sigma'$ is the 2-torsion section of $\EuT_0^{\mathsf{sm}} = \mathbb{G}_{\mathsf{m}} (\Z) \to \spec \Z$ and, as usual, $\sigma$ is the identity section. From this point, the argument is straightforward.
\end{rmk}

\subsection{A second identification of the central fiber of the mirror curve: Seidel's mirror map} 

\paragraph{The affine coordinate ring.}
Suppose given an abstract Weierstrass curve $(C, \sigma,\Omega)$ over $\spec S$. There is then a Weierstrass cubic embedding carrying $\sigma$ to $[0:1:0]$; the affine complement to the closure of $\im \sigma$ is $\spec R_C$, where $R_C$, the affine coordinate ring, is the ring of functions on $C$ with poles only at $\sigma$:
\[  R_C = \varinjlim_n {H^0(C, \EuO(n\sigma))}. \]
\subsubsection{The Dehn twist ring.}
We want to compute $R_{C_{\mathsf{mirror}}}$ by determining an equation that holds in $\h^0(C_{\mathsf{mirror}}, \EuO(6\sigma))$. (This is almost the method of \cite{Zas}, but we are concerned here with the affine rather than homogeneous coordinate ring.) We shall in fact carry out this computation only for $C_{\mathsf{mirror}}|_{q=0}$.

Since $\EuA$ split-generates $\EuF(T,z)$, the quasi-isomorphism $\EuA\to \EuB_{C_{\mathsf{mirror}}}$ extends to an $A_\infty$-functor $\psi \colon \EuF(T, z)\to \tw \VB(C_{\mathsf{mirror}})$. Consider the object $L_\infty^\#$. The Dehn twist $\tau=\tau_{L_\infty^\#}$, acting as an autoequivalence of $\EuF(T_0)$, is homotopic to the twist functor along the spherical object $L_\infty^\#$: this is elementary in the present case \cite{LP}, but is an instance of a general result of Seidel's \cite{SeiBook}. Now, $\psi(L_\infty^\#)=\EuO_\sigma$. The twist along the spherical object $\EuO_\sigma \in \ob \tw \vect(C_{\mathsf{mirror}})$ is homotopic to the functor $\EuO(\sigma) \otimes \cdot$ of tensoring with $\EuO(\sigma)$ (see \cite[(3.11)]{ST}; the argument is carried out over fields, but over $\Z\series{q}$, Seidel--Thomas's map $f$ must be a unit times the restriction map by base-changing to reduce to the case of fields). 
Thus $\psi$ induces an isomorphism $R_{C_{\mathsf{mirror}}} \cong R_\tau$, where
\[  R_\tau:=  \varinjlim_{n}{HF^*(L_0^\#, \tau^n(L_0^\#))}. \]
The `Dehn twist ring' $R_\tau$, needs explanation---neither the direct system, not the ring structure, is  obvious. The maps
\[ \sigma_{n,m+n}\colon HF^*(L_0^\#, \tau^n (L_0^\#))\to HF^*(L_0^\# , \tau^{m+n}(L_0^\#))\] 
which form the direct system are defined via holomorphic sections of a Lefschetz fibration over a strip; this interpretation is part of Seidel's analysis of Dehn twists \cite{SeiBook}. The ring structure is easier: the $m$th power of the Dehn twist defines a map 
\[ (\tau_*)^m \colon HF^*(L_0^\#,\tau^n(L_0^\#)) 
\to HF^*(\tau^m(L_0^\#),\tau^{m+n} (L_0^\#)) \]
which applies the Dehn twist to the intersection points between Lagrangians. The product in the ring is given by composing this with the triangle product $\cdot$,
\begin{align*} 
HF^*(L_0^\#, \tau^n(L_0^\#)) \otimes HF^*(L_0^\#,\tau^m(L_0^\#))  & \xrightarrow{\tau_*^m\otimes\id}
HF^*( \tau^m(L_0^\#), \tau^{m+n}L_0^\#) \otimes HF^*(L_0^\#,\tau^m(L_0^\#) \\
& \xrightarrow{\cdot} HF^*(L_0^\#,\tau^{m+n}(L_0^\#)).  
\end{align*}
Associativity of this product is easily seen, as is the fact that the unit element $e\in HF^*(L_0^\#,L_0^\#)$ is a 2-sided unit for the multiplication in the Dehn twist ring. However, commutativity is something that we learn from the isomorphism $R_\tau\cong R_{C_{\mathsf{mirrror}}}$. 

\subsubsection{Avoiding the direct system}
If one knew the maps $\sigma_{m,m+n}$ explicitly, one would be able to proceed by perfect analogy with the algebro-geometric side of the mirror, as follows. Take the unit $e\in HF^*(L_0^\#,L_0^\#)$ and its images $e_n =\sigma_{0,n}(e) \in HF^*(L_0^\#, \tau^n(L_0^\#))$. Take a basis $\{e_2,x \}$ for $HF^*(L_0^\#,\tau^2(L_0^\#))$ and a basis $\{ e_2, \sigma_{2,3}x, y\}$ for $HF^*(L_0^\#,\tau^3(L_0^\#))$. Identify the Weierstrass equation as the unique relation satisfied by the seven monomials
$\{y^2,x^3,xye_1, x^2 e_2, y e_3, x e_4, e_6\}$.

We can extract nearly complete information about  the maps $\sigma_{0,n}$ as follows. We choose $L_0$ and $L_\infty$ to have just one, transverse intersection point. Then $L_0\cap \tau (L_0)$ consists of a single point $z'$, and $ HF^*(L_0^\#,\tau (L_0^\#)) = \Z\series{q}z'$. Hence $e_1 = f(q)z'$ for some $f(q)\in \Z\series{q}^\times$.
From the isomorphism $R_\tau \cong R_{C_{\mathsf{mirrror}}}$, we see that $e_n\cdot e_m = e_{m+n}$, and hence that $e_n=(e_1)^n=f(q)^n z'^n$. Moreover, $\sigma_{m,m+1}(u) = e_1 u = f(q) z' \cdot u$. In practice, then, one obtains a cubic equation by picking $x'$ so as to make $\{ z'^2, x' \}$ a basis for
$HF^*(L_0^\#,\tau^2(L_0^\#))$ and $y'$ so as to make $\{ z'^3, z' x' , y'\}$ a basis for $HF^*(L_0^\#,\tau^3(L_0^\#))$. One computes the products
\[ \{y'^2,x'^3,x'y'z', x'^2 z'^2, y' z'^3, x' z'^4, z'^6\}\]
and identifies the unique (up to scale) relation that they satisfy. This relation is necessarily of form
\[ y'^2 - f(q) x'^3 = \dots, \]
hence it determines $f$. Now let $x=fx'$ and $y=fy'$. Then one has $y^2 - x^3 = \dots$,
i.e., these coordinates satisfy the Weierstrass equation.
\begin{rmk}
We learned something interesting \emph{en route} here, though we shall not pursue it:  the series $f$, which encodes information about sections of a Lefschetz fibration, can be computed.  
\end{rmk}

\subsubsection{A model for the Dehn twist} \label{model twist}

To compute with the Dehn twist ring, note that one can take for $\tau$ 
any compactly supported exact automorphism of $T_0$ that is isotopic to a Dehn twist along $L_\infty$. To obtain a convenient
model, start with the linear symplectomorphism of $T=\R^2/\Z^2$ given by 
\[\delta [x_1,x_2] = [x_1, x_2 - x_1]. \]
The fixed-point set of $\delta$ is the line $L_\infty = \{ x_1 =0 \}$.
We will take our basepoint to be $z=(\epsilon,\epsilon)$ where $\epsilon \in (0,1/4)$. One has $\delta(z)=(\epsilon, 0)$; let $D$ be the $\epsilon^2$-neighborhood of the line segment $[\delta(z),z]$.  Take $\rho \in \aut_c(D,\omega|_D)$ to be a symplectomorphism such that $\rho(\delta(z))=z$; extend $\rho$ to a symplectomorphism of $T$, still called $\rho$, which is trivial outside $D$. Let $\tau = \rho \circ \delta$.  Then $\tau(z)=z$; by adjusting $\rho$, we may assume that $\tau$ acts as the identity in some neighborhood of $z$.
Thus $\tau$ restricts to give $\tau_0\in \aut_c(T_0,\omega|_{T_0})$. Let $L_{(1,-n)} = \{[x_1,x_2] \in T:  nx_1 + x_2 =0\}$. When $0\leq n \leq (2\epsilon)^{-1}$, we have that $\tau_0^n(L_0) = L_{(1,-n)}$. 

\begin{lem} 
There is a primitive $\theta$ for $\omega|_{T_0}$ making $L_0$  an exact curve and $\tau_0$ an exact symplectomorphism.  
\end{lem} 
\begin{pf} 
We must exhibit a primitive $\theta$ for $\omega$ such that $[\tau_0^* \theta -\theta]=0 \in H^1_c(T_0;\R)$ and $\int_{L_0}{\theta}=0$. For the first requirement, it suffices to show that $\int_{\gamma}(\tau_0^* \theta -\theta)=0$ for curves $\gamma$ forming a basis for $H_1(T, \{z\} ;\R)$. Such a basis is given by $\{L_\infty, L_0\}$, with chosen orientations for these two curves. For any primitive $\theta$ for $\omega$, one has $\int_{L_\infty}(\tau_0^* \theta -\theta)=\int_{\tau_0(L_\infty)}\theta - \int_{L_\infty}\theta = 0$. It suffices to choose $\theta$, therefore, in such a way that $L_0$ and $\tau_0(L_0)=L_{(1,-1)}$ are both exact Lagrangians. It is easy to find a 1-form $\iota$ on $T_0$ such that $d\iota = \omega$ on a small regular neighborhood $N$ of $L_0\cup L_{(1,-1)}$ and such that $\int_{L_0}{\iota}=0 =\int_{L_{(1,-1)}}\iota$. Now, $\partial N$ is a circle, isotopic to a loop encircling $z$. There is therefore no obstruction to finding a 1-form $\kappa\in \Omega^1(T_0)$, supported outside $L_0\cup L_{(1,-1)}$, such that $\omega=d(\iota+ \kappa)$; then $\theta=\iota+\kappa$ is the required primitive.
\end{pf}

\begin{prop}\label{primitive}
For any natural number $N$, one can choose a primitive $\theta$ for $\omega|_{T_0}$ such that the curves $L_{(1,-n)}$ are exact for $n=0,\dots,N$.
\end{prop}
\begin{pf}
Choose $\epsilon < (2N)^{-1}$, and take $\theta$ as in the lemma. Since $L_0$ is exact, so too is $\tau^n_0(L_0)$ for any $n\in \Z$.  But for $0\leq n \leq N$ we have $\tau^n_0(L_0) = \delta^n(L_0) = L_{(1,-n)}$,
\end{pf}

\subsubsection{Computation in the exact case}
Use a $\theta$ as in the proposition, taking $N$ at least $6$. We have
\begin{align*}  
L_0 \cap L_{(1,-1)} & = \{ z'\}, && z'=[0,0]; \\
L_0 \cap L_{(1,-2)} & = \{ \zeta_0, \zeta_1\}, && \zeta_k =[k/2,0]; \\
L_0 \cap L_{(1,-3)} & = \{ \eta_0, \eta_1,\eta_2 \},&& \eta_k =[k/3,0]. 
\end{align*}
In calculating products in the Dehn twist ring $R_\tau$, immersed triangles in $T_0$ count with sign $+1$. To avoid repetition, we do not give the argument here but defer it to Section \ref{homog subsection}. With this understood, it is straightforward to calculate that in the Dehn twist ring of $\EuF(T_0)^{\exact}$, one has
\[z'^2  = \zeta_0 + 2 \zeta_1. \]
We put $x' = \zeta_1$; then $\{ z'^2, x'\}$ is a $\Z$-basis for $HF^*(L_0^\#,L_{(1,-2)}^\#)$.
Next, we compute
\begin{align*} 
z' \zeta_0 & = \eta_0 + \eta_1 + \eta_2,\\
z' \zeta_1 & = \eta_1 + \eta_2.
\end{align*}
We deduce that
\begin{align*} 
z'^3 & = z' (\zeta_0 + 2\zeta_1) = \eta_0 + 3\eta_1 + 3\eta_2, \\
z' x' & =  \eta_1 + \eta_2.
\end{align*} 
We put $y'= \eta_2$, and note that $\{ z'^3, z'x, y'\}$ is a basis for  $HF^*(L_0^\#,L_{(1,-3)}^\#)$. Further computations yield
\[  \eta_2^2 =\theta_4,\quad \eta_1\eta_2 = \theta_3, \quad \zeta_1^3 = \theta_3
\]
where $(\theta_0,\dots,\theta_5)$ are the intersection points $\theta_k=[k/6,0]\in L_0\cap L_{-6}$.
These relations imply that
\[ y'^2 + x'^3  =  x' y' z'. \]
Putting $x=-x'$, $y=y'$ and $z=z'$, we obtain the Weierstrass relation in the desired form
\[  y^2 - x^3 = - xyz. \] 

\section{The wrapped Fukaya category}

We restate clause (iv) of Theorem \ref{mainth}, in slightly refined form. The category $\tw \VB(\EuT_0)$ is a dg enhancement for $\perf(\EuT_0)$. Enlarge it to any dg enhancement $\widetilde{\der}^b \coh(\EuT_0)$ of the bounded derived category $\der^b \coh(\EuT_0)$ (over $\Z$). That is, $\widetilde{\der} ^b\coh(\EuT_0)$ is a dg category containing $\tw \VB(\EuT_0)$ as a full subcategory, with an equivalence of triangulated categories $H^0(\widetilde{\der}^b \coh(\EuT_0))\to \der^b\coh(\EuT_0)$ extending the canonical equivalence $H^0(\tw\VB(\EuT_0))\to \perf (\EuT_0)$. A standard method to construct such an enhancement would be to use injective resolutions for coherent sheaves, and then to exhibit equivalence of that approach to the approach via \v{C}ech complexes in the case of perfect complexes by the method of \cite[Lemma 5.1]{SeiQuartic}. However, there is no requirement for the enlargement to be of geometric origin.

\begin{thm}\label{psi wrap}
There is a $\Z$-linear $A_\infty$-functor $  \psi_{\mathsf{wrap}} \colon \EuW(T_0)\to \widetilde{\der}\coh(\EuT_0)$ which extends to a quasi-equivalence
\[  \tw \EuW(T_0) \to \widetilde{\der}^b \coh(\EuT_0), \]
and which restricts to $ \psi_0 \colon \EuF(T_0) \to \tw \VB(\EuT_0)$.
\end{thm}

\subsection{Rank 1 local systems} 

It will be helpful to make the quasi-equivalence $\psi_0$ more explicitly by describing its effect on Lagrangians with finite-rank local systems; these may be regarded as twisted complexes in $\EuF(T_0)^\exact$.  

As usual, let $L_\infty^\# \in \ob \EuF(T_0)$ be a Lagrangian brane of slope $(0,-1)$ with non-trivial double covering. For any $h\in \C^*$, let $L^h_\infty$ denote the brane $L_\infty$ equipped with a rank 1 $\BbK$-local system of holonomy $h\in \BbK^\times$. This is an object in a larger Fukaya category $\EuF(T_0)^{\exact}_{\mathsf{loc}}$ whose objects are exact Lagrangian branes with local systems of finite rank free $\Z$-modules.

For a $\BbK$-linear $A_\infty$-category $\EuC$, let $\modules \EuC$ denote the category of finitely generated projective $\EuC$-modules, assigning to each object a finite cochain complex of finitely generated projective $\BbK$-modules. With $\BbK=\Z$, the object $L^h_\infty$ defines a left Yoneda-module $\EuY_L(L^h_\infty) = \hom_{\EuF(T_0)}(L^h_\infty,\cdot) \in \ob \modules \EuF(T_0)$. 

Let $\phi_0\colon \tw \VB(\EuT_0) \to \EuF(T_0)$ be an $A_\infty$ functor quasi-inverse to $\psi_0$. Over $\Z$, the existence of such a functor is not quite trivial. However, our earlier analysis of \v{C}ech complexes implies that one can define a quasi-inverse (or even strict inverse) $\EuB_{\EuT_0}\to \EuA$ to $\psi_0|_\EuA\colon \EuA\to \EuB_{\EuT_0}$. We then define $\phi_0$ by extending the latter functor to twisted complexes. Module categories are contravariant, and so $\phi_0$ induces a functor
\[ \phi_0^* \colon \modules \EuF(T_0)\to \modules \VB (\EuT_0).\]

\begin{lem}\label{where local systems go}
Work over a base ring $R$ which is a commutative, unital, normal, noetherian domain. Identify the normalization of $\EuT_0$ with $\mathbb{P}^1(R)$ by sending the preimages of the nodal section of $\EuT_0$ to 
$\{ [0:1], [1:0]\}$ and $\sigma$ to $[1:1]$. Two such identifications exist, of which one has the property that for each $h\in R^\times$, the module $P^h:=\phi_0^* (L^h_\infty) \in \ob \modules \VB(\EuT_0)$ is represented by a locally free resolution of the skyscraper sheaf $\EuO_h$ at the section $h=[h:1] \colon \spec R\to \mathbb{P}^1(R)$. 
\end{lem}

\begin{pf}
There is, for each $h\in R^\times$, a rank 1 $\BbK$-local system $\EuE^h$ over $T_0$ for which $\hol(L_\infty)=h$ and $\hol(L_0)=1$. Moreover, $\EuE^{h_1} \otimes  \EuE^{h_2} = \EuE^{h_1h_2}$. The local system $\EuE^h$ induces a strict autoequivalence $\alpha^h$ of $\EuF(T_0)^{\exact}_{\mathsf{loc}}$: on objects: leave the Lagrangian brane unchanged but tensor the local system by the restriction of $\EuE^h$. On morphism-spaces, for each intersection point $x \in L \cap L'$, map $x$ to $\theta(h) x$, where $\theta\colon \BbK \to \End_{\BbK}\EuE^h(x)$ is the isomorphism that sends $1$ to $\id$. 

Define $\EuA^h=\alpha^h(\EuA)$, the full subcategory of $\EuF(T_0)^{\exact}_{\mathsf{loc}}$ on the two objects $L_0$ with its trivial local system and $L_\infty^h$ = $L_\infty$ with its local system of holonomy $h$. Then $\EuA^h$ is a minimal $A_\infty$-structure.  Moreover, $\alpha^h$ induces a trace-preserving isomorphism $H^*(\EuA)\cong H^*(\EuA^h)$. By Theorem \ref{Ainf comparison}, the $A_\infty$-structure $\EuA^h$ is gauge-equivalent to $\EuB_C$ for a unique Weierstrass curve $C$. As an abstract curve, we have $C =  \EuT_0$, but the Weierstrass data (basepoint, differential) might not be standard. We can think of $C$ as $\EuT_0$ with standard differential, but different basepoint $\sigma(h)$.  By Theorem \ref{Ainf comparison}, the isomorphism $\alpha^h\colon  \EuA \to \EuA^h$ arises from a Weierstrass isomorphism $(\EuT_0, \omega, \sigma)\to (\EuT_0, \omega, \sigma(h))$. Automorphisms of $(\EuT_0,\omega)$ are the same thing as automorphisms of $\mathbb{P}^1$ that map $\{ 0,\infty\}$ to $\{0,\infty\}$; thus they form a group $(\Z/2)\ltimes R^\times$, where $\Z/2$ acts as the antipodal involution. Our construction gives rise to a homomorphism 
\[ \beta_R \colon  R^\times \to \aut \EuT_0  = (\Z/2)\ltimes R^\times\]
mapping $h$ to the automorphism $\beta^h$ such that $(\phi_0)^*(L^h_\infty)$ is represented by $\beta_R(h)^*\EuO_\sigma = \EuO_{\beta_R(h)\circ \sigma}$. 
The homomorphisms $\beta_R$ are by construction compatible with base change $R\to R'$. We claim that they must map $R^\times$ to $R^\times$. Indeed, since $R^\times$ is normal in $(\Z/2)\ltimes R^\times $, there is a quotient map $\bar{\beta}_R\colon R^\times \to \Z/2$, also natural in $R$. One must have $\bar{\beta}_k=0$ when $k$ is an algebraically closed field, since then every $z\in k^\times$ is a square. By naturality, $\bar{\beta}_k=0$ for arbitrary fields $k$ (embed $k$ into an algebraic closure), and so for arbitrary domains $R$ of the sort specified in the statement (embed $R$ into its field of fractions).

In view of their compatibility with base change, we view the $\beta_R$ collectively as a natural transformation $  \beta \colon \mathbb{G}_m \Rightarrow \mathbb{G}_m,$
where $\mathbf{rings}$ is the category of rings satisfying the conditions listed in the statement, and $\mathbb{G}_m\colon \mathbf{rings} \to \mathbf{groups}$ is the multiplicative group functor: $\mathbb{G}_m(R)=R^\times$.

Now, $\mathbb{G}_m$ is co-represented by $S: = \Z[t,t^{-1}]$, meaning that $\mathbb{G}_m \cong \Hom_{\mathbf{rings}}(S,\cdot)$---the identification makes the set $\Hom_{\mathbf{rings}}(S,R)$ into a group. The isomorphism sends $r \in R^\times $ to the homomorphism $S\to R$, $f \mapsto f(r)$. We note that $S$ is indeed a noetherian normal domain! By Yoneda's lemma, $\beta$ must arise from some ring endomorphism $b\colon S\to S$. Moreover, $\Hom_{\mathbf{rings}}(S,S) \cong  S^\times = \{ \pm t^d: d\in \Z\}$. These endomorphisms give rise to the set-theoretic natural transformations $R^\times \owns r\mapsto \pm r^d  \in R^\times$, of which $r\mapsto r^d$ is a group homomorphism but $r\mapsto - r^d$ is not. We note also that $\beta_R$ must be injective for each $R$, which leaves us only with the two possibilities  $b(r)= r^{\pm 1}$. One of these possibilities is the right one, the other not; this is the ambiguity left in the statement of the lemma.
\end{pf}

\begin{rmk}
The imprecision in the previous lemma is easily resolved. To fix an identification of the normalization of $\EuT_0$ with $\mathbb{P}^1$ of the sort described, it suffices to describe its effect on the tangent space to 
$\EuT_0$ at $\sigma$. Now, $\sigma^*T\EuT_0\cong R$ canonically, via the 1-form $\omega$; and $T_{[1:1]}\mathbb{P}^1 = T_{1}\mathbb{A}^1 = R$. The correct identification is the one given in these terms by $\id_R$. That is because $\sigma^*T\EuT_0$ can be understood as $\ext^1(\EuO_\sigma,\EuO_\sigma)$, which is identified by $\psi$ with $HF^1(L_\infty^\#,L_\infty^\#)$. The latter identification is the one determined by the trace-maps. 
\end{rmk}

\subsection{Higher rank local systems}
The material in this subsection is not used elsewhere in the paper.

In Lemma \ref{where local systems go}, we established that over a base ring $R$, the functor $\phi_0^*\colon \mod \EuF(T_0)^{\exact} \to \mod \VB \EuT_0$ maps $L_\infty^h$, that is, $L_\infty^\#$ with a rank 1 local system of holonomy $h\in \mathbb{G}_{\mathsf{m}}(R)$, to the skyscraper sheaf located at the section $h$ of the smooth locus $V= \EuT_0^{\mathsf{sm}} = \spec R[t,t^{-1}] = \mathbb{G}_{\mathsf{m}}(R)$. Over $\Z$, rank 1 local systems are not very interesting, higher rank local systems more so. Working over $\Z$, take a local system on $L_\infty^\#$ with fiber $\Z^n$ and holonomy $\phi\in GL(\Z^n)$;  denote this object by $L_\infty^\phi\in \ob \EuF(T_0)^{\exact}_{\mathsf{loc}}$. It maps under $\phi_0^*$ to a module for $\VB(\EuT_0)$ co-represented by some quasi-coherent complex $\EuK_\phi^\bullet$. That is,  $\EuY_L(\EuK_\phi^\bullet) \cong \phi_0^*\EuY_L(L_\infty^\phi) \in H^0(\modules \VB(\EuT_0))$. 

\begin{thm}
Let $V$ be the smooth locus in $\EuT_0$; thus $V \cong \spec \Z [t,t^{-1}] \cong \mathbb{G}_{\mathsf{m}}(\Z)$. For $\phi\in GL(\Z^n)$, Let $\Z^n_\phi$ denote the $\Z[t,t^{-1}]$-module $\Z^n$, on which $t$ acts as $\phi$. Let $(\Z^n_\phi)^\sim$ be the associated quasi-coherent sheaf on $V$, and let $K_\phi$ be the push-forward of $(\Z^n_\phi)^\sim$ to $\EuT_0$. Then $\EuK_\phi^\bullet$ is quasi-isomorphic to the sheaf $K_\phi$.
\end{thm}

\begin{pf}
Notice that $L_\infty^\phi$ is quasi-isomorphic in $\twsplit \EuF(T_0)^{\exact}_{\mathsf{loc}}$ to an object of $\twsplit \EuA$. As such, it is a \emph{compact} object of the dg category of $\EuA$-modules (see for instance \cite{BFN}). Quasi-equivalences preserve compact objects; consequently $\EuK^\bullet_\phi$ is compact as an object of $\QC(\EuT_0)$. Therefore $\EuK^\bullet_\phi$ is quasi-isomorphic to a perfect complex \cite{Nee,BFN}. Since only its quasi-isomorphism class matters, we may assume that $\EuK^\bullet_\phi$ is a strictly perfect complex---a finite complex of locally free sheaves.

It will be helpful to work over base \emph{fields} $\BbK$. We then take the object $L_\infty^\phi \in \EuF(T_0)^{\exact}_{\mathsf{loc}}\otimes \BbK$ associated with 
$\phi \in GL_n(\BbK)$. 

\emph{Step 1}. Work over an algebraically closed field $\BbK$. We claim that $\EuK^\bullet_\phi$ is then quasi-isomorphic to a sheaf supported whose support is contained in the eigenvalue spectrum  $\mathsf{eval}\,\phi$ in $\BbK^* = \spec \BbK[t,t^{-1}]= V(\BbK)= V\times_{\spec \Z} \spec \BbK$. 

If $\lambda \in \BbK^*$, we have a hyperext spectral sequence
\[ E_2^{rs} =  \ext_{\EuT_0}^r(\EuH^{-s}(\EuK^\bullet_\phi), \EuO_\lambda) \Rightarrow \RHom^{r+s}_{\EuT_0}(\EuK^\bullet_\phi, \EuO_\lambda)
\cong HF^{r+s}(L^\phi_\infty, L^\lambda_\infty) .  \]
If $\lambda \in \mathsf{eval}\, \phi$ then (only) $HF^0$ and $HF^1$ are non-zero. If $\lambda \not \in \mathsf{eval}\, \phi$ then $HF^*=0$. One has $ \ext_{\EuT_0}^r(\EuH^{-s}(\EuK^\bullet_\phi), \EuO_\lambda)= \ext_V^r(\EuH^{-s}(j^*\EuK^\bullet_\phi), \EuO_\lambda)$, where $j\colon V\to \EuT_0$ is the inclusion. Since $V$ is the spectrum of a regular local ring of dimension 1, the $\ext^r$-modules vanish except for $r\in \{ 0 ,1 \}$.  Since it is supported in two adjacent columns, the spectral sequence degenerates at $E_2$. Therefore $\Hom_V( \EuH^{-s}(j^*\EuK^\bullet_\phi), \EuO_\lambda)=0$ except when $-s \in \{0,1\}$ and $\lambda \in \spec \phi$. Hence $j^*\EuH^{-s}(j^*\EuK^\bullet_\phi) = 0$ for $s \notin \{-1,0\}$. 

Next, consider the hyperext spectral sequence
\[ 'E_2^{rs} = \ext_{\EuT_0}^r(\EuO, \EuH^{s}(\EuK^\bullet_\phi)) \Rightarrow \RHom^{r+s}_{\EuT_0}(\EuO,\EuK^\bullet_\phi) \cong \BbK^n. \]
The cohomology sheaf $\EuH^{s}(\EuK^\bullet_\phi)$ is supported on the singular section together with $\mathrm{eval}\, \phi$. One has 
$\ext_{\EuT_0}^r(\EuO, \EuH^{s}(\EuK^\bullet_\phi))= \h^r( \EuH^{s}(\EuK^\bullet_\phi))$ which is zero for $r\neq 0$ because $\EuH^{s}(\EuK^\bullet_\phi)$  has affine support. Thus $'E_2^{rs}=0$ except when $r=0$. So the spectral sequence degenerates, and we see that $\EuH^s (\EuK^\bullet_\phi)=0$ for $s\neq 0$. We may therefore truncate the complex $\EuK^\bullet_\phi$, replacing it by $0$th cohomology $\EuH^0(\EuK^\bullet_\phi)$, to which it is quasi-isomorphic. Further, $\EuH^0(\EuK^\bullet_\phi)$ is a torsion sheaf since its stalks are generically zero.

Notice that if we have a short exact sequence of modules $\Z^n_\phi$, the corresponding sheaves $\EuH^0(\EuK^\bullet_\phi)$ form a long exact triangle. From that, and the fact that over the algebraically closed field $\BbK$, matrices are conjugate to upper triangular matrices, we see that $\EuH^0(\EuK^\bullet_\phi)$ is supported in $V(\BbK)$, and therefore in $\mathsf{eval}\, \phi \subset V(\BbK)$.

\emph{Step 2.} Over an arbitrary field $\BbK$, $\EuK^\bullet_\phi$ is quasi-isomorphic to a torsion sheaf $K'_\phi$, supported in $V(\BbK)$, such that $\h^0(K'_\phi)\cong \BbK^n$ canonically. 

Indeed, by (flat) base change from $\BbK$ to its algebraic closure, we see that the $s$th cohomology sheaf of $\EuK^\bullet_\phi$ vanishes for each $s\neq 0$. By truncation we may replace the complex by its zeroth cohomology $\EuF:= \EuH^0(\EuK^\bullet_\phi)$. Moreover, $\EuF$ vanishes at the generic point, and its stalk at the singular section is zero, so it torsion and supported in $V\otimes_\Z \BbK$. Moreover, $\h^0(\EuF)\cong \BbK^n$ via the spectral sequence $'E_*^{**}$ above.

\emph{Step 3.} Over $\Z$, $\EuK^\bullet_\phi$ is quasi-isomorphic to a torsion sheaf supported in $V$. Moreover,  $\h^0(\EuF)\cong \Z^n$ canonically.

Since $\EuK^\bullet_\phi$ is a perfect complex, its cohomology sheaves are coherent. Hence $\h^0(\EuH^s(\EuK^\bullet_\phi))$ is a finitely generated $\Z$-module for each $s$. By Step 2 and the compatibility of the construction with base change, one has $\h^0(\EuH^s(\EuK^\bullet_\phi))\otimes\BbK=0$ for every field $\BbK$ and every $s\neq 0$. Hence $\h^0(\EuH^s(\EuK^\bullet_\phi))=0$ for $s\neq 0$. On the other hand, $z^*(\EuH^s(\EuK^\bullet_\phi))$ is a coherent sheaf on $\spec \Z$, for each $s$; or in other words, it is a finitely generated abelian group $G_\phi$.  By the projection formula, and Step 2, $G_\phi\otimes\BbK=0$ for any field $\BbK$. Therefore $G_\phi=0$. Consequently $\EuH^s(\EuK^\bullet_\phi)$ is supported in $V$, and so is the module associated with its sections $\h^0(\EuH^s(\EuK^\bullet_\phi))$, which is $0$ if $s\neq 0$. Hence $\EuK^\bullet_\phi$ is quasi-isomorphic to its $0$th cohomology sheaf. One has $\h^0(\EuK^\bullet_\phi)=\Z^n$, again via the spectral sequence $'E^*_{**}$.  

\emph{Step 4}. Completion of the proof.

Since it is torsion and supported in $V$, the sheaf $\EuH^0(\EuK^\bullet_\phi)$ is the push-forward of a coherent sheaf on $V$. We think of this as the sheaf associated with a module $M_\phi$. We stress that $M_\phi$ is canonically identified with $\mathbb{\Z}^n$ as an $\Z$-module, so $M_\phi = \mathbb{\Z}^n_{\phi'}$ for a well-defined matrix $\phi'\in GL_n(\BbK)$. Moroever, the map $\phi \mapsto \phi'$ is compatible with conjugation of matrices (i.e. $(\chi \phi \chi^{-1})' = \chi \phi ' \chi^{-1}$).

These points apply over any $\Z$-algebra $\BbK$, and the construction is compatible with base change. Hence the map $\phi\to \phi'$ arises from a map of $\Z$-schemes $F\colon GL_n(\Z)\to GL_n(\Z)$. We claim that $F=\id$. It will suffice to show that $F =\id$ when we base-change to $\overline{\Q}$, an algebraic closure of the rationals. The induced map $F \colon GL_n(\overline{\Q})\to GL_n(\overline{\Q})$ is the identity on the diagonal matrices, and therefore, by compatibility with conjugation, also on the Zariski-open set of diagonalizable matrices. Therefore it is the identity map.
\end{pf}

\subsection{Generation of the wrapped Fukaya category and the bounded derived category}

Let $\Lambda^\# \in \ob \EuW(T_0)$ be an arc of slope $(0,-1)$, graded so that $HW(L_0,\Lambda^\#)$ lies in degree $0$, and oriented so that it runs into $z$.  Our functor $\psi_{\mathsf{wrap}}$ will carry $\Lambda^\#$ to $\EuO_s$, the skyscraper sheaf along the nodal section.

\begin{lem}
$\EuW(T_0)$ is generated by $\EuF(T_0)$ and $\Lambda^\#$.
\end{lem}
\begin{pf}
First, if orient $\Lambda$ in the opposite direction, we obtain an isomorphic object of $\EuW(T_0)$. This reflects the fact that a spin-structure on a Lagrangian (which trivializes $w_1$ and $w_2$) is more data than is needed; a Pin-structure, trivializing $w_2$, is sufficient \cite{SeiBook}.  Changing the orientation corresponds to an automorphism $(-1)^{\deg} \id$ of the object. The brane structure also involves a double covering of $\Lambda$, but that is necessarily trivial. With these points noted, we find that any object of $\EuW(T_0)$ whose Lagrangian is non-compact is quasi-isomorphic to a shift of an iterated Dehn twist of $\Lambda^\#$ along closed, exact curves equipped with non-trivial double coverings. These Dehn twists act on $\tw \EuW(T_0)$ by spherical twists \cite{SeiBook, LP}; hence the arcs are represented by twisted complexes in $\EuF(T_0)$ and $\Lambda^\#$.
\end{pf}

\begin{lem}
$\widetilde{\der}\coh(\EuT_0)$ is generated by $\VB (\EuT_0)$ and the skyscraper sheaf $\EuO_s$ along the nodal section $s$.
\end{lem}
\begin{pf}
It suffices to show that $G_0(\EuT_0)$ is generated, as an abelian group, by the image of $K_0(\EuT_0)$ and the class $[\EuO_s]$. Here $G_0$ denotes the Grothendieck group of coherent sheaves, while $K_0$ is the Grothendieck group of vector bundles. The quotient $G_0(\EuT_0)/\im K_0(\EuT_0)$ is certainly generated by coherent sheaves supported on $Z$, the closure of $\im z$. Thus it will suffice to show that $K_0(M_Z(\EuT_0))$, the Grothendieck group of the abelian category of coherent sheaves supported along $Z$, is generated by the class of $\EuO_s$. We now proceed as in \cite[ex. II  6.3.4]{WeiK}: $M_Z(\EuT_0)$ is the abelian category of finitely-generated modules $M$ for $R=\Z[x,y]/(y^2+xy-x^3)$ such that $I^n M =0$ for some $n$, where $I=(x,y)\subset R$. Such a module has a filtration $M\supset IM \supset I^2 M \supset \dots \supset I^n M = 0$ and therefore $K_0(M_Z(\EuT_0))$ is generated by the factors of such filtrations, i.e., by modules $N$ with $IN=0$. As sheaves, those are precisely the push-forwards of sheaves on $Z$, or equivalently those of the form $s_*\EuF$ for some coherent sheaf $\EuF$ on $\spec \Z$. Since $G_0(Z)=\Z$, the result follows.
\end{pf}

\subsection{\texorpdfstring{The functor $\psi_{\mathsf{wrap}}$}{The functor psi-wrap}}
We have a sequence of $A_\infty$-functors
\begin{equation}\label{functor sequence}
   \EuW(T_0) \xrightarrow{\EuY_L} \modules \EuA \xrightarrow{\phi_0^*} \modules\EuB_{\EuT_0}  \xleftarrow{\simeq}\modules \VB(\EuT_0) \xleftarrow{\simeq} \QC(\EuT_0).    
   \end{equation}
Here $\EuY_L$ is the (covariant) left Yoneda functor, $X\mapsto \hom_\EuA(X,\cdot)$. The functor  $\modules\EuB_{\EuT_0} \xleftarrow{\simeq}\modules \VB(\EuT_0)$ is restriction, which is a quasi-equivalence because $\EuB_C$ split-generates $\tw \VB(\EuT_0)$. We denote by $\QC(\EuT_0)$ a dg enhancement of the unbounded derived category of quasi-coherent complexes, as in \cite{Toe}; the last map, $\modules \VB(\EuT_0)\xleftarrow{\simeq} \QC(\EuT_0)$ is again a left Yoneda functor. 

It follows from \cite[Theorem 8.9]{Toe} or \cite{BFN} that $\modules \VB(\EuT_0)\xleftarrow{\simeq} \QC(\EuT_0)$ is a quasi-isomorphism.  Indeed, the dg category of modules over $\VB(\EuT_0)$ is, rather trivially, equivalent to the dg category of dg functors $[\tw \VB(\EuT_0), \tw \VB(\spec \Z)]$. There are further equivalences $[\tw \VB(\EuT_0), \tw \VB(\spec \Z)]\to [\QC(\EuT_0), \QC(\spec \Z)]_c$ where the $c$ denotes that the functors respect filtered colimits (think of quasi-coherent complexes as filtered colimits of perfect complexes); and $\QC(\EuT_0\times_{\spec \Z}\spec \Z)=\QC(\EuT_0) \to [\QC(\EuT_0), \QC(\spec \Z)]_c$ (mapping an integral kernel to its push-pull functor).

The conclusion of this discussion is as follows: 
\begin{lem}
There is a homotopy-commutative diagram of $A_\infty$-functors
\[\xymatrix{ 
\EuW(T_0) \ar[r] & \modules \EuB_{\EuT_0} & \QC(\EuT_0)\ar[l] \\
\EuF(T_0)^\exact\ar[u] \ar[rr]^{\simeq}   &&\tw \VB(\EuT_0)\ar[u]
} \]
\end{lem}

We now define $\EuZ^\bullet$ to be an object of $\QC(\EuT_0)$ corresponding under mirror symmetry to the arc $\Lambda^\#$. Precisely:
\begin{leftbar}
Define $\EuZ^\bullet$ to be a choice of quasi-coherent complex whose associated $\EuB_{\EuT_0}$-module $\EuY_L(\EuZ^\bullet)=\hom_{\EuB_{\EuT_0}}(\EuZ^\bullet,\cdot)$ is quasi-isomorphic to $\phi_0^*\EuY_L(\Lambda^\#)$. 
\end{leftbar}

We have the full subcategory  $\langle \EuZ^\bullet, \tw \VB(\EuT_0) \rangle \subset \QC(\EuT_0)$, and, by restricting the arrow $\leftarrow$ in the previous diagram, a functor-sequence
\[\EuW(T_0) \xrightarrow{r} \modules \EuB_{\EuT_0} \xleftarrow{s}  \langle \EuZ^\bullet, \tw \VB(\EuT_0)\rangle.   \]
The arrow $H^0(s)$ is an embedding, and the composite
\[ H^0(\psi_{\mathsf{wrap}}) := H^0(s)^{-1} \circ H^0(r) \colon H^0(\EuW(T_0)) \to H^0 \langle \EuZ^\bullet, \tw \VB(\EuT_0)\rangle \]
is well-defined. 

Similarly, we have quasi-embeddings $\QC(\EuT_0) \to \modules \VB(\EuT_0)\to  \modules\EuB_{\EuT_0}  \xrightarrow{\psi_0^*} \modules \EuA$ and a functor-sequence
\[   	 \langle \EuZ^\bullet,  \tw \VB(\EuT_0)\rangle \to  \modules \EuA  \leftarrow \EuW(T_0)  \]
whose effect on cohomology lifts to a functor
\[ H^0(\phi_{\mathsf{wrap}}) \colon H^0 \langle \EuZ^\bullet,  \tw \VB(\EuT_0)\rangle \to H^0(\EuW(T_0))\] 
mapping $\EuZ^\bullet$ to $\Lambda^\#$. The composite $H^0(\psi_{\mathsf{wrap}})\circ H^0(\phi_{\mathsf{wrap}})$ is the identity functor on $H^0\langle \EuZ^\bullet,  \tw \VB(\EuT_0)\rangle$ (this composite does not involve the potentially information-losing functor $\EuW(T_0)\to \modules \EuF(T_0)$). Hence
\begin{lem}
 $H^0(\psi_{\mathsf{wrap}})$ is full.
\end{lem}

A key point in the proof of mirror symmetry for the wrapped category will be the following assertion:

\begin{prop}\label{where the arc goes}
One has $\EuZ^\bullet \simeq \EuO_s$ in $\QC(\EuT_0)$, where as before, $\EuO_s$ is the push-forward of $\EuO_{\spec \Z}$ by the singular section $s \colon \spec \Z \to \EuT_0$.
\end{prop}

Taking the proposition for granted for the moment, we explain how to complete the proof of Theorem \ref{psi wrap}. First, it implies that $ \langle \EuZ^\bullet,  \tw \VB(\EuT_0)\rangle  \simeq \widetilde{\der}^b \coh(\EuT_0) \subset \QC(\EuT_0)$, where $\widetilde{\der}^b \coh(\EuT_0)$ is the dg category of bounded complexes with coherent cohomology. We now have a  functor sequence
\[  \EuW(T_0) \to \modules \EuB_{\EuT_0} \leftarrow  \widetilde{\der}^b \coh(\EuT_0) \]
Now, one can set up the full subcategory $\{L_0^\#, L_\infty^\# , \Lambda^\#\} \subset \EuW(T_0)$ so as to be a minimal $A_\infty$-category. (In particular, in the case of endomorphisms of $\Lambda^\#$, one can easily draw a perturbation $\Lambda'$ of $\Lambda$ so that they do not jointly bound any immersed bigons.) Hence there is a map $\EuW(T_0) \to H^0(\EuW(T_0))$ inducing the identity map on cohomology. The full subcategory $\{ \EuO,\EuO_\sigma, \EuO_\sigma \} \subset \widetilde{\der}^b \coh(\EuT_0)$ has projective hom-spaces, and thus there is a map $\der\coh(\EuT_0)  \to \widetilde{\der}^b \coh(\EuT_0)$ inducing the identity on cohomology. Using homological perturbation theory, we obtain a functor
\[ \psi_{\mathsf{wrap}} \colon \EuW(T_0) \to \widetilde{\der}^b \coh (\EuT_0) \] 
such that $H^0(\psi_{\mathsf{wrap}})$ is the functor previously so-denoted.
 
We must show that $H^0(\psi_{\mathsf{wrap}})$ is faithful. To do so, it will suffice to prove it on the hom-spaces $HW^*(\Lambda,X)$ and $HW^*(X,\Lambda)$, where $X$ runs through a split-generating set: $X=L_0^\#$ or $L_\infty^\#$ or $\Lambda$. The graded $\Z$-modules $HW^*(\Lambda,X)$ are in each case isomorphic to their mirrors $\ext^*(\EuO_s, \psi_{\mathsf{wrap}}(X))$; likewise $HW^*(X,\Lambda) \cong \ext^*(\psi_{\mathsf{wrap}}(X), \EuO_s)$. They are free abelian graded groups, of finite rank in each degree. The map $\psi_{\mathsf{wrap}}\colon  HW^*(\Lambda,X)\to \ext^*(\EuO_s, \psi_{\mathsf{wrap}}(X))$ is surjective and is therefore an isomorphism---and the same goes for $HW^*(X,\Lambda)$. This completes the proof of Theorem \ref{psi wrap}, modulo that of Proposition \ref{where the arc goes}. 

\begin{lem}
If we work over a base \emph{field} $\BbK$, then the quasi-coherent complex $\EuZ^\bullet_\BbK$ mirror to the arc is quasi-isomorphic to $\EuO_s$, the skyscraper at the nodal point $s$ of the $\BbK$-variety $\EuT_0(\BbK)$.
\end{lem}
\begin{pf}

Let $\EuH^k_\BbK$ be the $k$th cohomology sheaf of $\EuZ^\bullet_\BbK$. In general, if $\EuE^\bullet$ and $\EuF^\bullet$ are complexes of quasi-coherent sheaves over a scheme $X$, at least one of them bounded, one has a right half-plane spectral sequence arising from a filtration on a complex computing $\RHom_{\EuO_X}^{r+s}(\EuE^\bullet,\EuF)$, with
\begin{equation} \label{hyperext ss}
E_2^{rs} = \bigoplus_k { \ext^r_{\EuO_X}\left(\EuH^{k-s}(\EuE^{\bullet}), \EuH^k(\EuF^\bullet)\right) },
\end{equation} 
abutting to $\RHom_{\EuO_X}^{r+s}(\EuE^\bullet,\EuF^\bullet)$. In particular, we have a spectral sequence
\begin{equation} \label{hyperext at p}
E_2^{rs} = \ext^r_{\EuT_0}(\EuO_p, \EuH^{s}_\BbK) \quad \Rightarrow\quad \RHom_{\EuT_0}^{r+s}(\EuO_\sigma, \EuZ^\bullet_\BbK).
\end{equation}
where $p=[0:1:0]$. Since $H^0(\psi_{\mathsf{wrap}})$ is full and $HW(L_\infty^\#,\Lambda^\#)=0$, we have $\RHom_{\EuT_0}(\EuO_\sigma, \EuZ^\bullet_\BbK) = 0$.
This spectral sequence degenerates at $E_2$. Indeed, let $M^s$ denote the module of global sections of $j^*\EuH^s_\BbK$, where $j$ is  the open inclusion of the smooth locus $V(\BbK)\cong \spec \BbK[t,t^{-1}]$. Let $\mathfrak{p}=(t-1)\BbK[t,t^{-1}]$, a maximal ideal of $\BbK[t,t^{-1}]$, and let $A$ denote the localization of $\BbK[t,t^{-1}]$ at $\mathfrak{p}$. Then  $\ext^r_{\EuT_0}(\EuO_\sigma, \EuH^{s}_\BbK)\cong \ext^r_A (\BbK, \EuH^s_{\mathfrak{p}})$. Since $A$ is a regular ring of dimension 1, these Ext-modules vanish except for $r=0$ or $1$. Therefore the spectral sequence is concentrated in two adjacent columns, and so degenerates. Thus $\Hom_{A}(\BbK, \EuH^s_{\mathfrak{p}})=0$, and so $\EuH^s_{\mathfrak{p}}=0$. 

Consequently, $\EuH^k_\BbK$ is for each $k$ supported in the open set $U(\BbK)=\spec \BbK[x,y]/(y^2+xy-x^3)$. Next consider the spectral sequence
\begin{equation} 
'E_2^{rs} = \ext^r_{\EuT_0}(\EuO, \EuH^{s}_\BbK) \quad \Rightarrow\quad \RHom_{\EuT_0}^{r+s}(\EuO, \EuZ^\bullet_\BbK) \cong HF(L_0^\#, \Lambda^\#;\BbK )\cong \BbK.
\end{equation} 
Since $HF(L_0^\#, \Lambda^\#;\BbK )\cong \BbK$, one has $\RHom_{\EuT_0}^{k}(\EuO, \EuZ^\bullet_\BbK)=0$ for $k\neq 0$, while the space $\Hom_{\EuT_0}(\EuO, \EuZ^\bullet_\BbK)$ is at most 1-dimensional.
One has $\ext^r_{\EuT_0}(\EuO, \EuH^{s}_\BbK)=\h^r(\EuH^s_\BbK)$, but this cohomology module vanishes for each $r>0$. This spectral sequence therefore also degenerates, and so $\h^0(\EuH^s_\BbK)$ vanishes for all $s\neq 0$. Since $\EuH^s_\BbK$ is quasi-coherent with affine support in $U(\BbK)$, it follows that $\EuH^s_\BbK=0$. 

Because of the vanishing of the non-zero cohomology sheaves, we have a diagram of quasi-isomorphisms $\EuZ^\bullet \leftarrow \tau_{\leq 0}\EuZ^\bullet \to \EuH^0$, where $\tau_{\leq 0}$ denotes the truncation $\dots\to \EuZ^{-2} \to \EuZ^{-1}\to \ker \delta^0\to 0$. Therefore $\EuZ^\bullet_\BbK$ is quasi-isomorphic to the sheaf $\EuH^0_\BbK$, with which we may replace it.

Moreover, $\EuH^0_\BbK$ likewise has affine support in $U(\BbK)$, and so it is either zero, or is the push-forward from $U(\BbK)$ of the sheaf associated with the $\BbK[x,y]/(y^2+xy-x^3)$-module $\BbK$ (with some action of $x$ and $y$).  Thus $\EuH^0_\BbK$ is the skyscraper sheaf $\EuO_x$ at a closed point $x$ of $\EuT_0(\BbK)$. Now, if $x$ is a regular point, then we have by Lemma \ref{where local systems go} $\EuO_x \simeq \psi_0(L_\infty^h)$ for some $h\in \BbK^\times$, and thus $H^0(\psi_{\mathsf{wrap}})$ maps $\Lambda^\#$ and $L_\infty^h$ to isomorphic objects of $D\coh (\EuT_0)$. It is easy to check that the Yoneda functor $H^0(\EuW(T_0))\to H^0(\modules \EuF^\exact(T_0))$ reflects isomorphism, and sends no object to the zero object. Since it is a composite of that Yoneda functor and an embedding, $H^0(\psi_{\mathsf{wrap}})$ again reflects isomorphism. Thus it cannot map $L_\infty^h$ and $\Lambda^\#$ to isomorphic objects; nor can it map $\Lambda^\#$ to the zero-object. So $x$ is the unique singular point $s$.  
\end{pf}

\begin{pf}[Proof of Prop. \ref{where the arc goes}]
Again, let $j^*\EuH^k$ be the restriction of $\EuH^k$ to the smooth locus $j\colon V\to \EuT_0$. One has $V\cong \spec R$ where $R=\Z[t,t^{-1}]$, an identification under which the section $\sigma$ is defined by the homomorphism $R \to \Z$, $p(t)\mapsto p(1)$. Let $M^k$ be the $R$-module of sections $\Gamma(j^*\EuH^k , V)$. Since $\EuH^k$ is quasi-coherent, $j^*\EuH^k$ is the sheaf associated with $M^k$. Let $\mathfrak{p}=(t-1)R$, and let $R_{\mathfrak{p}}$ be the localization at $\mathfrak{p}$. Our task is to show that the $R_{\mathfrak{p}}$-module $M^k_{\mathfrak{p}}$ is zero. 

We observe that $M^k\otimes_\Z \Q = 0$. For this, consider the wrapped Fukaya category with $\Q$-coefficients; the arc now maps to $\EuZ^\bullet_\Q=\EuZ^\bullet \otimes_\Z \Q$. One has $\EuH^k(\EuZ^\bullet \otimes \Q) = \EuH^\bullet \otimes \Q$. Also, $j^*(\EuH^k(\EuZ^\bullet \otimes \Q)) = j^*\EuH^k \otimes \Q$. Thus, $j^*(\EuH^k(\EuZ^\bullet \otimes \Q))$ is the sheaf associated with its module of sections $M^k\otimes_\Z\Q$ over $R\otimes_\Z \Q$. The previous lemma then implies that $M^k\otimes_\Z\Q=0$.

Now, $R_{\mathfrak{p}}=S^{-1}R$, where $S=R\setminus \mathfrak{p}$; in particular, $S$ contains every prime of $\Z$, and hence $R_{\mathfrak{p}}$ is a $\Q$-algebra. Thus $M^k_{\mathfrak{p}}=0$ is a $\Q$-algebra, yet is torsion as an abelian group. Hence $M^k_{\mathfrak{p}}=0$. 

We now proceed on the same lines as the argument for the previous lemma. We now know that $\EuH^k$ is supported in $U=\spec \Z[x,y]/(y^2+xy-x^3)$. We deduce, just as before, that $\EuH^k=0$ for $k\neq 0$, and hence that $\EuZ^\bullet$ is quasi-isomorphic to the cohomology sheaf $\EuH^0$; and that $\EuH^0$ corresponds to the $\EuO_U$-module $\Z$, with some action of $x$ and $y$;  but $x$ and $y$ necessarily act as zero, because otherwise $\EuH^0$ would be the image under $H^0(\psi_{\mathsf{wrap}})$ of a skyscraper sheaf at a non-singular section, contradicting the fact that $H^0(\psi_{\mathsf{wrap}})$ reflects isomorphism. 
Therefore $\EuH^0=s_*\EuO_{\spec \Z}$.
\end{pf}

\section{\texorpdfstring{Identifying the mirror curve over $\Z\series{q}$: the Tate curve and toric geometry}{Identifying the mirror curve: the Tate curve and toric geometry}}
In this section, we offer a third proof that $C_{\mathsf{mirror}}|_{q=0}$ has the equation $y^2+xy=x^3$. More significantly, this proof extends to show that $C_{\mathsf{mirror}}$ is the Tate curve.

\subsection{The Tate curve}\label{Tate section}
The Tate curve was constructed by Raynaud using formal schemes; we have followed the expositions by Deligne--Rapoport \cite{DR} and Conrad \cite{Con}, and Gross's  reinterpretation of the construction in toric language \cite{Gro}. We review the construction.

\subsubsection{\texorpdfstring{Construction of a scheme with $\Z$-action $\EuT_\infty\to \spec \Z[t]$}{Construction of a scheme with Z-action}}

\begin{figure}[ht!]\label{fan}
\centering
\includegraphics[width=0.8\textwidth]{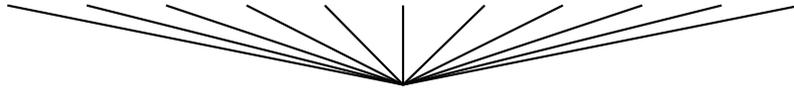}
\caption{Part of the fan $F$}\end{figure}

\paragraph{Toric construction of $\EuT_\infty$.} We begin with the toric fan picture. Fix a commutative ring $R$. Consider the rays $\rho_i = \Q_+(i,1)\subset \Q^2$, where $i\in \Z$ (Figure \ref{fan}). 
The convex hull of $\rho_i$ and $\rho_{i+1}$ is a cone $\sigma_{i+1/2}\subset \Q^2$; the collection of cones $\sigma_{i+1/2}$, their boundary faces $\rho_j$, and \emph{their} common endpoint $\{0\}$, form a rational fan $F$ in $\Q^2$. Each cone $c$ of $F$ has a dual cone 
\[ c^\vee= \{ x\in  \Q^2: \langle x ,s \rangle \geq 0 \; \forall s\in c \}  \]
(the use of inner products rather than dual spaces is an aid to visualization). To be explicit, $\rho_i^\vee$ is the half-plane $\{ (m,n)\in \Q^2: i m+n\geq 0 \}$, and $\sigma_{i+1/2}^\vee = \rho_i^\vee \cap \rho_{i+1}^\vee$. The dual cones $c^\vee$ have semigroup $R$-algebras $R[c^\vee\cap \Z^2]$ spanned by monomials $e^{\lambda}$ where $\lambda \in c^\vee\cap \Z^2$. We glue together the affine schemes $U_{i+1/2} = \spec R[\sigma_{i+1/2}^\vee\cap \Z^2]$ so as to form a toric scheme $\EuT_\infty\to \spec R$ with an affine open cover $\{U_{i+1/2}\}$.

For each $i$, there is a map of cones $\sigma_{i+1/2} \to \Q_+$, $(x_1,x_2)\mapsto x_2$. The dual to this map is the map of cones $\Q_+\to \sigma_{i+1/2}^\vee$ given by $1 \mapsto (0,1)$, which induces a map of semigroup $R$-algebras $R[\Q^+\cap \Z]\to R[\sigma^\vee_{i+1/2}\cap \Z^2]$ and hence a function $U_{i+1/2}\to \mathbb{A}^1(R)$. The maps $\sigma_{i+1/2} \to \Q_+$ assemble to form a map of fans $F \to \Q_+$, whence the functions $U_{i+1/2} \to \mathbb{A}^1(R)$ consistently define a morphism $\EuT_\infty \to \mathbb{A}^1(R)=\spec R[t]$. \emph{This morphism makes  $\EuT_\infty$ an $R[t]$-scheme.} We can regard $t$ as a regular function on $\EuT_\infty$; its restriction to $U_{i+1/2}$ is the monomial $t=e^{(0,1)}$.  

\paragraph{Explicit description of the gluing maps.}
We now describe the gluing construction of $\EuT_\infty$ in scheme-theoretic terms, without the toric language. The exposition follows Deligne--Rapoport's (\emph{op. cit.}); we recall it for convenience. The lattice points $\sigma_{i+1/2}^\vee \cap\Z^2$ are just the $\Z_{\geq 0}$-linear combinations of $(1, - i )$ and $(-1, i+1)$. In $R[\sigma_{i+1/2}^\vee]$, let $Y_{i+1}=e^{(-1,i+1)}$ and $X_{i}=e^{(1,-i)}$. Then
\[  U_{i+1/2} = \spec \frac{R[t][ X_{i}, Y_{i+1}]}{(X_i Y_{i+1} - t)}.\] 
Now let $V_i = U_{i-1/2}\cap U_{i+1/2}$ as subsets of $\EuT_\infty$. In abstract terms, we have
\begin{align*} 
& V_i \subset U_{i+1/2},  && V_i = U_{i+1/2}[X_i^{-1}] \stackrel{X_i\mapsto X_i}{\cong}  \spec R[t] [X_i,X_i^{-1}] \\
& V_i \subset U_{i-1/2},  && V_i = U_{i-1/2}[Y_i^{-1}]   \stackrel{Y_i\mapsto Y_i}{\cong}  \spec R[t] [Y_i^{-1}, Y_i], 
\end{align*}
and these two descriptions of $V_i$ are matched up by putting $X_i Y_i = 1$. The scheme obtained from the union of the open sets $U_{i+1/2}$ by gluing $U_{i-1/2}$ to $U_{i+1/2}$ is $\EuT_{\infty}$. No additional gluing is required because when $j-i>1$ one has $U_{i-1/2}\cap U_{j-1/2}\subset U_{i-1/2}\cap U_{i+1/2}$. 

In the scheme $\EuT_\infty[t^{-1}]\to \spec R[t,t^{-1}]$, all the open sets $U_{i+1/2}[t^{-1}]$ are identified with one another. Thus one has isomorphisms
$\EuT_\infty[t^{-1}] \cong U_{1/2}[t^{-1}]  \cong  \mathbb{G}_{\mathsf{m}}(R[t,t^{-1}])$.
The fiber $\EuT_\infty|_{t=0} = \EuT_\infty \times_{R[t]} R$ is the union of the subsets 
\[U_{i+1/2}|_{t=0} = \spec R[X_i,Y_{i+1}]/(X_iY_{i+1}).\] 

The sets $U_{i+1/2}|_{t=0}$ are disjoint when the indices $i$ and $j$ are not adjacent ($|i-j|>1$). Moreover, $U_{i-1/2}|_{t=0}$ is glued to $U_{i+1/2}|_{t=0}$ along $V_i|_{t=0}=R[X_i,X_i^{-1}]=R[Y_i,Y_i^{-1}]$. Thus $\EuT_\infty|_{t=0}$ is an infinite chain of $\mathbb{P}^1$'s, say $\mathbb{P}^1_i\subset U_{i-1/2}|_{t=0}\cap U_{i+1/2}|_{t=0}$.  

 The group $\Z$ acts on $\EuT_\infty$ covering the trivial action on $\spec R[t]$. The action is induced by a $\Z$-action on $\Q^2$ preserving the fan $F$, given by $n\cdot(x_1,x_2)=(x_1+n, x_2)$. One has $n\cdot U_{i+1/2} = U_{i+n+1/2}$; the action identifies $X_i$ with $X_{i+n}$ and $Y_i$ with $Y_{i+n}$.

\subsection{Quotients along the thickened zero-fiber} 
We wish to form quotients of $\EuT_\infty$ by the group $d\Z\subset \Z$, but find ourselves unable to do so in the category of $R[t]$-schemes. One can, however, construct quotients by $d\Z$ of thickened neighborhoods of the zero-fiber, i.e. of $\EuT_\infty|_{t^k  = 0}$, by virtue of the following fact:

\begin{leftbar}
The $U_{i+1/2}|_{t^k=0}$ form a chain: they are disjoint when $i$ and $j$ are not adjacent.
\end{leftbar}

Indeed, when  $j-i>1$, the set $U_{i-1/2}\cap U_{j-1/2}\subset \EuT_\infty$ lies over $\spec R[t,t^{-1}]$. 

We write down a concrete model for the quotient $\EuT_\infty|_{t^k=0}/(d\Z)$, when $d>1$. We take the affine schemes $U_{i+1/2}|_{t^k=0}$ for $i=0,\dots, d-1$, and we identify the open subset $V_d$ of $U_{d-1/2}|_{t^k=0}$ with the open subset $V_0$ of $U_{1/2}|_{t^k=0}$ in just the same way as we identify $V_d$ with an open subset of $U_{d+1/2}$ to form $\EuT_\infty$. The result is a proper scheme $\EuT^d_k \to \spec R[t]/(t^k)$ whose specialization to $t=0$ is a cycle of $d$ $\mathbb{P}^1$'s. We find it convenient to rename $t$ as $q^{1/d}$; so we have $\EuT^d_k \to \spec R[q^{1/d}]/(q^{k/d})$. These schemes form an inverse system in $k$; passing to the inverse limit as $k\to \infty$, we obtain a proper formal scheme $\hat{\EuT}^d \to \mathsf{Spf} \, R\series{q^{1/d}}$.

There are \'etale quotient maps $\hat{\EuT}^{d_1d_2} (q^{1/d_1}) \to \hat{\EuT}^{d_1}$, where by $\hat{\EuT}^{d_1d_2} (q^{1/d_1})$ we mean the formal base-change from $ \mathsf{Spf} \, R\series{q^{1/d_1d_2}}$ to  $\mathsf{Spf} \, R\series{q^{1/d_1}}$ given by $q^{1/d_1d_2}\mapsto q^{1/d_1}$. The construction of $\hat{\EuT}^d$ does not work in quite the same way when $d=1$, because one is then gluing $V_0$ to itself, and the Zariski-open cover by the $U_{i+1/2}$ becomes merely an \'etale cover. Nonetheless, we can define $\hat{\EuT}^1\to \mathsf{Spf}\, R\series{q}$ as the \'etale quotient $\hat{\EuT}^2/(\Z/2)$. This quotient is the locally-ringed space whose functions are the $\Z/2$-invariant functions of $\hat{\EuT}^2$. To see that it is a formal scheme, we need only note that $\hat{\EuT}^1$ can be covered by two $\Z/2$-invariant formal-affine open sets, which is straightforward to check.

Notice also that $\EuT_\infty\to \spec R[t]$ has sections $\sigma_{i+1/2}\colon \spec R[t]\to U_{i+1/2}$, defined by $X_i=1$ and $Y_{i+1}=t$. The $\Z$-action intertwines these sections, and consequently, $\hat{\EuT}^{1}\to \mathsf{Spf}\, R\series{q}$ has a distinguished section $\sigma$.

\subsubsection{Polarization} One now wants to polarize the formal scheme $\hat{\EuT}^1$, that is, to identify an ample line-bundle $\hat{\EuL} \to \hat{\EuT}^1$. A very ample power of $\hat{\EuL}$ will then define a projective embedding of $\hat{\EuT}^1$, and in doing so, will `algebraize' $\hat{\EuT}^1$, refining it to a true scheme over $R\series{q}$.\footnote{An analogous situation more familiar to geometers accustomed to working over $\C$ is that a projective embedding of a complex manifold makes this analytic object algebraic.}

That ample line bundles algebraize formal schemes is a general setting is a theorem of Grothendieck, but here, as noted by Gross, it is quite concrete. We will find a line bundle $\hat{\EuL}\to \hat{\EuT}^1$, by which we mean a sequence of line-bundles $\EuL_k\to \EuT_\infty|_{t^k=0}/\Z$ and isomorphisms $\EuL_k|_{t^{k-1}=0}\cong \EuL_{k-1}$. We will find a basis
$\{ \theta_{3, m/3} \}_{m=0,1,2}$ for $\h^0(\hat{\EuT}^1, \hat{\EuL}^{\otimes 3})$---by this we mean bases $\{\theta^{(k)}_{3, m/3} \}_{m=0,1,2}$ for $\h^0(\EuT_\infty|_{t^k=0}/\Z, (\EuL_k)^{\otimes 3})$, carried one to another by the maps in the inverse system---defining a plane embedding of $\hat{\EuT}^1$. In this way we will see that $\EuT_\infty|_{t^k=0}/\Z$ is cut out from $\mathbb{P}^2(R[t]/(t^k))$ by an equation which reduces modulo $t^{k-1}$ to that cutting out $\EuT_\infty|_{t^{k-1}=0}/\Z$. Passing to the limit, we see that $\hat{\EuT}^1$ is cut out by an equation from $\mathbb{P}^2(R\series{q})$, and hence arises from a projective scheme $\EuT$. 

\begin{figure}[ht!]\label{polytope}
\centering
\labellist
\small\hair 2pt 
\pinlabel $w_{-1/2}$ at 350 40
\pinlabel $w_{-3/2}$ at 420 100
\pinlabel $w_{-5/2}$ at 460 230
\pinlabel $w_{1/2}$ at 285 40
\pinlabel $w_{3/2}$ at 215 100
\pinlabel $w_{5/2}$ at 155 230

\endlabellist
\includegraphics[width=0.5\textwidth]{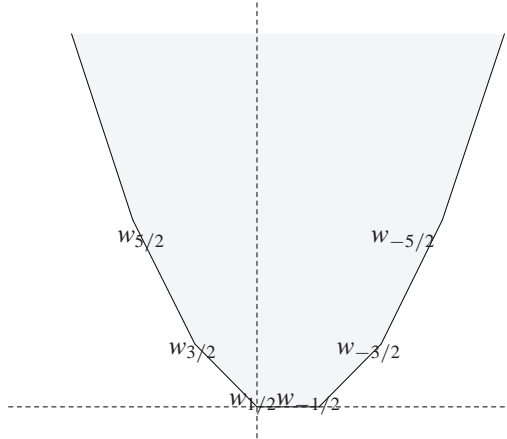}
\caption{The convex polygon $\Delta$}\end{figure}

As Gross explains, one obtains a line bundle $\EuL=\EuL_\Delta \to \EuT^{\infty}$ by observing that $F$ is the fan dual to an unbounded convex polygon $\Delta\subset \Q^2$, namely, the convex hull in $\Q^2$ of a sequence of points $w_{i+1/2} \in \Z^2$, $i\in \Z$, given by $w_{1/2}=(0,0)$ and $w_{i-1/2}-w_{i+1/2}=(1,-i)$ (Figure \ref{polytope}). 

The cone $\sigma_{i+1/2}^\vee$ is the tangent wedge $T_{i+1/2}$ at $w_{i+1/2}$. For each $j\in \Z + \half$, we define $\EuO_\Delta |_{U_j}$ to be the free $\EuO_{U_j}$-module $\EuO_{U_j}\cdot z^{w_j}$ on a generator $z^{w_j}$. Here $z$ is a formal symbol. We assemble the $\EuO_\Delta |_{U_j}$ into an invertible sheaf $\EuL=\EuO_\Delta \to \EuT_\infty$ by declaring that the transition function from $\EuL|_{U_{i+1/2}}$ to $\EuL|_{U_{i-1/2}}$ is multiplication by $z^{w_{1/2-i}-w_{1/2+i}}$. Each lattice point $\lambda \in \Delta\cap \Z^2$ defines a section $z^\lambda\in H^0(\EuL)$: over $U_{i+1/2}$, we view $\lambda$ as a point in $(w_{i+1/2}+\sigma_{i+1/2}^\vee)\cap \Z^2$, and so assign to it a section $w^\lambda|_{U_{i+1/2}} = X_i^a Y_{i+1}^b z^{i+1/2}$. These local sections agree on overlaps and so define a global section.  

The line bundle $\EuL$ is \emph{ample}: this is a general feature of line bundles associated with lattice polytopes. 

The $\Z$-action on $\EuT_\infty$ lifts to a $\Z$-action on $\EuL$ (we continue to follow Gross). Define $\tau \in SL_2(\Z)$ by
\[  \tau(x_1,x_2) = (x_1 + 1, x_1 + x_2).\] 
We have $\tau(w_{i+1/2})=w_{i-1/2}$, from which it is easy to see that $\tau(\Delta)=\Delta$. Thus $\tau$ generates an action of $\Z$ on $\Delta$. We lift the action of $m \in \Z$ to $\EuL$ as follows. Take the map $\phi_m \colon U_{j} \to U_{m+j}$ which defines the $\Z$-action (namely, $\phi_m(X_j)=X_{j+1}$ and $\phi_m(Y_{j+1})=Y_{j+m+1}$) and lift it to a map
\[ \widetilde{\phi}_m \colon \EuL\to \EuL, \quad z^\lambda \mapsto z^{\tau^m(\lambda)} ,\quad \lambda \in \Delta \cap \Z^2. \]
(we have specified in particular the effect of $\widetilde{\phi}_d$ on $z^{w_j}$).

As a result, $\EuL$ descends to a line bundle $\hat{\EuL}$ over each scheme $\EuT^d_k$, and indeed over the formal scheme $\hat{\EuT}^d$. The projective embeddings 
$\hat{\EuT}^d\to \mathbb{P} \h^0(\EuL^{\otimes N})^\vee$ for $N\gg 0$ cut out $\hat{\EuT}^d$ as a projective scheme $\EuT^d\to \spec \Z\series{q}^{1/d}$, lifting the formal scheme structure. This is even true for $d=1$. The scheme $\EuT= \EuT^1\to \spec \Z\series{q}$ is the Tate curve. It has its distinguished section $\sigma$.

The line bundle $\EuL\to \EuT$ has a global section $\theta := \sum_{k \in \Z}{ z^{w_k}}$. This formula is to be interpreted initially on $U_{i+1/2}|_{t^k=0}$, where it is a finite sum. It descends to a section of $\EuL$ over $\EuT_\infty|_{t^k=0}/(d\Z)$, and thereby a section over $\hat{\EuT}^d$ for each $d$. It is instructive to write $\theta$ in the open set $V_0$ in the following forms:
\[ \theta|_{V_0} = \sum_{k\in \Z}{q^{k(k-1)/2}  (z^{(1,0)})^k } = \sum_{k\in \Z}{(-1)^k q^{k(k-1)/2} \zeta^k} , \]
where $\zeta= -z^{(1,0)}$. Up to a factor of $iq^{1/4} \zeta$, $\theta(q,\zeta)$ is exactly the Fourier expansion for the classical theta-function $\vartheta_{1,1}$, written in terms of $q=e^{2\pi i\tau}$ and $\zeta = e^{2\pi i x}$, where $x$ is the coordinate on $\C/\langle 1, \tau \rangle$. 

Observe (i) that $\theta|_{V_0\cap \{ q=0\}}$ vanishes precisely where $\zeta = 1$, whence $\EuL$ has degree 1 on the geometric fibers of $\EuT\to \spec R\series{q}$; (ii) that $\theta(q, \zeta^{-1})= - \zeta^{-1} \theta(q,\zeta)$, so $\theta$ vanishes where $\zeta=1$; and hence (iii) that $\theta$ vanishes precisely where $\zeta=1$, i.e. along $\sigma$. Thus $\EuL\cong \EuO(\sigma)$. Using Riemann--Roch, we see that $\EuL^{\otimes 3}=\EuO(3\sigma)$ is very ample relative to the morphism $\EuT\to \spec R\series{q}$, and hence embeds $\EuT$ as a Weierstrass cubic in $\mathbb{P}^2(R\series{q})$. Thus $\EuT$ becomes a Weierstrass curve, with a canonical differential $\omega$. Moroever, this perspective makes clear that $\EuT$ coincides, as a Weierstrass curve, with the curve described in the introduction to this paper. (However, we shall not review the derivation of the Fourier expansions of $a_6(q)$ and $a_6(q)$).

\subsection{Lattice-points and multiplication of theta-functions}\label{lattice points and mult}
We have a canonical $\Z$-basis $\{z^\lambda : \lambda \in \Z^2\cap N\Delta \}$ for $\h^0(\EuT_\infty; \EuL^{\otimes N})$. The action of $q$ is given by $q\cdot z^\lambda = z^{\lambda+(0,1)}$. To get a basis over $\Z\series{q}$, it suffices to consider $\{ z^\lambda: \lambda \in \Z^2 \cap \partial(N \Delta)\}$, since we then obtain the remaining lattice points by multiplying by powers of $q$. There is also a canonical $\Z\series{q}$-basis $\beta_N$ for the $\Z$-invariant part $\h^0(\EuT_\infty; \EuL^{\otimes N})^\Z= \h^0(\EuT;  \EuL^{\otimes N})$,
\[ \beta_N = \{ \theta_{N,p} : p\in  C_N \},   \]
where $p$ runs over the cyclic group $C_N:=\left(\frac{1}{N} \Z\right) /\, \Z$. 
To define the `theta-functions' $\theta_{N,p}$, let $\phi\colon \Q \to \Q$ be the piecewise-linear function whose graph is the boundary of $\Delta$. Define automorphisms $\tau_N$ of the dilated polygon $N\Delta$, sending vertex $N w_j$ to the adjacent vertex $Nw_{j+1}$, by
\[ \tau_N (x_1,x_2) =  (x_1 + N, x_1 +x_2). \]
Then $(Np,N\phi(p))$ lies on the boundary of $N\Delta$, and hence so does $\tau^k_N(Np, N\phi(p))$ for each $k\in \Z$. Put
\[  \theta_{N,p} = \sum_{k \in \Z}{ z^{\tau^k_N( Np, N\phi(p))}}. \]
Again, the formula for $\theta_{N,p}$ makes sense, as a finite formal sum, on $U_{i+1/2}|_{t^k=0}$, hence on each $\hat{\EuT}^d$, and so finally on $\EuT$. It is clear that $\beta_N$ is a basis for $\h^0(\EuT, \EuL)$. One has $\theta_{1,0}=\theta$.

\begin{figure}\label{triangles}
\centering
\labellist 
\small\hair 2pt 
\pinlabel $T(p_1,p_2)$ at -50 300
\pinlabel $(p_1,0)$ at -20 430
\pinlabel $(E(p_1,p_2),0)$ at 160 430
\pinlabel $(p_2,-n_1(p_2-p_1))$ at 120 0
\pinlabel {\text{slope} $-n_1$} at 40 200
\pinlabel {\text{slope} $-(n_1+n_2)$} at 220 200
\endlabellist
\includegraphics[width=0.2\textwidth]{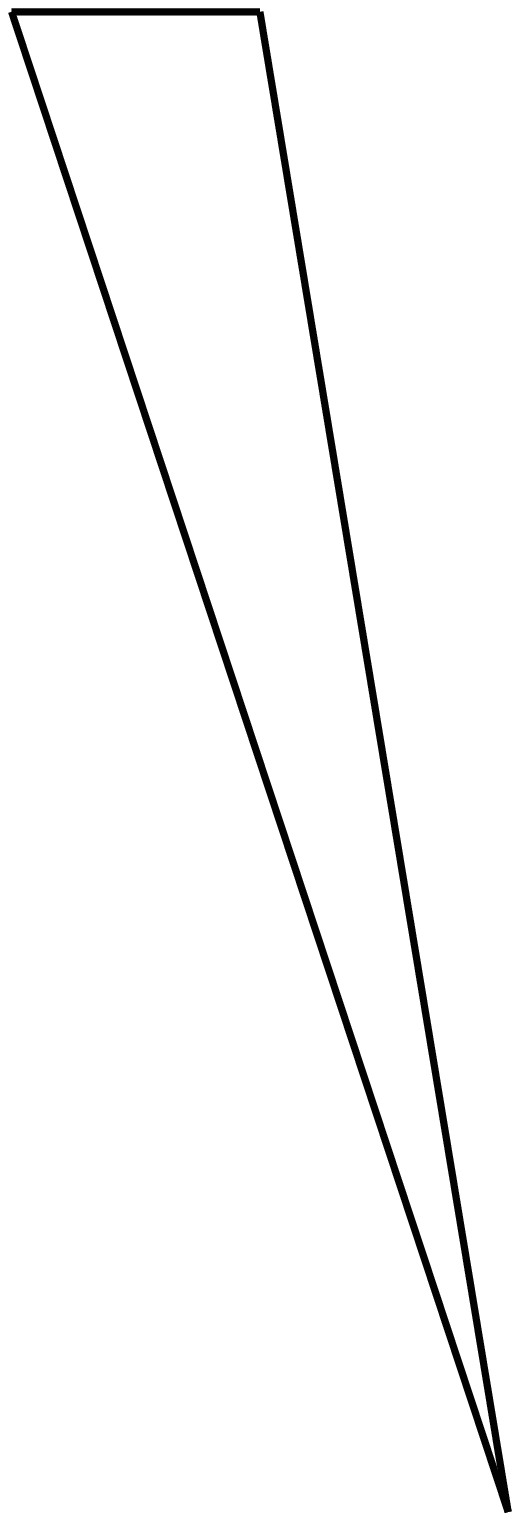}
\end{figure}

Explicit multiplication rules for classical theta-functions are standard. Whilst these theta-functions are not precisely identical to the classical theta functions which give canonical bases for $\h^0(\EuO(Np))$, they do obey a very similar multiplication rule \cite{Gro}:
\begin{equation} \label{structure constants}
\theta_{n_1, p_1} \cdot \theta_{n_2, p_2} = \sum_{j\in \Z} {q^{ \lambda(p_1,p_2+j)}  \theta_{n_1+n_2, E(p_1,p_2+ j ) }  } .
\end{equation}
Here $E(p_1,p_2)$ is the weighted average with respect to a distribution determined by $n_1$ and $n_2$,
\[  E(a,b) = \frac{n_1a + n_2 b}{n_1+n_2}, \] 
and
\begin{align} 
\lambda(p_1,p_2) 
& = n_1 \phi(p_1) + n_2\phi(p_2)  -   (n_1+n_2) \phi(E(p_1,p_2)) \\
& = (n_1+n_2) \left \{ E \left( \phi(p_1), \phi(p_2)\right) - \phi(E(p_1,p_2)) \right\}. \notag
\end{align}
Three points about $\phi$ are noteworthy:
\begin{itemize}
\item 
For any convex function $f$, the quantity $\Delta_f (p_1,p_2) : = E(f(p_1),f(p_2)) - f(E(p_1,p_2))$ is non-negative, by Jensen's inequality. One has $\lambda = (n_1+n_2) \Delta_\phi$. 
\item
$\phi$ is a piecewise-linear approximation to the quadratic function 
\begin{equation}\label{def psi}
\psi(x)=\half x(x-1);
\end{equation} 
indeed, $\phi(n)=\psi(n)$ for $n\in \Z$, and $\phi$ is affine-linear on intervals $[n,n+1]$. 
\item
The expression
\begin{align}
a(p_1,p_2) & := (n_1+n_2) \Delta_\phi(p_1,p_2)\\
& = n_1 \psi(p_1) + n_2\psi(p_2)  -   (n_1+n_2) \psi(E(p_1,p_2)) \notag \\
& = (n_1+n_2) \left \{ E \left( \psi(p_1), \psi(p_2)\right) - \psi(E(p_1,p_2)) \right\} ,\notag 
\end{align}
defined analogously to $\lambda$, is the area of the triangle $T(p_1,p_2)$ with vertices
\[  (p_1,0), \quad (p_2, -n_1(p_2-p_1)),\quad ( E(p_1,p_2) ,0)   \] 
(see Figure \ref{triangles}). 
\end{itemize}

The last fact is key to Gross's derivation of `classical' HMS \cite[section 8.4.2]{Gro}.\footnote{A similar multiplication rule, but for \emph{classical} theta-functions, plays an analogous role in Polishchuk--Zaslow's proof of cohomology-level mirror symmetry for elliptic curves \cite{PZ}. For those theta-functions, $a(p_1,p_2)$ is the relevant exponent.} For a fixed $n$, the theta-functions $\theta_{n,p}$ are \emph{almost} mirror to the intersection points of $L_0$ and $L_{(1,-n)}$. We say ``almost''  because correction factors are required, since it is $\lambda(p_1,p_2)$ and not $a(p_1,p_2)$ which appears in the product rule for theta-functions.

 It turns out that $\lambda$ has a similar geometric interpretation, one which will be equally central in our derivation of arithmetic HMS. We point out that whilst subsection \ref{Tate section} was a review, and the material of subsection \ref{lattice points and mult} has so far also been standard, this interpretation is to our knowledge original:

\begin{prop}\label{lattice points}
Fix $\epsilon>0$, and say a point in $\R^2$ is a \emph{perturbed lattice point} if it is congruent to $(\epsilon,\epsilon)$ modulo $\Z^2$. Then, if $\epsilon \ll (n_1+n_2)^{-1}$, the number of perturbed lattice points inside $T(p_1,p_2)$ is equal to $\lambda(p_1,p_2)$.
\end{prop}
Alas, we have not found an elegant proof of this proposition; the proof we give is an elementary but somewhat lengthy calculation:
\begin{pf}

Let $T_1$ be the right triangle with vertices $(p_1,0)$, $(p_2,0)$ and $(p_2  ,-n_1(p_2-p_1))$. Thus there is a horizontal edge, a vertical edge, and an edge which is a segment of the line $y=-n(x-p_1)$.  
Write $n_3=n_1+n_2$ and $p_3=E(p_1,p_2)$. Let $T_2$ be the right triangle with vertices $(p_3,0)$, $(p_2,0)$ and $(p_2  ,-n_3(p_2-p_3))$; thus $T_2$ is of the same form as $T_1$, substituting $n_3$ for $n_1$ and $p_3$ for $p_1$.  The number $\Lambda$ of perturbed lattice points in $T(p_1,p_2)$ is the difference
\[ \Lambda= \Lambda_1-\Lambda_2,\]
where $\Lambda_i$ is the number of perturbed lattice points in $T_i$ (see Figure \ref{triangles2} (left)).
\begin{figure}[ht!]\label{triangles2}
\centering
\labellist 
\small\hair 2pt 
\pinlabel $T_1$ at 0 200
\pinlabel $T_1$ at 300 200
\pinlabel $T_2$ at 132 280
\pinlabel $T_1'$ at 350 280

\endlabellist
\includegraphics[width=0.3\textwidth]{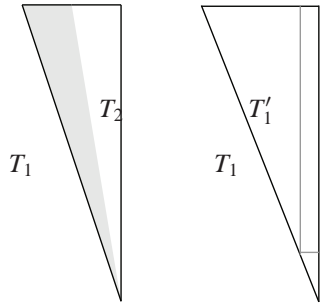}
\caption{Left:  $T_1$ divided into $T_2$ and $T(p_1,p_2)$. Right: shaving off a trapezium from $T_1$.}
\end{figure}

We shall compute $\Lambda_1$ by counting the points row by row, and shall then apply our formula to $T_2$ to obtain $\Lambda_2$.

Suppose, for $i=1$, $2$, we have $p_i = q_i + r_i/n_i$, where $q_i,r_i\in \Z$ and $0\leq r_i < n_i$. 
We have
\begin{align}  
\lambda(p_1,p_2) =  & \frac{n_1 }{2}q_1(q_1-1) + \frac{n_2}{2} p_2(p_2-1) - \frac{n_3}{2} q_3(q_3-1)+  r_1 q_1 + r_2q_2 -  r_3 q_3 
\end{align}
and we wish to show that $\Lambda=\lambda(p_1,p_2)$.

\paragraph{Case where $p_2\in \Z$.} In this case, $p_2=q_2$ and $r_2=0$. We calculate
\begin{align*}
\Lambda_1 & = \sum_{m=1}^{ \lfloor n_1(p_2-p_1)\rfloor}{(p_2 - \lceil p_1 + m/n_1 \rceil)}\\
&= \sum_{m=1}^{n_1(p_2- p_1)}{(p_2 - q_1 -  \lceil  (r_1+m)/n_1 \rceil)}\\
&= n_1(p_2-q_1)(p_2-p_1) - \sum_{m=1}^{n_1(p_2-p_1)} {\lceil (r_1+ m)/n_1 \rceil}\\
& =  n_1(p_2-q_1)(p_2-p_1)  - n_1 \left \{1+2+\dots +(p_2-q_1)\right \}  + r_1  \\
& =  n_1(p_2-q_1)(p_2-p_1)  - \frac{1}{2}n_1(p_2-q_1) (p_2-q_1+1)  + r_1  \\
& = \frac{1}{2}n_1(p_2-q_1) \left(p_2 - q_1 - \frac{2r_1}{n_1} -1 \right)  +  r_1  \\
& = \frac{1}{2}n_1q_1^2 + \frac{1}{2}n_1 p_2^2 + r_1 q_1  - n_1p_2q_1 -r_1p_2 -\frac{1}{2}n_1p_2 + \frac{1}{2}n_1q_1 + r_1.
\end{align*}
Hence, writing $p_3 =  q_3 + r_3/n_3$ with $0\leq r_3 \leq n_3$, we have
\begin{align*}
\Lambda_2  = \frac{1}{2}n_3q_3^2 + \frac{1}{2}n_3 p_2^2 + r_3 q_3  - n_3p_2q_3 -r_3p_2 -\frac{1}{2}n_3p_2 + \frac{1}{2}n_3q_3 + r_3.
\end{align*}
We now compute the difference $\Lambda$ by inputting the relation $r_1-r_3 = n_3q_3 - n_1q_1 -n_2 p_2$:
\begin{align*}
\Lambda = &  \Lambda_1 - \Lambda_2 \\
 =  & \frac{1}{2}n_1q_1^2   - \frac{1}{2}n_3q_3^2  + r_1 q_1 -r_3 q_3  + \frac{1}{2}n_1 p_2^2 - n_1p_2q_1 \\
 & -\frac{1}{2}n_1p_2 + \frac{1}{2}n_1q_1 +  \frac{1}{2}n_3p_2 - \frac{1}{2}n_3q_3 -  \frac{1}{2}n_3 p_2^2  \\
&+ n_3p_2q_3  +( n_3q_3 - n_1q_1 -n_2 p_2) - (n_3q_3 - n_1q_1 -n_2 p_2) p_2  \\
 =  & \frac{1}{2}n_1q_1^2 - \frac{1}{2}n_1q_1 
 + \frac{1}{2}n_2 p_2^2-\frac{1}{2}n_2p_2 
 - \frac{1}{2}n_3q_3^2  + \frac{1}{2}n_3q_3  + r_1 q_1   -r_3 q_3  \\
 = & \lambda(p_1,p_2).
 \end{align*}
\paragraph{General case.} We drop the assumption that $p_2$ is an integer. In this case, we can shave off the right-hand edge of the triangle $T_1$, so that its vertical edge is at $x=q_2$ instead of $x=p_2$ (Figure \ref{triangles2}). Our previous formula applies to this shaved triangle $T_1'$. The number of perturbed lattice points in the trapezium which we shaved off $T_1$ is computed as follows: the trapezium is a rectangular strip together with a triangle at the bottom. The triangle contains no lattice points, while the strip contains $n_1(q_2- p_1)$ lattice points.

Similarly, we shave off the right-hand edge of $T_2$; the trapezium which we shave off contains $n_3(q_2 - p_3)$ perturbed lattice points. Note that the difference between the counts of points in these trapezia is $n_1(q_2- p_1)-n_3(q_2 - p_3) = r_2$.

We deduce that the numbers of perturbed lattice points in $T_1$ and $T_2$ are, respectively,
\begin{align*}
\Lambda_1  & = \frac{1}{2}n_1q_1^2 + \frac{1}{2}n_1 q_2^2 + r_1 q_1  - n_1q_1q_2  -r_1q_2 -\frac{1}{2}n_1q_2 + \frac{1}{2}n_1q_1 + r_1 
 + n_1(q_2-p_1),\\
\Lambda_2  &=  \frac{1}{2}n_3q_3^2 + \frac{1}{2}n_3 q_2^2 + r_3 q_3  - n_3q_2q_3 -r_3q_2 -\frac{1}{2}n_3q_2 + \frac{1}{2}n_3q_3 + r_3  + n_3(q_2 - p_3) .
\end{align*}
The difference is
\begin{align*}
\Lambda =  & \Lambda_1-\Lambda_2\\
= &  \frac{1}{2}n_1q_1^2 + \frac{1}{2}n_1 q_2^2 + r_1 q_1  - n_1q_1q_2   -r_1q_2 -\frac{1}{2}n_1q_2 + \frac{1}{2}n_1q_1 + r_1\\
& - \frac{1}{2}n_3q_3^2 - \frac{1}{2}n_3 q_2^2 - r_3 q_3  + n_3q_2q_3 +  r_3q_2 +\frac{1}{2}n_3q_2 - \frac{1}{2}n_3q_3 - r_3 + r_2\\
=& \frac{1}{2}n_1q_1^2  - \frac{1}{2}n_3q_3^2 + \frac{1}{2}n_1 q_2^2- \frac{1}{2}n_3 q_2^2 \\
&  + r_1 q_1- r_3 q_3 - n_1q_1q_2 + n_3q_2q_3 -\frac{1}{2}n_1q_2+\frac{1}{2}n_3q_2+ \frac{1}{2}n_1q_1- \frac{1}{2}n_3q_3 + r_2\\
&+( n_3q_3 - n_1q_1 -n_2 p_2)  - q_2( n_3q_3 - n_1q_1 -n_2 p_2) \\
=& \frac{1}{2}n_1q_1^2  - \frac{1}{2}n_3q_3^2 - \frac{1}{2}n_2 q_2^2\\
&  + r_1 q_1- r_3 q_3 + \frac{1}{2}n_2q_2 - \frac{1}{2}n_1q_1+ \frac{1}{2}n_3q_3 + r_2  -n_2 p_2  + q_2 n_2 p_2 \\
=& \frac{1}{2}n_1q_1^2  - \frac{1}{2}n_3q_3^2 - \frac{1}{2}n_2 q_2^2 \\
&  + r_1 q_1- r_3 q_3 + \frac{1}{2}n_2q_2 - \frac{1}{2}n_1q_1+ \frac{1}{2}n_3q_3  - n_2 q_2 +  n_2q_2^2 +  r_2q_2 \\
=& \frac{1}{2}n_1q_1^2  - \frac{1}{2}n_3q_3^2 + \frac{1}{2}n_2 q_2^2  \\
&  + r_1 q_1 +    r_2 q_2 - r_3 q_3  - \frac{1}{2}n_1q_1 - \frac{1}{2}n_2q_2 + \frac{1}{2}n_3q_3   \\
& = \lambda(p_1,p_2).
\end{align*}
\end{pf}

\subsection{Homogeneous coordinate rings} \label{homog subsection}
We earlier discussed the affine coordinate ring $\varinjlim_{N} {\h^0(\EuO(Np))}$. To work with this ring, one must understand the direct system as well as the multiplication of sections of powers of $\EuO(p)$. On any Weierstrass curve $C\to \spec R$, the divisor $D=3p$ is very ample. It defines a homogeneous coordinate ring $R_C^0 = \bigoplus_{N\geq 0} {\h^0(\EuO(ND))}$, a graded $R$-algebra whose isomorphism determines the curve. The homogeneous coordinate ring has the advantage over the affine one that it is directly determined by the composition maps in $\perf C$. For a Weierstrass curve $(C,\sigma,\omega)$ over $\spec S$, let $R^*_C$ denote the extended homogeneous coordinate ring $\bigoplus_{n\geq 0} {\h^*(\EuO(nD))}$. This bigraded ring differs from $R_C^0$ only by the presence of the summand $\h^1(\EuO)$. This summand has its trace map $\tr_\omega \colon \h^1(\EuO)\to S$. The ring $R^*_C$ determines $C$, and the trace map then determines the differential $\omega$.  

The truncated ring $tR_C^*=R_C^*/I$, where $I=\bigoplus_{N>3}{\h^0(\EuO(ND))}$, with its bigrading and trace $\h^1(\EuO)\to S$, already determines $(C,\omega)$, because from it a defining cubic equation (in Hesse form) can be read off in $\h^0(\EuO(3D))$.  However, the regular section $\sigma$ might not be fully determined by $tR_C^*$, or even $R_C^*$, inasmuch as one could replace $\sigma$ by a different regular section $\sigma'$ such that the divisor $3\sigma'$ is linearly equivalent to $3\sigma$.  

\begin{pf}[Proof of Theorem \ref{mainth}, clauses (i), (ii)]
Let $C_{\mathsf{mirror}}$ be an abstract Weierstrass curve over $\Z\series{q}$ whose category $\EuB= \EuB_{\mathsf{mirror}}$ admits a quasi-isomorphism with $\EuA$ such that the induced isomorphism $H\EuA \to H\EuB$ is the standard one. We must show that $C_{\mathsf{mirror}}\cong \EuT$ as Weierstrass curves. We shall initially prove that they are isomorphic as curves with differential. For this it is sufficient to show that $tR_{C_{\mathsf{mirror}}}^*\cong tR^*_{\EuT}$ by a trace-respecting bigraded ring isomorphism.  
We recall that $L_{(1,-n)}^\#$ denotes an oriented exact Lagrangian in $T_0$ of slope $-n$, equipped with its non-trivial double covering, and graded in such a way that $HF^*(L_0^\#,L_{(1,-n)}^\#)=HF^0(L_0^\#,L_{(1,-n)}^\#)$ for $n\neq 0$ (this is not quite a complete specification of the grading). We have 
\[  R_{C_{\mathsf{mirror}}}^*\cong \bigoplus_{N \geq 0} {HF^*(L_0,L_{(1,-3N)})} \cong HF^1(L_0,L_0) \oplus \bigoplus_{N \geq 0} {HF^0(L_0,L_{(1,-3N)})}. \]

We use Prop. \ref{primitive} to set up the basepoint $z$ and the 1-form $\theta$ on $T_0$ in such a way that the exact curve $L_{(1,-n)}$ is the image in $\R^2/\Z^2$ of a straight line through the origin in $\R^2$---this for $0\leq n \leq 9$. We realize the non-trivial double covering of $L_{(1,-n)}$ by selecting the point $\star=\star_{-n}=(\epsilon, -n\epsilon)\in L_{(1,-n)}$, where $0<\epsilon\ll 1$, and declaring the double covering $\tilde{L}_{(1,-n)}\to L_{(1,-n)}$ to be trivial away from $\star$ and to exchange the sheets at $\star$.

Take $1\leq n\leq 9$. The differential on $CF^0(L_0^\#,L_{(1,-n)}^\#)$ is zero, since there are no immersed bigons bounding $L_0$ and $L_{(1,-n)}$. Thus we have a basis $B_n = \{x_{n,p} : p\in C_{n}\}$ for $HF^0(L_0^\#,L_{(1,-n)}^\#)$, where $C_n= \frac{1}{n}\Z / \Z$, and $x_{n,p}=[ p,0]\in \R^2/\Z^2 = T$. We have $L_{(1,-n)}=\tau^n(L_0)$ for $0\leq n\leq 9$, where $\tau$ is the nearly-linear Dehn twist along $L_\infty$ set up at (\ref{model twist}); hence $HF(L_{(1,-n_1)}^\#, L_{(1,-n_1-n_2)}^\#)$ has basis $\tau^{n_1}(B_{n_2})$ when $0<n_1<n_1+n_2\leq 9$. We claim that the Floer product $HF(L_{(1,-n_1)}^\#, L_{(1,-n_1-n_2)}^\#)\otimes_{\Z[[q]]}   HF(L_0^\#,L_{(1,-n_1)}^\#) \to HF(L_0,L_{(1,-n_1-n_2)})$ is given by
\[ \tau^{n_1}(x_{n_2,p_2})   \cdot  x_{n_1,p_1}  = \sum_{j\in \Z} { q^{ \lambda(p_1, p_2 + j) } x_{n_1+n_2, E(p_1, p_2 + j)} }, \] 
an expression formally identical to the multiplication rule for theta-functions (\ref{structure constants}); $E$ and $\lambda$ are as defined there. The contributions to $\tau^{n_1}(x_{n_2,p_2})\cdot x_{n_1,p_1}$ are immersed triangles, the images of embedded triangles in $\R^2$. The first vertex is a lift of $x_{n_1,p_1}$ to $\R^2$, which we assume is the point $ (p_1,0) \in \R^2$. The second vertex is a lift of $\tau^{n_1}(x_{n_2,p_2})=[p_2,n_1 p_2]$ which lies on the line of slope $-n_1$ through $(p_1,0)$; thus it is of form $(p_2+j, -n_1(p_2+j-p_1))$ where $j\in \Z$. The third vertex is then at $(0, E(p_1, p_2 + j))$. 

The contribution of the triangle just described is $\varepsilon q^{\lambda(p_1,  p_2 + j)}$, where the exponent is the number of perturbed lattice points in the triangle, as computed in Proposition \ref{lattice points}, and the sign $\varepsilon = \pm 1$ depends on the orientations, double coverings and the intersection signs of the corners determined by the orientations. A formula for $\varepsilon$ is given in  \cite{SeiGen2}; it is simplest when, as here, all the corners have intersection number $-1$ (and therefore even index for Floer \emph{co}homology). In that case, $\varepsilon = (-1)^s$, where where $s$ is the number of stars on the boundary. There are $\lceil E(p_1,p_2+j) \rceil - \lceil p_1 \rceil$ stars on the $L_0$ boundary, $\lceil p_2 + j \rceil -\lceil p_1\rceil$ on the $L_{(1,-n_1)}$ boundary, and 
$\lceil p_2 + j \rceil - \lceil E(p_1,p_2+j) \rceil $ on the $L_{(1,-n_1-n_2)}$ boundary, so $s=2 (j+ \lceil p_2 \rceil - \lceil p_1 \rceil)$ and 
$\varepsilon=+1$.  This justifies our claim.

Define a $\Z\series{q}$-linear map $t\psi\colon tR_{C_{\mathsf{mirror}}}^* \to tR_{\EuT}^*$ by linearly extending the assignment $t\psi(x_{3N,p}) = \theta_{3N,p}$ for $N\in \{1,2,3\}$ together with the canonical ring-isomorphism $t\psi\colon HF^*(L_0,L_0)\to H^*(\EuO)$ defined by the Weierstrass differential $\omega$. In view of (\ref{structure constants}) and Prop. \ref{lattice points}, and the fact that the unit of $HF^*(L_0,L_0)$ is also the unit of $tR_{C_{\mathsf{mirror}}}^*$, the map $t\psi$ is a bigraded ring isomorphism preserving the trace. Therefore, there is an isomorphism $\iota\colon (C_{\mathsf{mirror}}, \sigma_{\mathsf{mirror}}, \omega_{\mathsf{mirror} } ) \to (\EuT,\omega,\sigma')$, where $\omega$ is the standard Weierstrass differential on $\EuT$ and $\sigma'$ is some section of $\EuT^{\mathsf{sm}}\to \spec \Z\series{q}$, the regular locus in $\EuT$. 

Since $\EuT$ is a \emph{generalized elliptic curve} (see \cite{DR} or \cite[Def. 2.1.4]{Con}), one has a homomorphism $\Gamma \to \aut(\EuT/ \Z\series{q} )$, from the group $\Gamma$ of sections of $\EuT^{\mathsf{sm}}\to \Z\series{q}$, to the automorphism group of $\EuT$ fixing the base: sections act on $\EuT$ by fiberwise translations. Consequently, all regular sections are related by automorphisms of $\EuT$. Hence 
$(C_{\mathsf{mirror}}, \sigma_{\mathsf{mirror}}, \omega_{\mathsf{mirror} })$ is isomorphic to $(\EuT,\omega,\sigma)$. 

At this point we have proved clauses (i) and (ii) of Theorem \ref{mainth} apart from the uniqueness statement in (i). That is easily taken care of: $\EuA$ split-generates $\tw^\pi \EuF(T,z)$, and the functor $\psi$ is determined, up to natural quasi-equivalence, by its effect on a  full, split-generating subcategory.
\end{pf}
The proof just given also gives another identification $C_{\mathsf{mirror}}|_{q=0}$, which may be used as part of the argument for clauses (iii) and (iv) of Theorem \ref{mainth}, which we have already proved. 

\begin{rmk}
In the argument above, we identified the mirror essentially by describing a canonical isomorphism of graded rings
\[ \bigoplus_{n \geq 0} {HF^0(L_0^\#, L_{(1,-3n)}^\#)} 	\to 	S:= \bigoplus_{n\geq 0}	{\h^0(\EuO(3n\sigma))},\] 
given, in degrees $3n\geq 3M$ where $L_{(-1,3n)}^{\#})$ is exact, by matching up the canonical bases: $x_{3n,p} \mapsto \theta_{3n,p}$.  The argument did not depend on this assignment being the one arising from the functor $\psi$, but it is natural to expect that this is so. The assignment $\psi'(x_{3n,p})=\theta_{3n,p}$, defined when $n \leq M$, respects products. Hence, taking $M\geq 3$, we see that $\psi'$ extends to a ring isomorphism, because $S$ is generated in degree $3$ and subject only to a relation in degree 9.  The functor $\psi$ defines another such isomorphism. Thus $\psi' \circ \psi^{-1}$ extends to an automorphism of $S$. This must come from a Weierstrass automorphism $\alpha$ of $\EuT$. The only possibilities for $\alpha$ are the identity and the hyperelliptic involution, but we do not identify $\alpha$ here.
\end{rmk}

\begin{pf}[Proof of Theorem \ref{comparing functors}]

We begin by clarifying how to construct $\Phi$. We now consider $L_0^\#$ and $L_\infty^\#$ as objects of the `absolute' Fukaya category $\EuF(T)$ over $\Lambda_\C$---so their exactness plays no role. They form a full subcategory $\EuA_T$, which we identify with $\EuB_{E}$ for some mirror Weierstrasss curve $E\to \spec \Lambda_\C$. Thus we get a functor $\phi \colon \EuF(T)\to \tw \VB(E)$ inducing a quasi-equivalence $\dersplit \EuF(T_0) \to \perf(E)$. 

We identify $E$ via the homogeneous coordinate ring for its cubic embedding, just as we did for $C_{\mathsf{mirror}}$. Thus we must count triangles according to their areas, not the number of lattice points they contain. By the argument of \cite{PZ} (see also \cite{AS}), we identify the homogeneous coordinate ring of $E$ as that of $\EuT\times(\Lambda_\C)$, with the bases of intersection points $B_n$ mapping to standard bases of $\theta$-functions. Thus $E\cong \EuT(\Lambda_\C)=\EuT\times_{\Z\series{q}} \Lambda_\C$ as Weierstrass curves. Since $E$ is smooth, every bounded coherent complex is quasi-isomorphic to a perfect complex, and so $\perf(E)\simeq \der^b \coh(E)$. Thus we can think of $\Phi$ as a functor 
\[\Phi\colon \EuF(T) \to \widetilde{\der}^b \coh(\EuT(\Lambda_\C)).\]

We must now compare the two functors $\EuF(T,z)\to \coh (\EuT(\Lambda_\C))$: first $\Phi\circ e$, and second, $\psi$ followed a base-change functor---call this composite $\Psi$. Both extend to functors defined on $\tw^\pi \EuF(T,z)$, and it will suffice to show that these are homotopic. The inclusion $\EuA\to \EuF(T,z)$ induces a quasi-equivalence $\tw^\pi \EuA \to \tw^\pi \EuF(T,z)$, and the definitions of $\Phi$ and $\Psi$ both depend on a choice of quasi-inverse $ \tw^\pi \EuF(T,z)\to \tw^\pi \EuA$. We take this to be the same in both cases, so that the functors have the same effect on objects---otherwise the assertion of homotopy makes no sense. 

Take $\EuA$ to be the $A_\infty$-subcategory of $\EuF(T,z)$ with objects $L_0^\#$ and $L_{\infty}^\#$, and $\EuA_T$ the corresponding subcategory of $\EuF(T)$. 

The construction of the functor $e$ depends on choices of functions $f_L$ for each exact Lagrangian $L$ such that $\theta_L=df_L$. Pick such a function $f_{L_0}$, and choose $f_{L_\infty}$ so that the unique intersection point $x\in L_0\cap L_\infty$ has symplectic action 0. The endomorphism space $CF(L_0^\#,L_0^\#)$ is defined via a Hamiltonian image $\phi_H(L_0)$; endow this with the function $ f_{L_0}\circ \phi_H^{-1}$; similarly for $L_\infty$.  Since $e$ acts as the identity on objects, it defines a quasi-isomorphism $e\colon \EuA(\Lambda_\C) \to \EuA_T$. It induces the canonical isomorphism $H^*\EuA\to H^*\EuA_T$, because the relevant intersection points have action zero. 

We have quasi-isomorphisms $e^*\EuA_T \to \EuB_{\EuT(\Lambda_\C)}$ and $\EuA \to  \EuB_{\EuT(\Lambda_\C)}$, both inducing the canonical isomorphisms on cohomology. That is, we have two minimal $A_\infty$-structures on the graded algebra $A$, and identifications of both of them with $\EuB_{\EuT(\Lambda_\C)}$. By Theorem \ref{Ainf comparison}, these two identifications must be induced by an automorphism of $\EuT(\Lambda_\C)$ respecting the Weierstrass data $(\sigma,\omega)$. This can only be the identity map. The result follows.
\end{pf}

\end{document}